\documentclass[onecolumn,secnumarabic,nobibnotes,superscriptaddress,prx]{revtex4-2}

\usepackage{times}
\usepackage{color, xcolor, colortbl}
\usepackage{graphicx}
\usepackage{epstopdf}
\ifpdf
\DeclareGraphicsExtensions{.eps,.pdf,.png,.jpg}
\else
\DeclareGraphicsExtensions{.eps}
\fi

\usepackage{geometry}
\usepackage{amsmath,amssymb, amsthm}
\usepackage{bm}
\usepackage[caption=false]{subfig}
\usepackage[title]{appendix}
\usepackage{multirow}
\usepackage{braket}
\usepackage[english]{babel}
\usepackage{hyperref}
\usepackage[capitalize]{cleveref}
\usepackage{natbib}
\setcitestyle{square,numbers}
\usepackage{adjustbox}
\usepackage{xspace}
\usepackage{comment}
\usepackage{booktabs}
\usepackage{makecell}
\usepackage{threeparttable}

\usepackage{pifont}
\newcommand{\cmark}{{\color{green}\ding{51}}}%
\newcommand{\xmark}{{\color{red}\ding{55}}}%

\newcommand{\bvec}[1]{\mathbf{#1}}

\newcommand{\valpha}{ {\bm{\alpha}} }
\newcommand{\vbeta}{ {\bm{\beta}} }

\newcommand{\vk}{\bvec{k}}

\newcommand{\vq}{\bvec{q}}
\newcommand{\vrr}{\bvec{r}}     

\newcommand{\vx}{\bvec{x}}
\newcommand{\vy}{\bvec{y}}
\newcommand{\vz}{\bvec{z}}

\newcommand{\vG}{\bvec{G}}

\newcommand{\vR}{\bvec{R}}

\newcommand{\I}{\mathrm{i}}

\newcommand{\abs}[1]{\left\lvert#1\right\rvert}

\newcommand{\ud}{\,\mathrm{d}}
\newcommand{\Or}{\mathcal{O}}

\newcommand{\RR}{\mathbb{R}}

\newtheorem{thm}{\protect\theoremname}
\theoremstyle{plain}
\newtheorem{lem}[thm]{\protect\lemmaname}
\theoremstyle{plain}
\newtheorem{rem}[thm]{\protect\remarkname}
\theoremstyle{plain}

\theoremstyle{plain}
\newtheorem{cor}[thm]{\protect\corollaryname}
\newtheorem{assumption}[thm]{\protect\assumptionname}
\newtheorem{defn}[thm]{\protect\definitionname}

\providecommand{\definitionname}{Definition}
\providecommand{\assumptionname}{Assumption}
\providecommand{\corollaryname}{Corollary}
\providecommand{\lemmaname}{Lemma}
\providecommand{\propositionname}{Proposition}
\providecommand{\remarkname}{Remark}
\providecommand{\theoremname}{Theorem}
\providecommand{\problemname}{Problem}

\numberwithin{figure}{section}
\numberwithin{table}{section}

\newcommand{\xsum}{\mathop{\sum\nolimits'}}

\crefname{equation}{Eq.}{Eqs.}
\crefformat{equation}{Eq.~#2(#1)#3}

\newcommand{\LL}[1]{\textcolor{brown}{[LL:#1]}}

\newcommand{\REV}[1]{\textcolor{black}{#1}}

\newcommand{\DeptMath}{Department of Mathematics, University of California, Berkeley, CA 94720, USA}
\newcommand{\LBLMath}{Applied Mathematics and Computational Research Division, Lawrence Berkeley National Laboratory, Berkeley, CA 94720, USA}

\begin{document}

\title{Inverse Volume Scaling of Finite-Size Error in Periodic Coupled Cluster Theory}

\author{Xin Xing}
\affiliation{\DeptMath}
\author{Lin Lin}
\affiliation{\DeptMath}
\affiliation{\LBLMath}
\date{\today}

\begin{abstract}
Coupled cluster theory is one of the most popular post-Hartree--Fock  methods for \textit{ab initio} molecular quantum chemistry.
The finite-size error of the correlation energy in periodic coupled cluster calculations for three-dimensional insulating 
systems has been observed to satisfy the inverse volume scaling, even in the absence of any correction schemes. 
This is surprising, as simpler theories that utilize only a subset of the coupled cluster diagrams exhibit much slower decay 
of the finite-size error, which scales inversely with the length of the system.
In this study, we review the current understanding of finite-size error in quantum chemistry methods for periodic systems. 
We introduce new tools that elucidate the mechanisms behind this phenomenon in the context of coupled cluster doubles calculations.  This reconciles some seemingly paradoxical statements related to finite-size scaling. 
Our findings also show that singularity subtraction can be a powerful method to effectively reduce finite-size errors in practical quantum chemistry calculations for periodic systems.

\end{abstract}

\maketitle


\section{Introduction}

In the past few decades, \textit{ab initio} methods for quantum many-body systems, such as density functional theory (DFT), quantum Monte Carlo methods, and quantum chemistry wave function methods, are becoming increasingly accurate and applied to ever larger range of systems~\cite{Martin2008,MartinReiningCeperley2016}. Unlike molecular systems,
periodic systems, including solids and surfaces, require calculating properties in the thermodynamic limit (TDL), a theoretical state in which the system size approaches infinity. However, the TDL cannot be directly accessed in practical applications. Finite-sized computational supercells are employed to approximate this limit, which introduces finite-size errors into the calculations. Finite-size errors can significantly affect the accuracy of calculations, even for systems with thousands of atoms. An extreme case is moir\'e system such as  magic angle twisted bilayer graphene (MATBG), where each computational unit cell consists of approximately 10,000 atoms, and the supercell needs to have more than 100,000 atoms to capture subtle correlation effects~\cite{BistritzerMacDonald2010,BultinckKhalafLiuEtAl2020}.
Directly tackling finite-size effects by enlarging the supercell size is very demanding, even for relatively inexpensive DFT calculations with modern-day supercomputers. For more accurate theories, this task is often computationally intractable. Understanding the finite-size scaling, i.e., the scaling of the finite-size error with respect to the system size, and employing finite-size error correction schemes are, therefore, crucial for obtaining accurate results using moderate-sized calculations.

The sources of finite-size errors in \textit{ab initio} calculations are multifaceted and complex~\cite{MakovPayne1995,FraserFoulkesRajagopalEtAl1996,GruberLiaoTsatsoulisEtAl2018}. These errors are influenced by numerous factors, including system characteristics such as whether it is insulating or metallic, or whether it is a three-dimensional bulk system versus a low-dimensional system. 
Calculations of electron kinetic energy, electron-ion interaction energy, Hartree energy, Fock exchange energy, and electron correlation energy can all contribute to finite-size errors.
The first four types are predominantly single particle in nature, while the electron correlation energy is significantly more complex. To a large extent, electronic correlation is short-ranged, and this characteristic has spurred the development of local correlation methods, whose computational cost may scale linearly with system size. However, for \textit{ab initio} methods to be accurate, they must also effectively account for van der Waals (vdW) interactions~\cite{HermannDiStasioJrTkatchenko2017}. In solids, the cumulative effect of weak van der Waals interactions can become a significant contributor to the energy. The convergence of vdW energy follows an inverse volume scaling, implying that the finite-size error is inversely proportional to the volume of the supercell. 
The origin of finite-size error is not solely confined to the characteristics of a physical system. The finite-size scaling can vary significantly across different theoretical frameworks even when applied to the same physical system. Sometimes, the finite-size scaling can manifest differently even when using the same theoretical framework for the same system, simply by altering the shape of the computational unit cell.
All these complexities require careful analysis under varying scenarios, and this is important for accurately extrapolating towards the TDL or, preferably, for devising improved correction schemes aimed at reducing finite-size errors.

In recent years, there has been a growing interest in employing quantum chemistry wave function methods, such as M{\o}ller-Plesset perturbation theory and coupled cluster (CC) theory~\cite{ShavittBartlett2009}, 
to compute ground-state and excited-state properties for periodic systems~\cite{MarsmanGruneisPaierKresse2009,GruneisMarsmanKresse2010,GruberGruneis2018, Gruneis2015,MuellerPaulus2012, BoothGruneisKresseAlavi2013, SchaeferRambergerKresse2017,McClainSunChanEtAl2017,GruberLiaoTsatsoulisEtAl2018, ZhangGruneis2019, GaoSunJasonMottaEtc2020, WangBerkelbach2020, ZhuChan2021, NeufeldYeBerkelbach2022, BintrimBerkelbach2022, BanerjeeSokolov2022, NeufeldBerkelbach2023}.
Originally developed to study nuclear physics~\cite{CoesterKuemmel1960,Kuemmel1991}, CC theory has become one of the most popular methods in quantum chemistry
that involves electronic correlation~\cite{Cizek1966}. The coupled cluster singles, doubles and perturbative triples (CCSD(T)) theory is often referred to as the ``gold standard'' in molecular quantum chemistry.
In these methods, there are two primary strategies to approximate the TDL. The first involves expanding the computational supercells within real space.
The second strategy involves performing calculations using a fixed unit cell and refining discretization of the Brillouin zone within reciprocal space using a $\vk$-point mesh, such as the Monkhorst-Pack mesh~\cite{MonkhorstPack1976}. 
This paper focuses on the latter approach, where the number of $\vk$ points is denoted by $N_{\vk}$. If the Monkhorst-Pack mesh  includes the Gamma point of the Brillouin zone, this approach is equivalent to using a supercell comprised of $N_{\vk}$ unit cells. 
The convergence toward the TDL can be studied by increasing this single parameter $N_{\vk}$ toward infinity. The computational expense of quantum chemistry methods can rise sharply with respect to the system size. 
For example, in coupled cluster singles and doubles (CCSD) calculations, the computational cost scaling is $\Or(N_{\vk}^6)$ in real-space implementations~\cite{ShavittBartlett2009} and $\Or(N_{\vk}^4)$ in reciprocal-space implementations~\cite{HirataPodeszwaTobitaEtAl2004}.
Therefore, even moderate-sized calculations with $N_{\vk}=2\times2\times 2$ or $3\times 3\times 3$ can already be computationally challenging.  
In CCSD(T) calculations, the computational cost of the real-space implementation can escalate to $\Or(N_{\vk}^7)$ in the worst-case scenario. This implies that even a moderate refinement of the Monkhorst-Pack grid by a factor of 2 along each dimension can lead to an increase in computational cost by $(2^3)^7\approx 2$  million folds. 
Given this computational challenge imposed by system size, it is mostly impractical to estimate the finite-size scaling using power-law fitting over calculations on increasingly large systems 
and then determine the TDL value by extrapolation. 
Instead, a more feasible approach is to acquire the exact finite-size scaling through rigorous mathematical analysis and subsequently utilize power-law extrapolation to estimate the TDL value.

Although numerous empirical studies have examined the finite-size scaling of some quantum chemistry methods, there has been a notable lack of rigorous analysis. 
To the best of our knowledge, the first rigorous analysis of finite-size error in  Hartree-Fock (HF) theory and second-order M{\o}ller-Plesset perturbation theory (MP2) for insulating systems has been conducted only recently~\cite{XingLiLin2023}. 
The principal findings in Ref.~\cite{XingLiLin2023} include the following.

\begin{enumerate}

\item In the absence of any corrections, the finite-size error in HF scales as $\Or(N_{\vk}^{-\frac13})$. 
As $N_{\vk}$ is proportional to the supercell volume, $N_{\vk}^{\frac13}$ is proportional to the length of the supercell (we always assume uniform refinement of the $\vk$-point mesh along all three dimensions for any shaped unit cell). 
This scaling is, thus, also referred to as inverse length scaling. 

\item  By applying the Madelung constant correction, the finite-size error in HF improves to $\Or(N_{\vk}^{-1})$, i.e., the inverse volume scaling.

\item The finite-size error of the MP2 correlation energy also satisfies the inverse volume scaling $\Or(N_{\vk}^{-1})$. 
\end{enumerate}

As can be seen from HF and MP2, the specific finite-size scaling depends on the level of computational theory and numerical treatment employed in the calculations. 
In addition to confirming the empirically observed scalings, this rigorous analysis offers significant additional insights. 
It elucidates the nature of the singularity of the electron repulsion integral (ERI) due to the Coulomb kernel and its impact on finite-size scaling; it explains why the finite-size error can depend on the shapes and symmetries of the unit cell; and it has led to the development of new correction schemes, such as the staggered mesh method, which can expedite the convergence of HF, MP2, and random phase approximation (RPA) calculations toward TDL~\cite{XingLiLin2021,XingLin2022}.

An innovative theoretical observation provided in Ref.~\cite{XingLiLin2023} is that the finite-size error can be largely comprehended if assuming HF single-particle orbitals can be acquired exactly at any given $\vk$ point in the Monkhorst-Pack mesh. 
This perspective helps to disentangle the contribution due to the relatively manageable single-particle effects from the collective and more complex electron correlation effects.
Based on this assumption, the finite-size error can be rigorously examined from a numerical quadrature perspective. 
Specifically, the value of a physical observable at the TDL can often be written as a multidimensional integral over the Brillouin zone. 
As the Monkhorst-Pack mesh forms a uniform grid discretizing the Brillouin zone, the analysis simplifies to investigating the quadrature error of certain trapezoidal rule in a periodic region, a topic widely discussed in numerical analysis literature. 
The novelty here lies in the recognition that the associated integrands  possess a unique singularity structure that is asymptotically of a specific fractional form.
Reference~\cite{XingLiLin2023}, therefore, develops a new Euler-MacLaurin type of analysis that facilitates the study of the finite-size error associated with HF and MP2 methods, taking into account this fractional form singularity.

\REV{Compared to HF and MP2, the  finite-size error analysis and finite-size error correction methods in CC methods remain nascent. 
Early works in \cite{LiaoGrueneis2016,GruberGruneis2018} focus on an important intermediate quantity called the structure factor, and develop finite-size correction methods based on the different types of polynomial interpolation of the structure factor near the Coulomb singularity.
Subsequent studies on the uniform electron gas system \cite{MihmYangShepherd2020, MihmWeilerShepherd2023} numerically examine the power-law scaling of the finite-size error in CC theory, offering new analytical perspective on the structure factor. These works also provide
new error correction methods, such as twist-angle techniques \cite{MihmMcIsaacShepherd2019,MihmBenschotenShepherd2021} and structure factor interpolation \cite{WeilerMihmShepherd2022} for these metallic systems. 
Very recently, Ref.~\cite{SchaferBenschotenShepherdGruneis2023} proposes to replace the Coulomb kernel in CC calculations with an averaged Coulomb kernel in each quadrature element to reduce the finite-size error for anisotropic systems.
Despite these developments, theoretical understanding of the finite-size error scaling in CC methods (with and without error correction methods) has been lacking. 
}

Recently, by applying a similar approach, we have expanded the finite-size error analysis to include the correlation energy of the third-order M{\o}ller-Plesset perturbation theory (MP3) and CC theory~\cite{XingLin2023}. The simplest CC theory is the coupled cluster singles (CCS) theory. However, due to Thouless's theorem~\cite{Thouless1960}, CCS only rotates the Slater determinant to another Slater determinant. 
For a non-degenerate closed-shell Hartree-Fock reference, the correlation energy from the CCS theory vanishes. Therefore, we focus on the coupled cluster doubles (CCD) theory, which is mathematically the simplest and representative form of CC theory. 
From a diagrammatic perspective, CC diagrams encompass all M{\o}ller-Plesset perturbation diagrams.
When the CCD amplitude equation is solved iteratively with $n$ fixed point iterations (referred to as the CCD$(n)$ scheme), the MP2 and MP3 diagrams can be generated from CCD$(1)$ and CCD$(2)$, respectively. 
More generally, CCD($n$) consists of a finite subset of M{\o}ller-Plesset perturbation diagrams.
It is worth pointing out that there is a method called CCSD(2) which uses second-order perturbation to rectify CCSD energies for multireference and open-shell systems~\cite{GwaltneyHeadGordon2001}, and its meaning should, therefore, be distinguished from our CCD($n$).
Our analysis uncovers that the finite-size error in MP3 and CC is fundamentally different from that in MP2. 
Specifically, theories such as MP2 and RPA incorporate only ``particle-hole'' types of diagrams. 
The integrands corresponding to these diagrams are singular, but the singularity is relatively mild. 
However, starting from MP3, additional perturbative terms, namely the ``particle-particle'' and ``hole-hole'' diagrams, must be considered. 
These terms introduce much stronger singularities, necessitating the development of new analytical tools.
Our quadrature error analysis adapts the Poisson summation formula in a new setting and aligns with a recently developed trapezoidal quadrature analysis for certain singular integrals~\cite{Izzo2022}. 
This approach provides more accurate estimates and can be more applicable than our previous quadrature analysis based on the Euler-MacLaurin formula for HF and MP2 in~\cite{XingLiLin2023}.
In the absence of finite-size correction schemes and assuming exact orbital energies at any $\vk$ point, the study in Ref.~\cite{XingLin2023} concludes that the finite-size errors of MP3 and CCD both satisfy the inverse length scaling $\Or(N_{\vk}^{-\frac13})$.

Interestingly, for CCD, earlier numerical calculations did not provide conclusive evidence regarding its finite-size scaling, with different studies suggesting either an inverse volume scaling \cite{LiaoGrueneis2016,BoothGruneisKresseAlavi2013} or  inverse length scaling \cite{McClainSunChanEtAl2017}.
More recent calculations demonstrate that the electron correlation energy in periodic coupled cluster calculations should follow an inverse volume scaling, even in the absence of finite-size correction schemes. 
This observation points to a significant gap in the theoretical understanding of the finite-size error, and 
prompts the \textit{central question} we address in this paper:\textit{
How can we reconcile the following seemingly paradoxical facts? 
}

\begin{enumerate}

\item Without finite-size corrections, the finite-size error in CCD exhibits inverse volume scaling.

\item Without finite-size corrections, the finite-size error in MP3 exhibits inverse length scaling. This rate is sharp and cannot be improved. 

\item All MP3 diagrams are encompassed within the CCD formulation.

\end{enumerate}

There are several often-cited physical justifications for expecting that the CCD method, and CC theory more generally, may exhibit superior behavior when applied to periodic systems. 
One such reason stems from the size extensivity of CC theory. 
A theory is \textit{size extensive} when the total energy of two noninteracting identical systems, calculated as a combined system, equates to twice the energy of one system computed independently. 
Unlike theories such as truncated configuration interaction methods, which are not size extensive, truncated CC theory (such as CCD) possesses this advantageous characteristic. 
Another reason is that CC theory can be formulated in such a way that it does not explicitly depend on orbital energies. One practical consequence of this is that upon the convergence of the coupled cluster iterations, the Madelung constant correction, often used to reduce  finite-size effects in many-body simulations, cancels out naturally. Therefore in this scenario, CC without the Madelung constant correction is equivalent to that with the Madelung constant correction.

We would like to clarify that neither size extensivity nor the cancellation of the Madelung constant alone is sufficient to  address the aforementioned question. While size extensivity is indeed a desirable property, many methods such as HF, MP2, and MP3, among others, are all size extensive. In fact, given that periodic systems are infinitely large, size extensivity should be viewed as a \textit{necessary condition} for the applicability of any numerical method to such systems. However, this property does not provide insight into the convergence rate of the finite-size error.
The cancellation of the Madelung constant plays a pivotal role here. Nevertheless, demonstrating that CCD calculation  (or even MP3) with the Madelung constant correction exhibits inverse volume scaling in finite-size error is itself a significant challenge. This requires the development of new technical tools not currently available in existing literature. Indeed, the development and application of these tools are the main technical contribution of this paper.

In this paper, we elucidate the origin of the inverse volume scaling behavior. Our analysis consists of two steps. 
First, we investigate the structure of the CCD amplitude equation. 
We show that the Madelung constant correction, commonly used to reduce finite-size errors in Fock exchange energy and orbital energies, can also be applied to reduce the finite-size error in ERI contractions within the CCD amplitude equation. 
We establish a connection between the Madelung constant correction and a quadrature error reduction technique known as the singularity subtraction method~\cite{DavisRabinowitz2007}. 
By subtracting the leading singular terms from the integrands in the numerical quadratures, the Madelung constant correction reduces the finite-size errors in both the ERI contractions and the orbital energies from $\Or(N_{\vk}^{-\frac13})$ to $\Or(N_{\vk}^{-1})$. 
Furthermore, we demonstrate that with the Madelung constant correction, the finite-size errors in CCD($n$) and converged CCD calculations satisfy the desired inverse volume scaling.

In the second step of our analysis, we observe that upon convergence of the CCD amplitude equations, the Madelung constant corrections to both orbital energies and ERI contractions 
perfectly cancel each other out for any finite-sized system. 
This cancellation ensures that the CCD correlation energy remains the same, regardless of whether the Madelung constant correction is applied. 
Combining this result with the first step, we conclude that the finite-size error of the correlation energy in converged CCD calculations satisfies the desired inverse volume scaling without the need for any additional correction schemes. 
However, prior to the convergence of the amplitude equations, this perfect cancellation does not occur, and the finite-size error of CCD($n$) calculations remains $\Or(N_{\vk}^{-\frac13})$. 
A similar lack of cancellation occurs when the orbital energies take their exact value at the TDL but the ERI contractions are not corrected, resulting in an $\Or(N_{\vk}^{-\frac13})$ finite-size error for both converged CCD and CCD($n$) calculations studied in Ref.~\cite{XingLin2023}.

To validate our theoretical analysis, we perform CCD calculations on a 3D periodic hydrogen dimer system using the PySCF software package~\cite{SunBerkelbachBluntEtAl2018}. 
Our numerical results support the conclusions drawn from the theoretical analysis and provide further evidence for the finite-size scaling behavior summarized in Table~\ref{tab:summary}. 

\begin{table}[htbp]
        \centering
        \caption{
        Finite-size scaling of different computational theories with and without corrections to orbital energies and ERI contractions.
        The first two rows refer to the finite-size scaling of the Hartree-Fock (HF) exchange energy, and the remaining rows refer to the finite-size scaling of the correlation energy (excluding the HF exchange).
        While we focus on the Madelung constant correction in Ref.~\cite{XingLiLin2023} and this work, any correction schemes that can reduce the 
        finite-size errors in orbital energies and ERI contractions to $\Or(N_\vk^{-1})$ can also be applied, and the conclusions drawn
        here remain valid. In particular, exact values of single particle orbital energies at the TDL satisfy the condition
        above. N/A means that Madelung constant correction does not apply within the theory.
        }
        \label{tab:summary}

        \def\arraystretch{1.5}
        \begin{tabular}{c@{\hskip 1em}c@{\hskip 1em}c@{\hskip 1em}c@{\hskip 1em}c}
        \toprule
        Theory &\makecell{Correction to\\
        orbital energies}
        & 
        \makecell{Correction to\\
        ERI contractions} & \makecell{Finite-size\\
        scaling}    & Reference             \\
        \midrule
        HF  & N/A & \xmark & $N_\vk^{-\frac13}$ & \cite[Thm 3.1]{XingLiLin2023}, \cite{McClainSunChanEtAl2017, FraserFoulkesRajagopalEtAl1996, ChiesaCeperleyMartinEtAl2006}\\
        HF  & N/A & \cmark & $N_\vk^{-1}$ &  \cite[Thm 5.1]{XingLiLin2023} \cite{ChiesaCeperleyMartinEtAl2006}     \\
        MP2 & \cmark      & N/A & $N_\vk^{-1}$ &  \cite[Thm 4.1]{XingLiLin2023} \cite{LiaoGrueneis2016, BoothGruneisKresseAlavi2013}    \\
        \makecell{RPA, SOSEX,\\drCCD}  & \cmark & N/A & $N_\vk^{-1}$ &   \cite{XingLin2022,LiaoGrueneis2016}    \\        
        MP3 & \cmark & \xmark & $N_\vk^{-\frac13}$ & \cite{XingLin2023} \\
        MP3  & \cmark & \cmark & $N_\vk^{-1}$ &    Thm \ref{thm:error_ccdn}   \\        
        CCD$(n)$/CCD & \cmark    & \xmark & $N_\vk^{-\frac13}$ / $N_\vk^{-\frac13}$ & \cite{XingLin2023} \\
          CCD$(n)$/CCD & \cmark    & \cmark & $N_\vk^{-1}$ / $N_\vk^{-1}$ & Thm \ref{thm:error_ccdn} / Thm \ref{thm:error_ccd_madelung}  \\ 
        CCD$(n)$/CCD & \xmark    & \cmark & $N_\vk^{-\frac13}$ / $N_\vk^{-\frac13}$ & Thm \ref{thm:error_ccdn} / Thm \ref{thm:error_ccd_madelung}  \\ 
        CCD$(n)$/CCD & \xmark    & \xmark & $N_\vk^{-\frac13}$ / $N_\vk^{-1}$ & Thm \ref{thm:error_ccdn} / Cor \ref{thm:error_ccd_inexact}\\
        \bottomrule
        \end{tabular}
\end{table}

The paper is organized as follows. 
Section 2 introduces background knowledge and basic notations. 
Section 3 first decomposes the finite-size error in CCD$(n)$ calculations into the errors in four basic components and then describes how 
the four components contribute to the overall finite-size error with the possible Madelung constant correction. 
Section 4 explains the key ideas of how the Madelung constant correction can reduce the finite-size error in orbital energies and ERI conttractions 
from a numerical quadrature perspective. 
Section 5 illustrates the numerical results that corroborate our error estimate. 
Lastly, Section 6 discusses the implication of our theoretical study and future directions. 

\section{Background}
Consider a unit cell and its Brillouin zone as $\Omega$ and $\Omega^*$, respectively.
Denote the associated real- and reciprocal-space lattices as $\mathbb{L}$ and $\mathbb{L}^*$. 
To model such a periodic system, the Brillouin zone $\Omega^*$ is commonly discretized using a uniform mesh $\mathcal{K}$ of size $N_\vk$ (known as the Monkhorst-Pack grid~\cite{MonkhorstPack1976}). The orbitals and orbital energies (also called bands and band energies) 
$\{\psi_{n\vk}, \varepsilon_{n\vk}\}$ indexed by orbital indices $n$ and momentum vectors $\vk\in\mathcal{K}$ can be solved by
the HF method.
As a common practice, $n\in\{i,j\}$ refers to an occupied orbital and $n\in\{a,b\}$ refers to an virtual orbital. 
Throughout this paper, we use the \textit{normalized} ERI:
\begin{equation}\label{eqn:eri}
        \braket{n_1\vk_1,n_2\vk_2|n_3\vk_3,n_4\vk_4}
        = 
        \frac{4\pi}{\abs{\Omega}} \xsum_{\vG\in\mathbb{L}^*}
        \frac{1}{\abs{\vq+\vG}^2}  
        \hat{\varrho}_{n_1\vk_1,n_3(\vk_1 + \vq)}(\mathbf{G}) \hat{\varrho}_{n_2\vk_2,n_4(\vk_2 - \vq)}(-\mathbf{G}),
\end{equation}
where $\vq = \vk_3 - \vk_1$ is the momentum transfer vector, the crystal momentum conservation $\vk_1 + \vk_2 -\vk_3 - \vk_4 \in \mathbb{L}^*$ is assumed implicitly, 
and $\hat{\varrho}_{n'\vk', n\vk}(\vG) = \braket{\psi_{n'\vk'} | e^{\I(\vk' - \vk - \vG)\cdot \vrr} | \psi_{n\vk}}$ is Fourier representation of the pair product.
The primed summation over $\vG$ means that the possible term with $\vq + \vG = \bm{0}$ is excluded in the numerical calculation. 
Using a finite mesh $\mathcal{K}$, the HF orbital energy without any correction is computed as 
\begin{equation}\label{eqn:orbital_energy}
        \varepsilon_{n\vk}^{N_\vk} = 
        \braket{n\vk | \hat{H}_0| n\vk} + 
        \dfrac{1}{N_\vk}\sum_{\vk_i \in \mathcal{K}}\sum_{i} \left(2\braket{i\vk_i, n\vk | i\vk_i, n\vk} - \braket{i\vk_i, n\vk | n\vk, i\vk_i}\right),
\end{equation}
where $\hat{H}_0$ refers to the single-particle component of the many body Hamiltonian. 

In this paper, we focus on three-dimensional insulating systems with an indirect gap, i.e., $\varepsilon_{a \vk_a}-\varepsilon_{i \vk_i}>0, \forall i,a,\vk_i,\vk_a$.
To simplify the analysis, we assume that the orbitals are exact at any $\vk$ point and the number of virtual orbitals are truncated to a finite number. 
In addition, we assume that the exact orbitals and orbital energies in the TDL are smooth and periodic with respect to their momentum vector index $\vk\in\Omega^*$. 
This assumption is a restriction in our current analysis. For systems free of topological obstructions~\cite{BrouderPanatiCalandraEtAl2007,MonacoPanatiPisanteEtAl2018}, 
these conditions may be replaced by weaker ones using techniques based on Green's functions, or Hamiltonians defined in the atomic orbital basis instead of the band basis~\cite{ScuseriaAyala1999}.

\subsection{CCD theory}
Based on the reference HF determinant $\ket{\Phi}$, the CCD theory represents the wavefunction as 
\begin{align*}
\ket{\Psi} 
& = e^{\mathcal{T}} \ket{\Phi},
\\
\mathcal{T} &= \dfrac{1}{N_\vk}\sum_{ijab}\sum_{\vk_i,\vk_j,\vk_a\in\mathcal{K}} T_{ijab}(\vk_i, \vk_j, \vk_a)a_{a\vk_a}^\dagger a_{b\vk_b}^\dagger a_{j\vk_j} a_{i\vk_i},
\end{align*}
where $a_{n\vk}^\dagger$ and $a_{n\vk}$ are creation and annihilation operators for $\psi_{n\vk}$, 
$T_{ijab}(\vk_i, \vk_j, \vk_a)$ (commonly denoted as $t_{i\vk_i,j\vk_j}^{a\vk_a, b\vk_b}$ in literature) is the \textit{normalized} CCD double amplitude, and 
$\vk_b \in \mathcal{K}$ is uniquely determined using the crystal momentum conservation $\vk_i + \vk_j - \vk_a - \vk_b\in \mathbb{L}^*$.
The amplitude tensor $T^{N_\vk}_* = (T_{ijab}(\vk_i, \vk_j, \vk_a))$ is defined as the root of a nonlinear amplitude equation that consists of constant, linear and quadratic terms.

In practice, the amplitude equation can be solved using a quasi-Newton method \cite{Schneider2009,ShavittBartlett2009}, which is equivalent to applying fixed point iteration to 
\begin{equation}\label{eqn:amplitude_finite}
        T
        = 
        \dfrac{1}{\varepsilon^{N_\vk}}
        \mathcal{A}_{N_\vk}(T).
\end{equation}
Here $\varepsilon^{N_\vk}$ denotes a diagonal operator with entries 
$\varepsilon_{i\vk_i,j\vk_j,a\vk_a,b\vk_b}^{N_\vk} = \varepsilon_{i\vk_i}^{N_\vk} + \varepsilon_{j\vk_j}^{N_\vk} - \varepsilon_{a\vk_a}^{N_\vk} - \varepsilon_{b\vk_b}^{N_\vk}$, 
and $1/\varepsilon^{N_\vk}$ gives the diagonal operator with $1/\varepsilon_{i\vk_i,j\vk_j,a\vk_a,b\vk_b}^{N_\vk}$.
The operator $\mathcal{A}_{N_\vk}(T)$ is referred to as the \textit{ERI-contraction map} (see the definition in \cref{appendix:ccd}). It 
consists of contractions between ERIs and $T$, and does not involve orbital energies.
Note that both $T$ and $\mathcal{A}_{N_\vk}(T)$ are tensors indexed by $(i,j,a,b)$ and $(\vk_i,\vk_j,\vk_a)\in\mathcal{K}\times\mathcal{K}\times\mathcal{K}$.
The CCD correlation energy is then defined as 
\begin{align*}
       E_\text{CCD}^{N_\vk} 
       & = \frac{1}{N_{\vk}^3}
       \sum_{\vk_i\vk_j \vk_a\in\mathcal{K}}\sum_{i j a b}
       W_{ijab}(\vk_i, \vk_j, \vk_a)
       T_{ijab}(\vk_i, \vk_j, \vk_a)
       := \mathcal{G}_{N_\vk}(T^{N_\vk}_*),
\end{align*}
where $W_{ijab}(\vk_i, \vk_j, \vk_a)$ denotes the antisymmetrized ERI, 
$2\braket{i\vk_i,j \vk_j |a\vk_a, b\vk_b}-\braket{i\vk_i,j \vk_j |b\vk_b, a\vk_a}$. 

In the TDL with $\mathcal{K}\rightarrow\Omega^*$, $T^{N_\vk}_* = \{T_{ijab}(\vk_i, \vk_j, \vk_a)\}$ with $(\vk_i, \vk_j, \vk_a) \in \mathcal{K}\times \mathcal{K}\times\mathcal{K}$
converges to $t_* = \{t_{ijab}(\vk_i, \vk_j, \vk_a)\}$, where each $t_{ijab}(\vk_i,\vk_j,\vk_a)$ is a function of $(\vk_i, \vk_j, \vk_a)$ in $\Omega^*\times\Omega^*\times\Omega^*$. 
The converged TDL amplitude $t_*$ satisfies a similar  amplitude equation
\begin{equation}\label{eqn:amplitude_tdl}
        t = 
        \dfrac{1}{\varepsilon^\text{TDL}}
        \mathcal{A}_\text{TDL}(t),
\end{equation}
where $\mathcal{A}_\text{TDL}$ can be obtained from $\mathcal{A}_{N_\vk}$ by taking $\frac{1}{N_\vk}\sum_{\vk\in\mathcal{K}} \rightarrow \frac{1}{|\Omega^*|}\int_{\Omega^*}\ud\vk$ with $N_\vk \rightarrow\infty$.
For example, the \textit{4-hole-2-particle (4h2p) linear term} in $[\mathcal{A}_{N_\vk}(T)]_{ijab}(\vk_i, \vk_j, \vk_a)$ converges in $\mathcal{A}_\text{TDL}(t)$ as 
\begin{align}
    & \dfrac{1}{N_\vk}
    \sum_{\vk_k \in \mathcal{K}}\sum_{kl}
    \braket{k\vk_k, l\vk_l | i \vk_i, j\vk_j}
    T_{klab}(\vk_k, \vk_l, \vk_a) 
    \nonumber \\
    & \xrightarrow{\mathcal{K}\rightarrow\Omega^*}
    \dfrac{1}{|\Omega^*|}
    \int_{\Omega^*}\ud\vk_k
    \sum_{kl}
    \braket{k\vk_k, l\vk_l | i \vk_i, j\vk_j}
    t_{klab}(\vk_k, \vk_l, \vk_a).
    \label{eqn:eri_4h2p}
\end{align}
The image of $\mathcal{A}_\text{TDL}$ is a set of functions of $(\vk_i,\vk_j,\vk_a)$ indexed by $(i,j,a,b)$. 
The CCD correlation energy in the TDL is defined in a similar way as 
\begin{align*}
       E_\text{CCD}^\text{TDL}
       = 
       \dfrac{1}{|\Omega^*|^3}
       \int_{\Omega^*\times \Omega^*\times \Omega^*}\ud\vk_i\ud\vk_j\ud\vk_a
       W_{ijab}(\vk_i, \vk_j, \vk_a)
       t_{ijab}(\vk_i, \vk_j, \vk_a)
       := \mathcal{G}_\text{TDL}(t_*).
\end{align*}

Applying $n$ steps of fixed point iteration over \cref{eqn:amplitude_finite} and \cref{eqn:amplitude_tdl} with zero initial guess (i.e., $t=0,T=0$), we obtain the CCD$(n)$ amplitude and the CCD$(n)$ energy, 
\begin{equation}\label{eqn:ccdn}
    \begin{array}{ll}
        & E_{\text{CCD}(n)}^{N_\vk} = 
    \mathcal{G}_{N_\vk}(T_{n}^{N_\vk})
    \\
    & 
    \hspace{1em}
    \begin{array}{l@{\hskip 0.5em}l}
                T_{m}^{N_\vk} 
                &=
                (\varepsilon^{N_\vk})^{-1}
                \mathcal{A}_{N_\vk}(T_{m-1}^{N_\vk})
                \\
                T_{0}^{N_\vk} & = \bm{0}
    \end{array}
    \end{array}
    \longrightarrow
    \begin{array}{ll}
        & E_{\text{CCD}(n)}^\text{TDL} = \mathcal{G}_\text{TDL}(t_n)
        \\        
        &
        \hspace{1em}
        \begin{array}{l@{\hskip 0.1em}l}
                t_m &= 
                (\varepsilon^\text{TDL})^{-1}
                \mathcal{A}_\text{TDL}(t_{m-1})
                \\
                t_0 & = \bm{0}
        \end{array}.
    \end{array}
\end{equation}
CCD$(n)$ is related to the perturbative description of CCD and consists of a subset of finite order perturbation energy terms in the M{\o}ller-Plesset perturbation theory. 
For example, CCD(1) can be identified with MP2, and CCD(2) contains all the terms in MP2 and MP3, and a subset of terms in MP4.

One main result of this paper is the rigorous analysis of the finite-size error in CCD($n$) calculation with any fixed $n > 0$. 
If the fixed point iterations in both the finite and TDL cases converge to the CCD amplitudes with $n\rightarrow \infty$, i.e., $T_n^{N_\vk}\rightarrow T_*^{N_\vk}$ and $t_{n} \rightarrow t_*$ 
(the technical definition of this convergence is provided in \cref{app:error_ccd}),
the finite-size error analysis for CCD($n$) also applies to the converged CCD calculation. 
In other words, we analyze the finite-size error in CCD calculation using a perturbative approach based on the analysis on CCD$(n)$. 

\subsection{Madelung constant correction}

Ref.~\cite{XingLin2023} shows that the finite-size errors in CCD($n$) and CCD both scale as 
$\Or(N_\vk^{-\frac13})$ when assuming that exact orbital energies are used in the amplitude equation \cref{eqn:amplitude_finite}. 
The same finite-size scaling also appears in Fock exchange energy and occupied orbital energy calculations.
One common correction to reduce the finite-size errors in the latter two calculations is to add a Madelung-constant shift 
\cite{FraserFoulkesRajagopalEtAl1996,ChiesaCeperleyMartinEtAl2006,DrummondNeedsSorouriEtAl2008} to the Ewald kernel. 
This shift introduces a correction to all involved ERIs in the calculations as
\begin{equation*}\label{eqn:eri_correction}
        \braket{n_1\vk_1,n_2\vk_2|n_3\vk_3,n_4\vk_4} - \delta_{n_1n_3}\delta_{n_2n_4}\delta_{\vk_1\vk_3}\delta_{\vk_2\vk_4} N_\vk\xi. 
\end{equation*}
Such a correction is triggered only in ERIs which have fully matched orbital indices, i.e., $n_1 = n_3, n_2 = n_4$, and zero momentum transfer, i.e., $\vk_1 = \vk_3$.
The Madelung constant $\xi$ is defined uniquely by the unit cell and the $\vk$-point mesh $\mathcal{K}$ as 
\begin{equation}\label{eqn:madelung_practice}
        \xi
        = 
        \left(
        \dfrac{1}{N_\vk}\sum_{\vq\in\mathcal{K}_\vq}
        - 
        \dfrac{1}{|\Omega^*|}\int_{\Omega^*}\ud\vq
        \right)
        \xsum_{\vG \in \mathbb{L}^*}\dfrac{4\pi}{|\Omega|}\dfrac{e^{-\sigma|\vq+\vG|^2}}{|\vq + \vG|^2}
        -\dfrac{4\pi \sigma}{N_\vk}
        + 
        \xsum_{\vR \in \mathbb{L}_{\mathcal{K}}}\dfrac{\text{erfc}\left(\sigma^{-1/2}|\vR|/2\right)}{|\vR|},
\end{equation}
where $\mathcal{K}_\vq$ is a uniform mesh that is of the same size as $\mathcal{K}$ and contains $\vq = \bm{0}$ and $\mathbb{L}_{\mathcal{K}}$
is the real-space lattice associated with the  the reciprocal-space lattice $\vq + \vG$ with $\vq \in \mathcal{K}_\vq, \vG\in\mathbb{L}$. 
Note that $\xi$ does not vary with respect to $\sigma> 0$ \cite{FraserFoulkesRajagopalEtAl1996} and this parameter $\sigma$ is commonly 
tuned to control the lattice cutoffs in the summation over $\mathbb{L}^*$ and $\mathbb{L}_{\mathcal{K}}$ in \cref{eqn:madelung_practice} when numerically computing $\xi$.

For finite-size orbital energy calculation in \cref{eqn:orbital_energy}, the Madelung constant correction gives  
\begin{equation}\label{eqn:madelung_orbitalenergy}
        \varepsilon_{n\vk}^{N_\vk, \xi} =
        \begin{cases}
                \varepsilon_{n\vk}^{N_\vk} + \xi & n \text{ is occupied},\\
                \varepsilon_{n\vk}^{N_\vk} & n \text{ is virtual}.
        \end{cases} 
\end{equation}

For the ERI contractions in $\mathcal{A}_{N_\vk}(T)$, this correction when applied will be triggered in six linear amplitude terms (see \cref{proof:error_amplitude_contraction}). 
For example, the 4h2p linear term in \cref{eqn:eri_4h2p} can trigger the correction and be modified to 
\begin{equation}\label{eqn:4h2p_xi}
        \dfrac{1}{N_\vk} \sum_{\vk_k\in\mathcal{K}}\sum_{kl}\braket{k\vk_k, l\vk_l | i\vk_i, j\vk_j} T_{klab}(\vk_k, \vk_l, \vk_a)
        - \xi T_{ijab}(\vk_i, \vk_j, \vk_a).
\end{equation}
Collecting the corrections to all the six terms together, $\mathcal{A}_{N_\vk}(T)$ is modified to  
\begin{equation}\label{eqn:madelung_ericontraction}
        \mathcal{A}_{N_\vk,\xi}(T) = 
        \mathcal{A}_{N_\vk}(T) + 2\xi T. 
\end{equation}
In the finite-size CCD$(n)$ and converged CCD calculations, the Madelung constant correction
can be applied to the orbital energies, ERI contractions, or both in the amplitude equation \cref{eqn:amplitude_finite}. 
As a result, we may have three correction schemes compared to the standard calculation without any correction. 

Particularly, applying the Madelung constant correction to both components gives the amplitude equation
\begin{equation}\label{eqn:amplitude_finite_xi}
        T = \dfrac{1}{\varepsilon^{N_\vk} + 2\xi} (\mathcal{A}_{N_\vk}(T) + 2\xi T) 
                = \dfrac{1}{\varepsilon^{N_\vk, \xi}} \mathcal{A}_{N_\vk,\xi}(T).
\end{equation}
It can be easily verified that this amplitude equation has the same solution as \cref{eqn:amplitude_finite} without any correction.
Its CCD solution is thus identical to the one without correction.
However, the associated CCD$(n)$ calculation differs and can be interpreted 
as solving the original amplitude equation \cref{eqn:amplitude_finite} using a quasi-Newton method that  
has a $2\xi$-diagonal shift to the Jacobian matrix.

\section{Main statements}

We start our error analysis with the CCD$(n)$ calculation with a fixed number of iterations $n>0$ and then generalize the analysis to the fully converged CCD calculation. 
In CCD$(n)$, the finite-size error is quantified by  $E_{\text{CCD}(n)}^\text{TDL} - E_{\text{CCD}(n)}^{N_\vk}$. 
According to \cref{eqn:ccdn}, this error can be traced back to two sources: the difference in energy operators, $\mathcal{G}_\text{TDL}$ versus $\mathcal{G}_{N_\vk}$, and the difference in amplitudes, $t_n$ versus $T_n^{N_\vk}$. 
Recall that $t_n$ and $T_n^{N_\vk}$ are the amplitudes of the system in TDL and finite-size cases, respectively.
Let us consider the \textit{evaluation map}, denoted by $\mathcal{M}_\mathcal{K}$. 
This map evaluates a tensor-valued function, initially defined on the product space $\Omega^*\times \Omega^* \times \Omega^*$, on a finite-sized grid  $\mathcal{K}\times\mathcal{K}\times\mathcal{K}$. 
Consequently, the values of the TDL amplitude $t_n$ on this finite-size grid, $\mathcal{K}\times\mathcal{K}\times\mathcal{K}$, are given by $\mathcal{M}_\mathcal{K} t_n$ and are approximated by the finite-size amplitude $T_{n}^{N_\vk}$.

By applying the triangle inequality, we can decompose the finite-size error into two distinct sources: the errors arising from the discretized energy calculation using the exact amplitude $\mathcal{M}_\mathcal{K}t_n$, and the errors stemming from the amplitude calculation itself.
\begin{align}
        \left|
                E_{\text{CCD}(n)}^\text{TDL} - E_{\text{CCD}(n)}^{N_\vk}
        \right|
        & \leqslant 
        \left|
                \mathcal{G}_\text{TDL}(t_n) - \mathcal{G}_{N_\vk}(\mathcal{M}_{\mathcal{K}}t_n)
        \right|
        +
        \left|
                \mathcal{G}_{N_\vk}(\mathcal{M}_{\mathcal{K}}t_n) - 
                \mathcal{G}_{N_\vk}(T_{n}^{N_\vk})
        \right|
        \nonumber \\
        & \leqslant 
        \left|
                \mathcal{G}_\text{TDL}(t_n) - \mathcal{G}_{N_\vk}(\mathcal{M}_{\mathcal{K}}t_n)
        \right|
        +
        \dfrac{C}{N_\vk^3}\sum_{ijab}\sum_{\vk_i,\vk_j,\vk_a\in\mathcal{K}}
        \left|
                [\mathcal{M}_\mathcal{K}t_n - T_{n}^{N_\vk}]_{ijab}(\vk_i, \vk_j, \vk_a)
        \right|.
        \label{eqn:decomposition_energy}
\end{align}
Here we use the fact that $\abs{W_{ijab}(\vk_i, \vk_j,\vk_a)}$ can be upper bounded uniformly by a constant.

To further break down the error in amplitude calculation, we note that $t_n$ and $T_n^{N_\vk}$ are recursively 
constructed by \cref{eqn:ccdn} with initial values $t_0 = \bm{0}$ and $T_0^{N_\vk} = \bm{0}$. 
As a result, the error in the CCD$(n)$ amplitude calculation can also be recursively decomposed using the same strategy above as
\begin{align}
        \mathcal{M}_{\mathcal{K}}t_{n} - T_{n}^{N_\vk}
        & = \dfrac{1}{\varepsilon^\text{TDL}}
        \left[
        \mathcal{M}_{\mathcal{K}}
        \mathcal{A}_\text{TDL}(t_{n-1})
        -
        \mathcal{A}_{N_\vk}(\mathcal{M}_{\mathcal{K}}t_{n-1})
        \right]
        \nonumber\\
        & \hspace*{1em}
        + 
        \mathcal{A}_{N_\vk}(\mathcal{M}_{\mathcal{K}}t_{n-1})
        \left[
        \dfrac{1}{\varepsilon^\text{TDL}}
        -
        \dfrac{1}{\varepsilon^{N_\vk}}
        \right]
        \nonumber\\
        & \hspace*{1em}
        +
        \dfrac{1}{\varepsilon^{N_\vk}}
        \left[
        \mathcal{A}_{N_\vk}( \mathcal{M}_{\mathcal{K}}t_{n-1})
        -
        \mathcal{A}_{N_\vk}(T_{n-1}^{N_\vk})
        \right].
        \label{eqn:decomposition_amplitude}
\end{align}
The three error terms from this dissection can be interpreted as the errors in ERI contractions using exact CCD($n-1$) amplitudes, orbital energies, and CCD$(n-1)$ amplitude calculation composed with $\mathcal{A}_{N_\vk}$, respectively. 
For the last term, it can be shown that entries in $\mathcal{A}_{N_\vk}( \mathcal{M}_{\mathcal{K}}t_{n-1})-\mathcal{A}_{N_\vk}(T_{n-1}^{N_\vk})$ 
have the same scaling with respect to $N_\vk$ as those in $\mathcal{M}_{\mathcal{K}}t_{n-1} - T_{n-1}^{N_\vk}$.  
Replacing the last term by $\mathcal{M}_{\mathcal{K}}t_{n-1} - T_{n-1}^{N_\vk}$ and applying the dissection recursively, 
we find that the error in the CCD$(n)$ amplitude calculation is determined by those in the ERI contractions and orbital energies, i.e., the first two terms in \cref{eqn:decomposition_amplitude}.

Overall, the finite-size error of CCD$(n)$ calculation can be decomposed into errors in three basic factors: 
(1) energy calculation using exact CCD$(n)$ amplitude, (2) ERI contractions using exact CCD$(n-1)$ amplitude, and 
(3) orbital energies. 
This error decomposition is also valid when applying the Madelung constant correction to orbital energies (\cref{eqn:madelung_orbitalenergy}) or ERI contractions (\cref{eqn:madelung_ericontraction}). 
By analyzing these three error sources with or without corrections separately, we can obtain the finite-size error estimate for CCD$(n)$ with various correction schemes.


The Madelung constant correction can reduce the finite-size error in orbital energies from $\Or(N_\vk^{-\frac13})$ to $\Or(N_\vk^{-1})$. 
This correction is at the HF level.
One \textit{main technical result} of this work is to show that the Madelung constant correction can also reduce the 
finite-size error in \textit{most (but not all) entries} of the ERI contraction (note that $\mathcal{M}_{\mathcal{K}}\mathcal{A}_\text{TDL}(t)-\mathcal{A}_{N_\vk}(\mathcal{M}_{\mathcal{K}}t)$ is a tensor)
from $\Or(N_\vk^{-\frac13})$ to $\Or(N_\vk^{-1})$, see \cref{sec:madelung_idea}.
As a result, when applying the Madelung constant correction to \textit{both} orbital energies \textit{and} ERI contractions in the CCD$(n)$ calculation, 
the overall finite-size error can be successfully reduced to $\Or(N_\vk^{-1})$. 
\begin{thm}\label{thm:error_ccdn}
        In CCD($n$) calculation, the finite-size error in the correlation energy scales as $\Or(N_\vk^{-\frac13})$ in each of the following scenarios 
        (1) there is no finite-size correction, (2) the Madelung constant correction is only applied to the ERI contraction $\mathcal{A}_{N_\vk}$, 
        and (3) the Madelung constant correction is only applied to the orbital energy $\varepsilon_{n\vk}^{N_\vk}$. 

        When the Madelung constant correction is applied to both $\mathcal{A}_{N_\vk}$ and $\varepsilon_{n\vk}^{N_\vk}$ in the CCD$(n)$ calculation, 
        the overall finite-size error scales as $\Or(N_\vk^{-1})$.
        
        As a special case, the same conclusion applies to MP3 calculations.
\end{thm}
\begin{proof}
        See \cref{app:error_ccdn}.
\end{proof}

Now what happens to the finite-size error of the CCD calculation as $n\to \infty$? 
For gapless and small-gap systems, it has been observed in practice that the fixed point iteration might not converge 
or the amplitude equation might have multiple solutions. 
While we focus on systems with an indirect gap, it is worth noting that even in such favorable scenarios, 
existence and uniqueness of solutions to the CCD amplitude equations in both finite and TDL cases remain an open question and are beyond the scope of this paper.

To study the finite-size error of CCD via the above results on CCD($n$), we make additional technical assumptions 
(see \cref{app:error_ccd}) that can guarantee the convergence of CCD$(n)$ to CCD.
Under these assumptions, we show that the finite-size scaling of CCD calculation is upper bounded by those of its associated converging CCD$(n)$ calculations. 
Numerical observations (see \cref{sec:numerics}) further show that this finite-size scaling estimate through CCD$(n)$ is asymptotically sharp for 
CCD calculation with the Madelung constant correction applied to orbital energies, ERI contractions, or both.

\begin{thm}\label{thm:error_ccd_madelung}
        (Informal) Under additional conditions on the convergence of CCD$(n)$ to CCD, 
        the finite-size error of the CCD correlation energy scales as $\Or(N_\vk^{-\frac13})$
        in each of the following scenarios 
        (1) the Madelung constant correction is only applied to the ERI contraction $\mathcal{A}_{N_\vk}$ 
        (2) the Madelung constant correction is only applied to the orbital energy $\varepsilon_{n\vk}^{N_\vk}$.

        When the Madelung constant correction is applied to both $\mathcal{A}_{N_\vk}$ and $\varepsilon_{n\vk}^{N_\vk}$ in the CCD calculation, 
        the overall finite-size error scales as $\Or(N_\vk^{-1})$.
\end{thm}

Compared to \cref{thm:error_ccdn}, \cref{thm:error_ccd_madelung} does not address the finite-size error of the CCD calculation 
when no finite-size correction is applied. The natural conclusion from \cref{thm:error_ccdn} is that this finite-size scaling should be  $\Or(N_\vk^{-\frac13})$ as well. However, this error estimate is loose and inconsistent with the numerical observations which suggest an inverse volume scaling.
To obtain a tight estimate, we now use the observation that the CCD calculation without any correction is equivalent to the 
CCD calculation with the Madelung constant correction applied to both ERI contractions and orbital energies. 
Specifically, when applying the Madelung constant corrections \cref{eqn:madelung_orbitalenergy} and \cref{eqn:madelung_ericontraction}, the CCD 
amplitude equation in \cref{eqn:amplitude_finite_xi} can be formulated as 
\begin{equation*}
        \left(
                \varepsilon^{N_\vk} + 2\xi
        \right) 
        T 
        = 
        \mathcal{A}_{N_\vk}(T) + 2\xi T 
        \quad
        \Longleftrightarrow
        \quad
        \varepsilon^{N_\vk}
        T 
        = 
        \mathcal{A}_{N_\vk}(T).
\end{equation*}
With $2\xi T$ on both sides cancelling each other, this reformulation is exactly reduced to the original amplitude equation \cref{eqn:amplitude_finite} without corrections.
Therefore, the roots of the two amplitude equations with and without the corrections are the same, 
and the correlation energy in \textit{converged} CCD calculation without any corrections is the same as that with the Madelung constant correction. 
In other words, when investigating the finite-size error of CCD calculation without corrections, we should apply \cref{thm:error_ccd_madelung} with the Madelung constant correction applied to both $\mathcal{A}_{N_\vk}$ and $\varepsilon_{n\vk}^{N_\vk}$.
This yields a sharp estimate $\Or(N_\vk^{-1})$ and explains the origin of the inverse volume scaling of the finite-size error.


\begin{cor}\label{thm:error_ccd_inexact}
        Under the same additional conditions as in \cref{thm:error_ccd_madelung}, the finite-size error of the CCD correlation energy without finite-size correction scales as $\Or(N_\vk^{-1})$.
\end{cor}

We provide the proof of \cref{thm:error_ccd_madelung} in \cref{app:error_ccd}, and \cref{thm:error_ccd_inexact} follows directly from \cref{thm:error_ccd_madelung} and the reasonings above.

\section{Key ideas}

As demonstrated in the previous section, the finite-size errors in CCD$(n)$ and CCD calculations can be 
reduced to the errors in three simpler basic calculations: energy calculation using exact amplitudes, 
ERI contraction using exact amplitudes, and orbital energies. 
A key observation is that the finite-size errors in all the three calculations can be interpreted and analyzed from a numerical quadrature perspective. 

For a function $g$ over a hypercube $V$, we denote a (generalized) trapezoidal rule using a uniform mesh $\mathcal{X}$ of $V$ and its quadrature error as
\begin{align*}
        \mathcal{Q}_{V}(g,\mathcal{X}) & = \dfrac{|V|}{|\mathcal{X}|}\sum_{\vx\in\mathcal{X}}g(\vx),
        \\
        \mathcal{E}_{V}(g,\mathcal{X}) & = \int_{V}\ud\vx g(\vx) - \mathcal{Q}_{V}(g,\mathcal{X}).
\end{align*}
Under the assumption of exact orbitals at any $\vk$ point in finite-size calculation, the errors in 
(1) orbital energy, (2) energy, and (3) ERI contraction $\mathcal{A}_{N_\vk}$ 
can be respectively formulated by their definitions as
\begin{equation}
\begin{split}
    \varepsilon_{n\vk}^\text{TDL} - \varepsilon_{n\vk}^{N_\vk} 
    & = 
    \dfrac{1}{|\Omega^*|}
    \mathcal{E}_{\Omega^*}
    \left(
        \sum_{i}
        W_{inin}(\vk_i, \vk, \vk_i), 
        \mathcal{K}
    \right),
    \\
    \mathcal{G}_\text{TDL}(t) - \mathcal{G}_{N_\vk}(\mathcal{M}_{\mathcal{K}}t)
    & = 
    \dfrac{1}{|\Omega^*|^3}
    \mathcal{E}_{\Omega^*\times\Omega^*\times\Omega^*}
    \left(
    \sum_{ijab}
    W_{ijab} t_{ijab}(\vk_i, \vk_j, \vk_a), 
    \mathcal{K}\times\mathcal{K}\times\mathcal{K}
    \right), 
    \\
    \mathcal{M}_{\mathcal{K}}
    \mathcal{A}_\text{TDL}(t)
    -
    \mathcal{A}_{N_\vk}(\mathcal{M}_{\mathcal{K}}t)
    & = 
    \dfrac{1}{|\Omega^*|}
    \mathcal{E}_{\Omega^*}
    \left(
            \sum_{kl} \braket{k\vk_k, l\vk_l | i\vk_i, j\vk_j} t_{klab}(\vk_k, \vk_l, \vk_a), 
            \mathcal{K}
    \right) + \cdots.
\end{split}
\label{eqn:quadrature_basic}
\end{equation}
For ERI contraction in the last line, the error is denoted by a tensor with indices $(i, j, a, b, \mathbf{k}_i, \mathbf{k}_j, \mathbf{k}_a)$, 
and the error for the 4h2p linear term \cref{eqn:eri_4h2p} is detailed with $\mathbf{k}_k$ being the integration variable. 
Meanwhile, terms not shown account for errors from computing other linear and quadratic amplitude terms. 
From a numerical quadrature perspective, all the three finite-sized calculations approximate the corresponding integrals in the TDL 
using the trapezoidal rule and a finite mesh $\mathcal{K}$ to discretize $\Omega^*$. 
Therefore, their finite-size errors can be estimated systematically by quadrature error analysis.

In general, the quadrature error associated with a trapezoidal rule is influenced by the integrand's smoothness and boundary conditions. 
If we take $\Delta h$ as the mesh size along each dimension, the quadrature error for a smooth integrand is typically of order $\Or(\Delta h^2)$. 
Interestingly, if the integrand is also periodic, the error diminishes much more rapidly than it does in a nonperiodic scenario. 
The decay rate is faster than any finite power of $\Delta h$, showcasing a super-algebraic decay~\cite{TrefethenWeideman2014}. 
However, for an integrand periodic but marked by singularities, its quadrature error tends to taper off at a slower rate.

For quadrature errors in \cref{eqn:quadrature_basic}, the involved integrands are all periodic across their integration domains. 
However, many of these integrands have singularities within the domains, affecting the scaling of their quadrature errors.  
As to be demonstrated next, these integrands include point singularities arising from both ERIs and amplitudes, 
resulting in low-order power-law decay of the corresponding quadrature errors as $N_\vk$ increases.

\subsection{Singularity structure and quadrature error estimate}


        
        

All the quadrature errors in \cref{eqn:quadrature_basic} have integrands comprising of either ERIs or contractions between ERIs and exact CCD$(n)$ amplitudes. 
The integration variables are momentum vectors sampled in the Brillouin zone $\Omega^*$. 
Consequently, understanding the singularity structure of ERIs and exact CCD$(n)$ amplitudes with respect to their momentum vector indices is crucial in 
comprehending the singularity structure of these integrands and ultimately estimating the quadrature error in the three basic calculations.

First consider a generic ERI $\braket{n_1\vk_1,n_2\vk_2|n_3\vk_3,n_4\vk_4}$ with fixed band indices $(n_1,n_2,n_3,n_4)$ and treat it as a function of $\vk_1, \vk_2$, and 
$\vq = \vk_3 - \vk_1$ in $\Omega^*$. 
By its definition, the ERI can be separated as
\begin{equation}\label{eqn:eri_expansion}
        \braket{n_1\vk_1,n_2\vk_2|n_3\vk_3,n_4\vk_4}  
        = 
        \frac{4\pi}{\abs{\Omega}} 
        \frac{
                \hat{\varrho}_{n_1\vk_1,n_3(\vk_1 + \vq)}(\bm{0}) \hat{\varrho}_{n_2\vk_2,n_4(\vk_2 - \vq)}(\bm{0})
        }{
                |\vq|^2
        }
        + 
        \frac{4\pi}{\abs{\Omega}}
        \sum_{\vG\in\mathbb{L}^*\setminus \{\bm{0}\}}
        \dfrac{\cdots}{|\vq+\vG|^2}.
\end{equation} 
Since we assume all orbitals $\psi_{n\vk}$ periodic and smooth with respect to $\vk \in \Omega^*$, this ERI has a point singularity 
at $\vq = \bm{0}$ only due to the first fraction term. 
Specifically, this term is of fractional form $f(\vk_1, \vk_2, \vq) /|\vq|^2$ with a smooth numerator $f$. 
As can be verified by direct calculation, such a fraction term can have its point singularity at $\vq = \bm{0}$ described by the following general concept called 
algebraic singularity. 

\begin{defn}\label{def:fractional}
        (Informal) 
        A function $f(\vx)$ has \textbf{algebraic singularity of order} $\gamma\in\RR$ at $\vx_0 \in \mathbb{R}^d$ if
        everty $l$-th order derivatives near $\vx_0$ is bounded asymptotically by $|\vx - \vx_0|^{\gamma - l}$, i.e., 
        \[
        \left|
        \dfrac{\partial^\valpha}{\partial \vx^{\valpha}} f(\vx)  
        \right|
        \leqslant 
        C |\vx - \vx_0|^{\gamma - |\valpha|},
        \ \forall \valpha \geqslant 0,
        \]
        where $\valpha$ denotes a non-negative $d$-dimensional derivative multi-index.  
        For brevity, $f$ is said to be singular at $\vx_0$ with order $\gamma$.
        See \cref{def:fractional_app} in Appendix for the rigorous mathematical definition.
\end{defn}


A representative example of such a singular function with order $\gamma$ is $g(\vx)/|\vx|^2$ where $g(\vx)$ is smooth and scale as $\Or(|\vx|^{2+\gamma})$ 
near $\vx = \bm{0}$.
Using orbital orthogonality, the generic ERI exhibits singularities at $\vq = \bm{0}$ with order $0$ when band indices mismatch ($n_1\neq n_3$, $n_2 \neq n_4$), 
$-1$ when partially match ($n_1=n_3, n_2\neq n_4$ or $n_1\neq n_3, n_2 = n_4$), and $-2$ when fully match ($n_1=n_3, n_2=n_4$). 

We can now characterize the singularity structure of integrand $\sum_{i}W_{inin}(\vk_i, \vk, \vk_i)$ defined by orbital energy $\varepsilon_{n\vk}^\text{TDL}$ 
in \cref{eqn:quadrature_basic}. 
Fixing $(n, \vk)$, the leading singularity in each $W_{inin}(\vk_i, \vk, \vk_i)$ (as a function of $\vk_i$) comes from the exchange term $\braket{i\vk_i, n\vk | n\vk, i\vk_i}$ which has algebraic singularity at $\vk_i = \vk$
with order $0$ when $i \neq n$ and $-2$ when $i =n$. 
In computing $\varepsilon_{n\vk}^\text{TDL}$, the overall integrand is thus singular at \(\vk_i = \vk\)
with order \(0\) for a virtual band \(n\) and \(-2\) for an occupied band \(n\).

For such periodic functions with one point of algebraic singularity, the conventional textbook analysis of the trapezoidal rule is overly pessimistic. 
A key technical tool in this work is \cref{lem:quaderror} below, which provides a rigorous and sharp quadrature error estimate linking its error scaling to the singularity order.
(See \cref{lem:quaderror1} in Appendix for a more general statement.)

\begin{lem}
        \label{lem:quaderror}
        Let $f(\vx)$ be periodic with respect to $V = [-\frac12,\frac12]^d$ and smooth everywhere except at $\vx = \bm{0}$ with order $\gamma \geqslant -d+1$.     
        At $\vx = \bm{0}$, $f(\vx)$ is set to $0$.
        The quadrature error of a trapezoidal rule using an $m^d$-sized uniform mesh $\mathcal{X}$ that contains $\vx = \bm{0}$ can be estimated as 
        \[
                \left|\mathcal{E}_V(f, \mathcal{X})\right| \leqslant C m^{-(d + \gamma)}.
        \]
        If $f(\bm{0})$ is set to an $\Or(1)$ value in the calculation,
        it introduces additional $\Or(m^{-d})$ quadrature error.
\end{lem}
Combining this error estimate with the integrand singularity structure for orbital energies (with $m = N_\vk^{\frac13}$, $d = 3$, and $\gamma = -2, 0$ for occupied and virtual orbitals, respectively),
we obtain the finite-size error estimate for orbital energy calculation as 
\begin{equation}\label{eqn:error_estiamte_orbitalenergy}
        \left|
        \varepsilon_{n\vk}^\text{TDL} - \varepsilon_{n\vk}^{N_\vk} 
        \right|
        \leqslant
        C 
        \left\{
        \begin{array}{ll}
                N_\vk^{-\frac13} & n \text{ is occupied} \\
                N_\vk^{-1} & n \text{ is virtual} 
        \end{array}
        \right.
        .
\end{equation}

To adapt the above approach for examining energy and ERI contraction calculations with exact amplitude, 
we need to characterize the singularity structure of the exact CCD$(n)$ amplitude entries. First we note that the exact CCD$(1)$ amplitude is just the MP2 amplitude 
\[
\left\{\braket{a\vk_a, b\vk_b | i\vk_i, j\vk_j} (\varepsilon_{i\vk_i,j\vk_j,a\vk_a,b\vk_b}^\text{TDL})^{-1}\right\}.
\] 
As a function of $(\vk_i, \vk_j, \vk_a)$, each CCD$(1)$ amplitude entry indexed by $(i,j,a,b)$ has the same singularity structure as the included 
ERI term which is singular at $\vk_a - \vk_i = \bm{0}$ with order $0$.  
It turns out that the exact CCD$(n)$ amplitude with any $n > 0$ all has the same singularity structure as $\braket{a\vk_a, b\vk_b | i\vk_i, j\vk_j}$.  
\begin{lem}[Singularity structure of the amplitude, Lemma 4 in \cite{XingLin2023}]
        \label{lem:thm1_nonsmooth}
        In CCD$(n)$ calculation with $n > 0$, each entry of the exact double amplitude $t_n = \left\{[t_n]_{ijab}(\vk_i, \vk_j, \vk_a)\right\}$ belongs to the following function space
        \begin{align*}
                \mathbb{T}(\Omega^*) = \big\{
                f(\vk_i, \vk_j, \vk_a):\  
                & f \text{ is periodic with respect to } \vk_i, \vk_j, \vk_a \in \Omega^*,\\
                & f \text{ is smooth everywhere except at } \vk_a = \vk_i 
                \text{ with algebraic singularity of order } 0,\\
                & f \text{ is smooth with respect to } \vk_i, \vk_j \text{ at the singularity }\vk_a = \vk_i
                \big\}.
        \end{align*}
\end{lem}

Combining the above singularity structure characterizations of ERIs and exact CCD$(n)$ amplitudes, 
we are able to analyze the integrands for energy and 
ERI contraction calculations in \cref{eqn:quadrature_basic}. 
Take the CCD$(n)$ exchange energy term as an example whose finite-size error can be formulated as 
\begin{equation*}
    -\dfrac{1}{|\Omega^*|^3}
    \mathcal{E}_{\Omega^*\times\Omega^*\times\Omega^*}
    \left(
    \sum_{ijab}
    \braket{i\vk_i, j\vk_j | b\vk_b,a\vk_a} t_{ijab}(\vk_i, \vk_j, \vk_a), 
    \mathcal{K}\times\mathcal{K}\times\mathcal{K}
    \right), 
\end{equation*}
with integration variables $(\vk_i, \vk_j, \vk_a)$.
For each set of $(i,j,a,b)$, both the associated ERI and exact amplitude in the integrand exhibit an algebraic singularity of order $0$ at $\vk_a - \vk_j = \bm{0}$ 
and $\vk_a - \vk_i = \bm{0}$, respectively. 
For such a product of two functions with algebraic singularities, we also provide a rigorous quadrature error estimate similar to \cref{lem:quaderror}. 
Specifically for all the integrands in the energy and ERI contraction calculations, their quadrature errors are determined 
by the most singular product components in the integrand. 
The quadrature error scaling still scales as $m^{-(d+\gamma)}$ similar to \cref{lem:quaderror} but with $\gamma$ denoting 
the minimum algebraic singularity order of all ERIs and exact amplitudes.

For the exchange term above, all involved ERIs and exact amplitudes have point singularities of order $0$, and similar for the direct term. Thus, we can get the 
finite-size error estimate for the energy calculation using exact CCD$(n)$ amplitude as 
\begin{equation}
        \left|
                \mathcal{G}_\text{TDL}(t) - \mathcal{G}_{N_\vk}(\mathcal{M}_{\mathcal{K}}t)
        \right|
        \leqslant 
        C N_\vk^{-1}.
\end{equation}

Similar analysis is also applicable to ERI contractions using exact amplitudes. 
The key distinction lies in the fact that the ERI contractions involve integrands formed by ERIs with stronger singularities, such as those defined by particle-particle or hole-hole diagrams. 
A prominent example is the 4h2p linear term in \cref{eqn:eri_4h2p}, where the involved ERI $\braket{k\vk_k, l\vk_l | i\vk_i, j\vk_j}$ exhibits an algebraic singularity of order 
$-2$ at $\vk_k = \vk_i$ when $k = i$ and $l = j$.
Consequently, as per the above analysis, the finite-size error in the 4h2p linear term calculation alone scales as $\Or(N_\vk^{-\frac13})$. 
This error turns out to dominate the overall finite-size error in the ERI contraction calculation, and we have 
\begin{equation}
        \left|
                \left[
                \mathcal{M}_{\mathcal{K}}
                \mathcal{A}_\text{TDL}(t)
                -
                \mathcal{A}_{N_\vk}(\mathcal{M}_{\mathcal{K}}t)
                \right]_{ijab}(\vk_i, \vk_j, \vk_a)
        \right|
        \leqslant 
        C N_\vk^{-\frac13},
        \quad 
        \forall 
        i,j,a,b, \ 
        \forall 
        \vk_i,\vk_j,\vk_a \in \mathcal{K}. 
\end{equation}

\subsection{Madelung constant correction as a quadrature error reduction method}\label{sec:madelung_idea}
    
To reduce the quadrature error for a singular integrand, one common numerical quadrature technique is the singularity subtraction method~\cite{DavisRabinowitz2007}. 
Essentially, this method involves constructing an auxiliary function 
$h$ that possesses the same leading singularity as the integrand 
$g$. The integral is then approximated as:
\begin{align}
        \int_{V}\ud\vx g(\vx) 
        & \approx \frac{|V|}{|\mathcal{X}|}\sum_{\vx \in \mathcal{X}}(g - h)(\vx) + \int_{V}\ud\vx h(\vx)
        \nonumber\\
        & = \frac{|V|}{|\mathcal{X}|}\sum_{\vx \in \mathcal{X}}g(\vx) + \left(\int_{V}\ud\vx  - \dfrac{|V|}{|\mathcal{X}|}\sum_{\vx \in \mathcal{X}}\right)h(\vx).
        \label{eqn:singularity_subtraction}
\end{align}
This approximation consists of the numerical quadrature of $g - h$ and the exact integral of $h$ (which can be computed analytically or numerically with high precision). 
It is also equivalent to adding a correction $\mathcal{E}_{V}(h, \mathcal{X})$ to the numerical quadrature of $g$. 
By this correction, the quadrature error changes from $\mathcal{E}_V(g, \mathcal{X})$ to $\mathcal{E}_{V}(g - h, \mathcal{X})$.
Since $h$ removes the leading singularity of $g$ in the subtraction $g - h$, 
$\mathcal{E}_{V}(g - h, \mathcal{X})$ can be asymptotically smaller than $\mathcal{E}_V(g, \mathcal{X})$.

The Madelung constant defined in \cref{eqn:madelung_practice} can be reformulated using an arbitrary $\sigma > 0$ as 
\begin{equation}\label{eqn:madelung_sigma}
        \xi =
        -\dfrac{1}{|\Omega^*|}
        \left(
        \int_{\Omega^*}\ud\vq
        -
        \dfrac{|\Omega^*|}{N_\vk}\sum_{\vq\in\mathcal{K}_\vq}
        \right)
        \sum_{\vG \in \mathbb{L}^*}\dfrac{4\pi}{|\Omega|}\dfrac{e^{-\sigma|\vq+\vG|^2}}{|\vq + \vG|^2}
        + \Or(N_\vk^{-1}).
\end{equation}
Compared to \cref{eqn:singularity_subtraction}, this representation connects $\xi$ to the singularity subtraction correction defined by an auxiliary function
\[
        h_\sigma(\vq)
        =
        \sum_{\vG \in \mathbb{L}^*}\dfrac{4\pi}{|\Omega|}\dfrac{e^{-\sigma|\vq+\vG|^2}}{|\vq + \vG|^2}.
\] 
The effectiveness of the Madelung constant correction to reduce the finite-size error can be rigorously explained by this connection.


Taking the occupied orbital energy $\varepsilon_{n\vk}$ as an example, its exchange portion leads to the dominant finite-size error and the associated Madelung constant correction 
modifies the calculation as (with a change of variable $\vk_i \rightarrow \vk - \vq$)
\[
        \dfrac{1}{N_\vk}\sum_{\vq\in\mathcal{K}_\vq} \left(\sum_{i}\braket{i(\vk - \vq), n\vk | n\vk i(\vk - \vq)}\right) - \xi. 
\]
Comparing this calculation with \cref{eqn:singularity_subtraction} and \cref{eqn:madelung_sigma}, the Madelung constant correction exactly uses the auxiliary function 
$h_\sigma(\vq)$ to remove the leading singularity (i.e., $4\pi/|\Omega||\vq|^{-2}$) of the target integrand, and thus reduces the associated finite-size error asymptotically to $\Or(N_\vk^{-1})$.

One major technical contribution of this paper is to rigorously prove the effectiveness of the Madelung constant correction for reducing the finite-size error in the ERI contractions 
calculations following the same singularity subtraction interpretation.
For instance, consider the 4h2p linear term calculation in \cref{eqn:eri_4h2p} with any fixed entry index $(i,j,a,b,\vk_i,\vk_j,\vk_a)$. 
Using the change of variable $\vk_k \rightarrow \vk_i - \vq$, this term with the Madelung constant correction can be detailed as 
\begin{equation}\label{eqn:4h2p_linear_xi}
        \dfrac{1}{N_\vk} \sum_{\vq\in\mathcal{K}_\vq}\left(\sum_{kl}\braket{k(\vk_i-\vq), l(\vk_j+\vq) | i\vk_i, j\vk_j} t_{klab}(\vk_i-\vq, \vk_j+\vq, \vk_a)\right)
        - \xi t_{ijab}(\vk_i, \vk_j, \vk_a).
\end{equation}
The leading singularity of the integrand comes from the product with $(k,l) = (i,j)$.  
In this product, the ERI is singular at $\vq = \bm{0}$ with order $-2$ and the amplitude is singular at $\vq = \vk_i - \vk_a$ with order $0$. 
The correction in \cref{eqn:4h2p_linear_xi} defines a singularity subtraction with auxiliary function $h_\sigma(\vq)t_{ijab}(\vk_i,\vk_j,\vk_a)$. 
This auxiliary function shares exactly the same leading singularity as the integrand at $\vq = \bm{0}$ due to the ERI term, i.e.,  
\[
        \dfrac{4\pi}{|\Omega|}\dfrac{1}{|\vq|^2}t_{ijab}(\vk_i, \vk_j, \vk_a).
\]
Similar to the orbital energy analysis, the finite-size error in this 4h2p linear entry can be reduced to $\Or(N_\vk^{-1}$). 
However, the key difference here is that this error reduction is the case only for \textit{most but not all} the amplitude entries.
The exception is when the amplitude singularity $\vq = \vk_a - \vk_i$ is close or equal to the ERI singularity $\vq = \bm{0}$, i.e., when computing an ERI contraction entry 
whose momentum vector indices $\vk_i$ and $\vk_a$ are close or identical.  
In the worst case when $\vk_i = \vk_a$, the two singularities overlap and the finite-size error of such 4h2p linear entry can be shown 
to be still of scale $\Or(N_\vk^{-\frac13})$.
Similar analysis also applies to other terms in the ERI contraction calculation. 

To summarize, the Madelung constant correction does not uniformly reduce the finite-size errors in the ERI contraction tensor. 
This is the case for the 4h2p linear term and also for other terms in the ERI contraction. 
More precisely, we have the following technical error estimate (see \cref{lem:error_amplitude_contraction_app} in Appendix for the general statement and proof).
\begin{lem}[Error in ERI contractions]
                The finite-size error in the ERI contractions using exact CCD$(n)$ amplitude $t_n$ 
                with the Madelung constant correction in 
                \cref{eqn:madelung_ericontraction} satisfies 
                \begin{equation*}
                        \left|
                                \left[\mathcal{M}_{\mathcal{K}}\mathcal{A}_\text{TDL}(t_n) - \mathcal{A}_{N_\vk,\xi}(\mathcal{M}_{\mathcal{K}}t_n)\right]_{ijab}(\vk_i, \vk_j, \vk_a)
                        \right|
                        \leqslant 
                        C
                        \begin{cases}
                                \frac{1}{|\vq_{ia}|^2} N_\vk^{-1} & \vq_{ia} \neq \bm{0}\\
                                N_\vk^{-\frac13} & \vq_{ia} = \bm{0}\\
                        \end{cases},  
                \end{equation*}
                where $\vq_{ia} = \vk_ a - \vk_i + \vG_0$ with $\vG_0\in \mathbb{L}^*$ chosen such that $\vq_{ia}\in\Omega^*$. 
        \end{lem}

As a result of this nonuniform error reduction, the maximum entrywise finite-size error in the CCD$(n)$ amplitude calculation 
$\mathcal{M}_\mathcal{K}t_n - T_{n}^{N_\vk}$ with the Madelung constant correction remains $\Or(N_\vk^{-\frac13})$. However, the \textit{average} entrywise error satisfies the bound 
\[
\dfrac{1}{N_\vk^3}\sum_{ijab}\sum_{\vk_i,\vk_j,\vk_a\in\mathcal{K}}
\left|
        [\mathcal{M}_\mathcal{K}t_n - T_{n}^{N_\vk}]_{ijab}(\vk_i, \vk_j, \vk_a)
\right| 
\leqslant C N_\vk^{-1}.
\]
From the error decomposition in \cref{eqn:decomposition_energy}, such a refined bound is sufficient for our finite-size error analysis. Therefore 
the Madelung constant correction to the ERI contraction and the orbital energies can effectively reduce the finite-size error 
in the overall CCD$(n)$ energy calculation to $\Or(N_\vk^{-1})$.

\section{Numerical examples}\label{sec:numerics}

To validate the above theoretical analysis, we conduct CCD and CCD$(n)$ calculations on a 3D periodic system of hydrogen dimers. 
One hydrogen dimer is positioned at the center of each cubic unit cell with an edge length of $6$ Bohr in the $x$-direction, and a separation distance of $1.8$ Bohr.
For each uniform mesh $\mathcal{K}$, we perform an HF calculation  on $\mathcal{K}$ to obtain the orbitals and orbital energies.
We then perform CCD and CCD$(n)$ calculations under four distinct settings, each with a different combination of the Madelung constant correction to the orbital energies and the ERI contractions in the amplitude equation.
Specifically, we compute CCD and CCD$(n)$ with both corrections, the correction to the orbital energies alone,the correction to the ERI contractions alone, and without 
any corrections.
All the calculations are carried out using the PySCF package~\cite{SunBerkelbachBluntEtAl2018} with a minimal basis set \textsf{gth-szv}.

\begin{figure}[tbh]
        \centering
        \subfloat[CCD(1)]{
                \includegraphics*[width=0.23\textwidth]{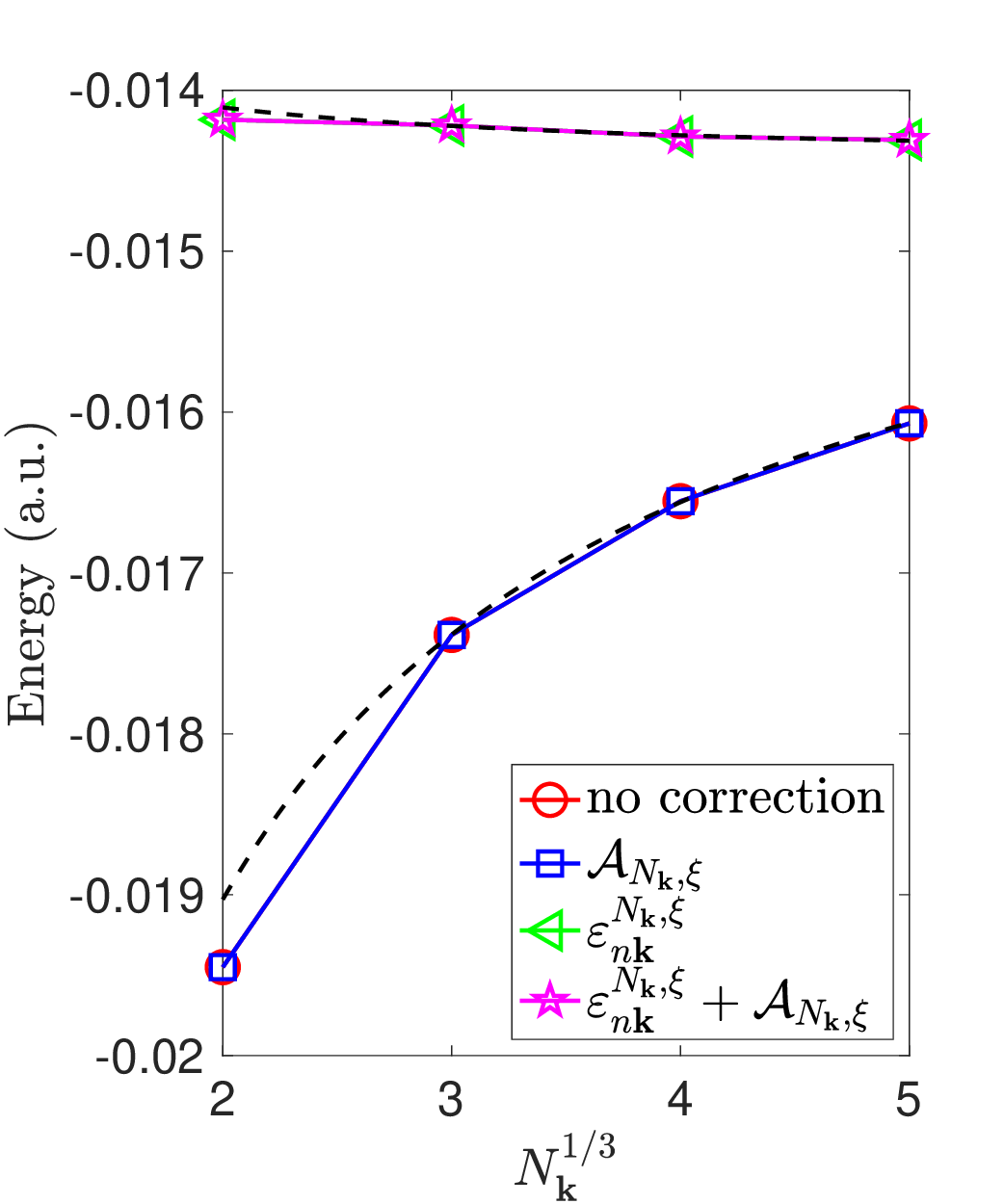}
        }
        \subfloat[CCD(2)]{
                \includegraphics*[width=0.23\textwidth]{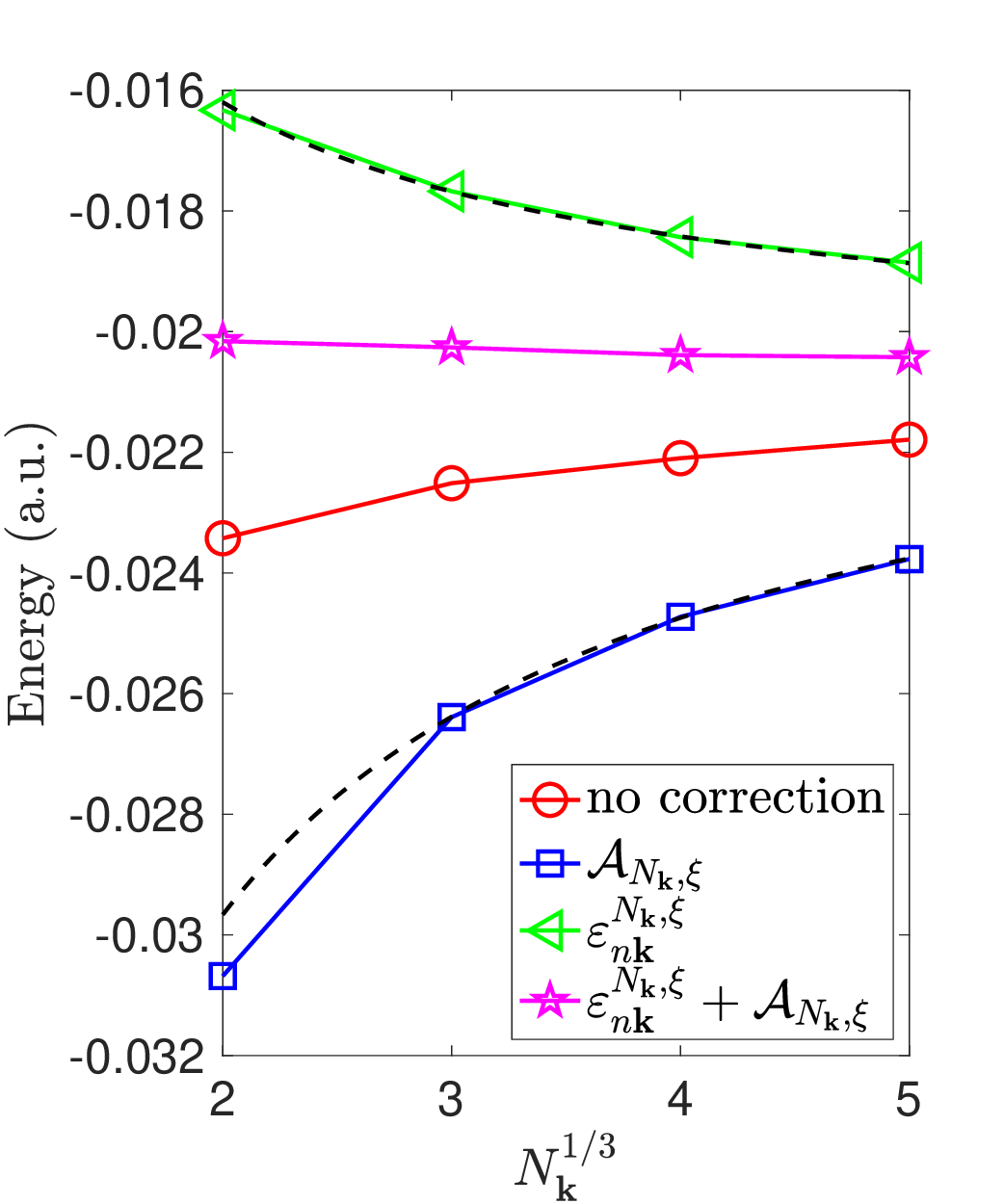}
        }
        \subfloat[CCD(3)]{
                \includegraphics*[width=0.23\textwidth]{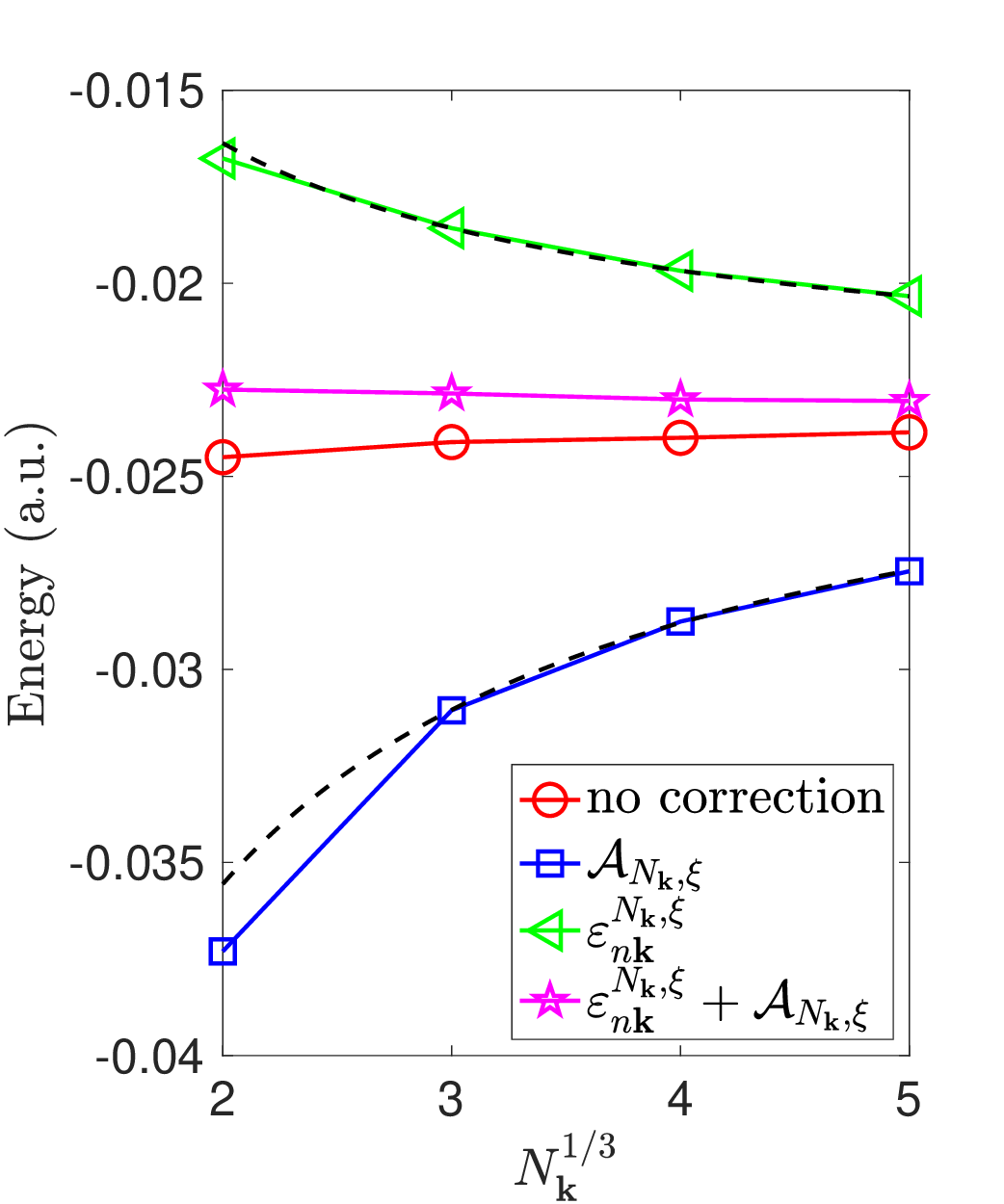}
        }
        \subfloat[CCD]{
                \includegraphics*[width=0.23\textwidth]{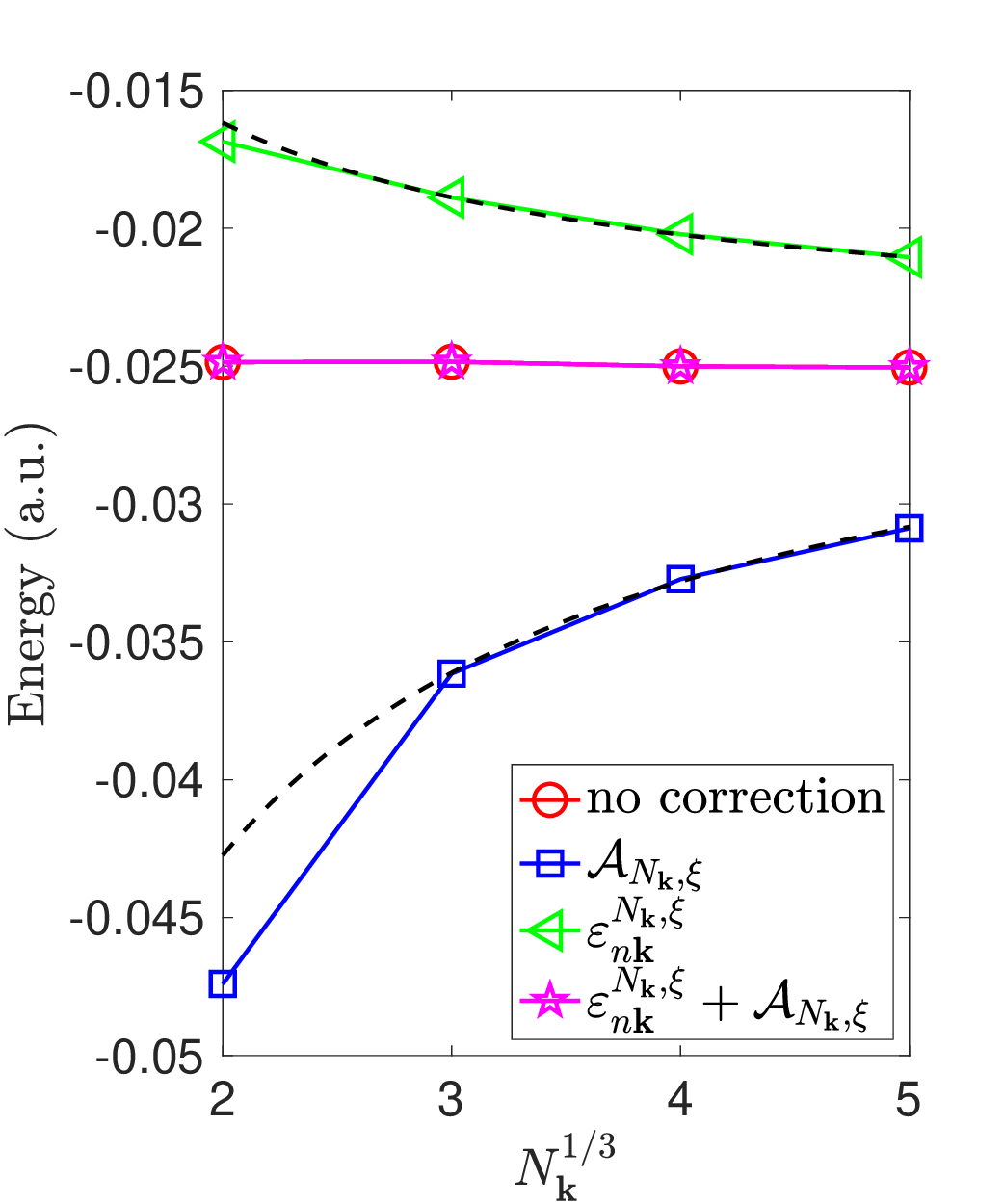}
        }

        \subfloat[CCD(1)]{
                \includegraphics*[width=0.23\textwidth]{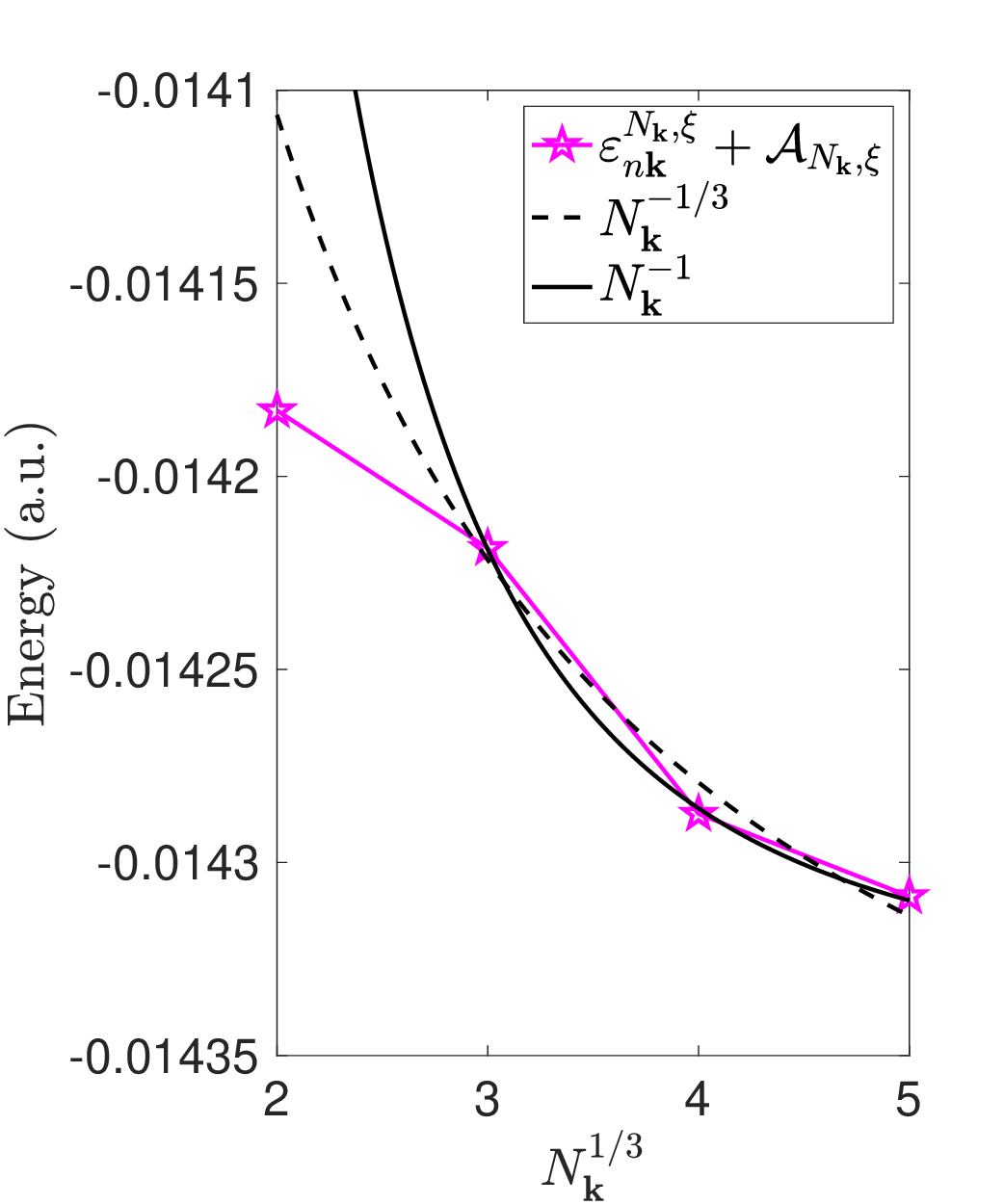}
        }
        \subfloat[CCD(2)]{
                \includegraphics*[width=0.23\textwidth]{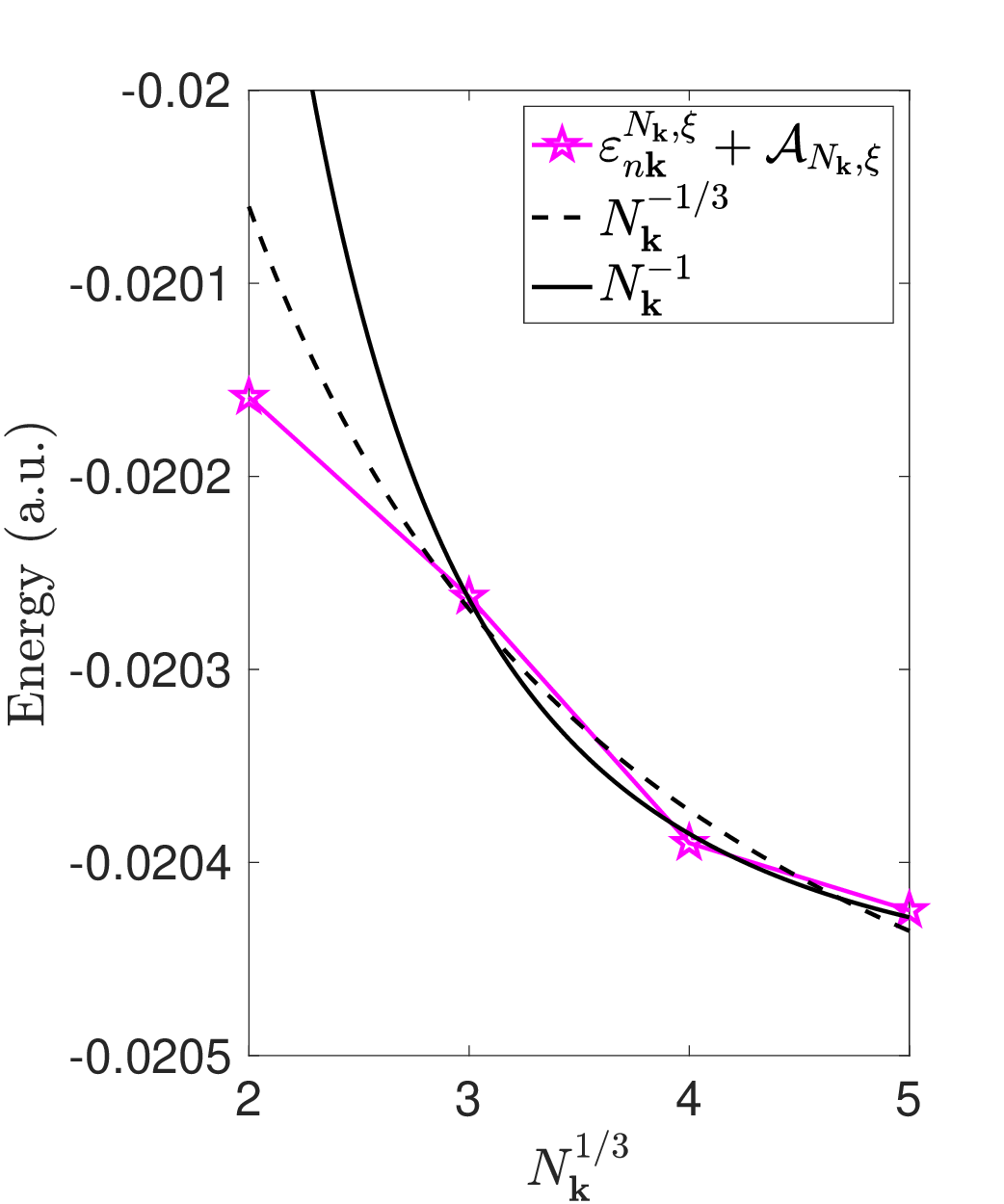}
        }
        \subfloat[CCD(3)]{
                \includegraphics*[width=0.23\textwidth]{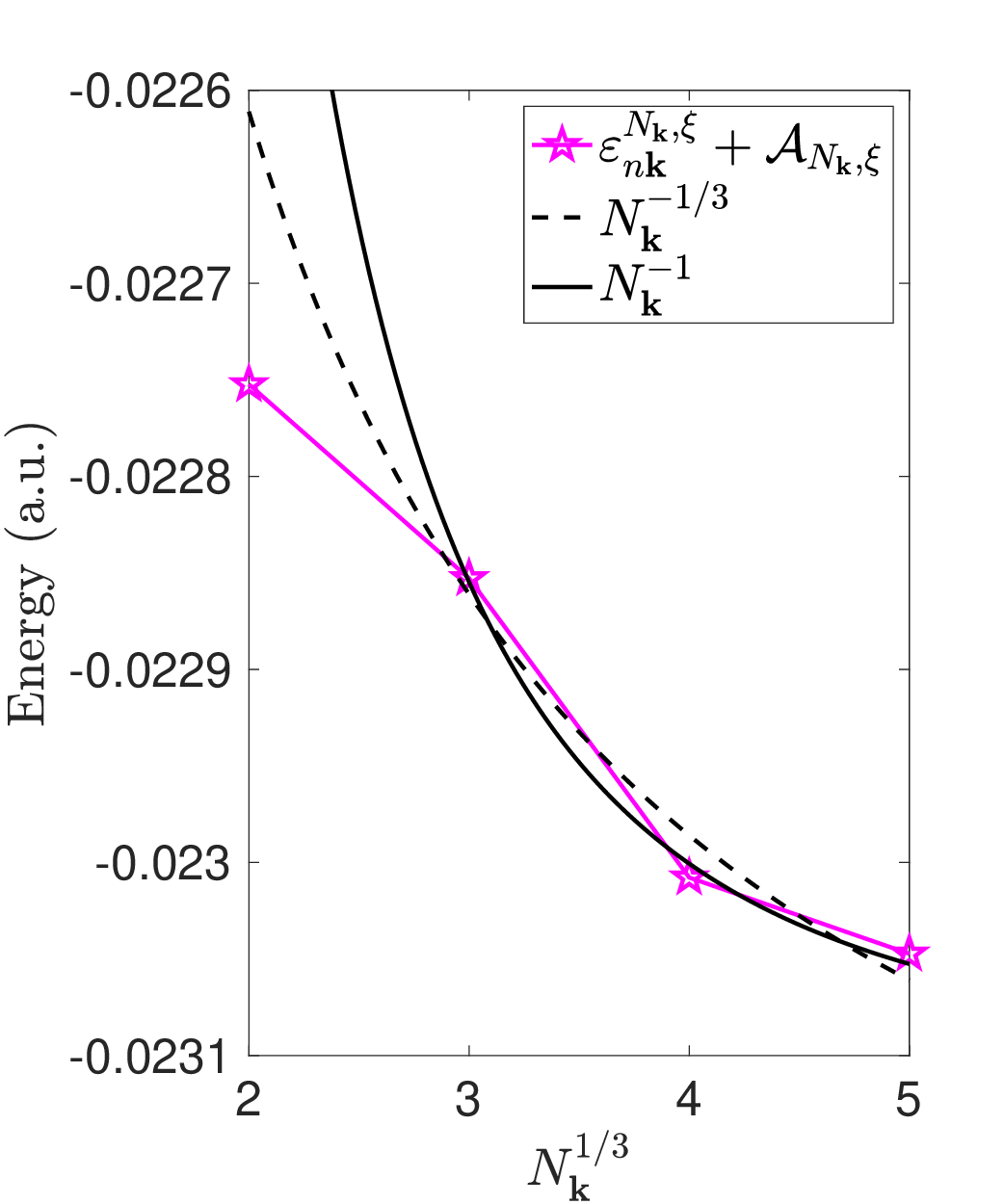}
        }
        \subfloat[CCD]{
                \includegraphics*[width=0.23\textwidth]{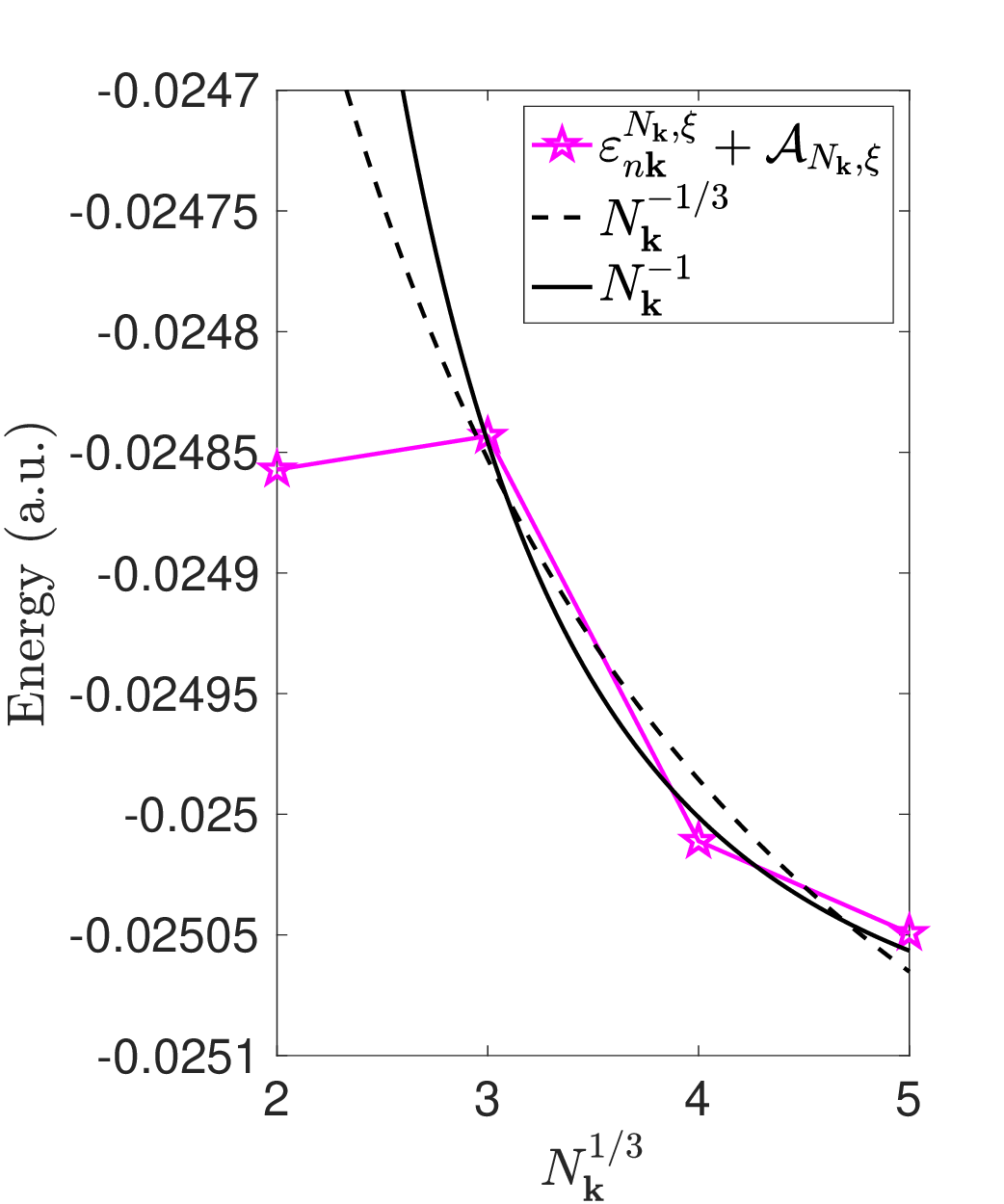}
        }
        \caption{
                Convergence of the CCD(1), CCD(2), CCD(3), and CCD correlation energies for a 3D periodic system of hydrogen dimers with increasing $N_\vk$. 
                The settings are distinguished by the presence or absence of the Madelung constant correction to 
                the orbital energies (``$\varepsilon_{n\vk}^{N_\vk,\xi}$'') and the ERI contraction (``$\mathcal{A}_{N_\vk,\xi}$'').
                In (a)-(d), the dashed curves show the power-law fitting using $C_0 + C_1N_\vk^{-\frac13}$ and data points $N_\vk^{\frac13} = 3,4,5$ for the two cases with partial Madelung constant correction.
                Subfigures (e)-(h) plot the curve fittings using $N_\vk^{-\frac13}$ and $N_\vk^{-1}$ over the calculations with 
                correction to both components, numerically corroborating the inverse volume scaling. 
        }
        \label{fig:ccd_h2}
\end{figure}

\Cref{fig:ccd_h2} illustrates the numerical results of the CCD(1), CCD(2), CCD(3), and converged CCD calculations.
For CCD(1), which is identical to MP2, we have $T_{0} = \bm{0}$ and the Madelung constant correction to the ERI contractions has no effect. 
As a result, the curves for the calculations with and without this correction are identical.
For CCD(2) and higher, the four correction settings produce distinct curves. 
Only the CCD($n$) calculation with both corrections exhibits convergence rate of $\Or(N_\vk^{-1})$, 
while the other three calculations have convergence rates of $\Or(N_\vk^{-\frac13})$. 
This highlights the importance of taking into account the Madelung constant correction to both the orbital energies and ERI contractions in higher-level CCD($n$) calculations.
As CCD$(n)$ converges to CCD, the difference between the calculation with both corrections and the one without any corrections gradually diminishes to zero, and the finite-size error satisfies the inverse volume scaling.

\section{Discussion}\label{sec:discussion}


Recent years have witnessed significant progresses in leveraging wavefunction methods to study solids. 
One important driver of this trend is the potential of these methods to provide systematically improvable results.
Among the simplest post-HF methods, MP2 has been applied to increasingly large periodic systems along with the development of 
efficient implementations utilizing density fitting and parallelization techniques~\cite{SchaeferRambergerKresse2017,BintrimBerkelbach2022}. 
Notably, the next order method, MP3, has also been applied to periodic systems recently for the first time~\cite{BanerjeeSokolov2022}.
Though computationally much more expensive, CC theory has also gained momentum in this context. 
It is being increasingly employed to compute ground state and band structures for a variety of solid materials including insulators and metallic systems 
\cite{McClainSunChanEtAl2017,ZhangGruneis2019,GaoSunJasonMottaEtc2020,NeufeldYeBerkelbach2022,NeufeldBerkelbach2023}. 
Furthermore, related theories such as equation-of-motion coupled-cluster theory, GW method, and algebraic diagrammatic construction theory 
also start their utilizations in computing excited state properties of solids \cite{WangBerkelbach2020,ZhuChan2021,BanerjeeSokolov2022}.
Due to the omnipresence of finite-size effects, as wavefunction methods continue to be increasingly applied to periodic systems, 
the importance of developing a thorough understanding of finite-size effects becomes even more pronounced. 

In this work, we fill a gap in understanding the finite-size error in periodic coupled cluster calculation for insulating systems. 
Notably, we reveal an unexpected inverse volume scaling of this error in CCD calculation. 
This behavior manifests even in the absence of any finite-size correction schemes, owing to an error cancellation.
Our findings, together with the methodologies employed in this study, provide valuable insights for practitioners, method developers, and theorists.

For practitioners, when applying computational quantum chemistry methods to periodic systems, reducing finite-size errors using techniques such as power-law extrapolation 
requires an in-depth understanding of the error scaling. 
This is particularly important when calculations are constrained to small-sized systems due to the steep increase of the computational cost with respect to the system size and limited resources. 
Many production-level quantum chemistry packages that support periodic systems use certain finite-size corrections to reduce errors in the HF exchange energy calculation and the HF orbital energies. 
For instance, the truncated Coulomb correction scheme \cite{SundararamanArias2013,spencer08} can be applied to insulating systems, 
so that the finite-size error in orbital energies decays super-algebraically with respect to $N_{\vk}$.  
However, our analysis shows that if this correction is only applied to the orbital energy but without any correction to ERI contractions, 
the finite-size error of the correlation energy deteriorates from inverse volume scaling to inverse length scaling. 
This is also a finding consistent with numerical observations.

For method developers, a key aspect of our result is the connection between the Madelung constant correction and the singularity subtraction method. 
This relationship serves not just as a crucial element in our theoretical proof, but also points towards new methods for further finite-size error reduction. 
The Madelung constant correction operates as a one-shot correction, while the singularity subtraction method can be systematically improved to reduce the finite-size error. 
Our analysis and correction schemes can also be extended to more advanced coupled cluster theories such as coupled cluster singles and doubles (CCSD), and coupled cluster singles, doubles and triples (CCSDT). 
Exploring the finite-size error in lower dimensional systems, such as magic angle twisted bilayer graphene (MATBG) that requires a different form of the Coulomb kernel~\cite{Ismail-Beigi2006,BultinckKhalafLiuEtAl2020}, can lead to different finite-size scaling patterns and novel correction schemes. Going beyond finite-size corrections, the singularity structure of the amplitude, as outlined \cref{lem:thm1_nonsmooth}, may be of independent interest. 
The singularity structure imposes important analytic constraints that should be taken into account when developing efficient numerical methods for compressing the CC amplitude tensor.

For theorists, several critical questions persist: How can finite-size error analysis be integrated with the study of complete basis limits to tackle basis set dependence? How should the finite-size error behavior in metallic systems be analyzed, especially when the orbital energy difference in the denominator could vanish? Can the scope of finite-size analysis be broadened to include more complex systems like disordered systems and finite-temperature alloys? These questions present a fertile ground for future research.%

\vspace{1em}
\noindent\textbf{Acknowledgement:}

This material is based upon work supported by the U.S. Department of Energy, Office of Science, Office of Advanced Scientific Computing Research and Office of Basic Energy Sciences, Scientific Discovery through Advanced Computing (SciDAC) program under Award Number DE--SC0022364 (X.X.). This work is also partially supported by the Applied Mathematics Program of the US Department of Energy (DOE) Office of Advanced Scientific Computing Research under contract number DE-AC02-05CH1123 (L.L.).
L.L. is a Simons Investigator in Mathematics.  We are grateful to Timothy Berkelbach, Garnet Chan, and Alexander Sokolov for
insightful discussions, \REV{and to the anonymous reviewers for their helpful suggestions, which have been instrumental in enhancing the presentation and clarity of our paper.}

\bibliography{mp2}

\newpage
\appendix
\renewcommand\thesubsection{\thesection.\arabic{subsection}}
\crefalias{section}{appendix}

\section{CCD amplitude equation}\label{appendix:ccd}

In this section, we use a capital letter to denote an index pair consisting of the orbital index and the $\vk$-index. 
For instance, $I=(i, \vk_i)$, $J=(j, \vk_j)$, $A=(a, \vk_a)$ etc. 
We use $P\in \{I, J, K, L\}$ to refer to occupied orbitals 
(a.k.a., holes) and $P\in\{A,B,C,D\}$ to refer to unoccupied orbitals (a.k.a., particles). 
Any summation $\sum_P$ refers to summing over all occupied or virtual orbital indices $p$ and all momentum vectors $\vk_p\in \mathcal{K}$ while 
the crystal momentum conservation is enforced according to the summand.  
This notation is used only in this section to simplify the notation and also to connect the equations to those in the molecular case \cite{ShavittBartlett2009} for better readability.

Using a finite mesh $\mathcal{K}$ of size $N_\vk$, the normalized CCD amplitude $T_*^{N_\vk} = \{T_{ijab}(\vk_i, \vk_j, \vk_a)\}  := \{t_{IJ}^{AB}\}$ 
with $\vk_i,\vk_j,\vk_a \in \mathcal{K}$ is defined as the solution of the amplitude equation
\begin{align}
        t_{IJ}^{AB}
        & 
        = 
        \frac{1}{\varepsilon_{IJAB}^{N_\vk}} 
        [\mathcal{A}_{N_\vk}(T_*^{N_\vk})]_{IJAB}
        \nonumber\\
        & = 
        \frac{1}{\varepsilon_{IJAB}^{N_\vk}} 
        \Bigg[
        \braket{AB|IJ} + 
        \mathcal{P}
        \left(
        \sum_C\kappa_C^A t_{IJ}^{CB}
        -\sum_K\kappa_I^K t_{KJ}^{AB} 
        \right)
        \nonumber       
        \\
        & \hspace*{5em}
        + \dfrac{1}{N_\vk}\sum_{KL}\chi_{IJ}^{KL} t_{KL}^{AB} 
        + \dfrac{1}{N_\vk}\sum_{CD}\chi_{CD}^{AB} t_{IJ}^{CD} 
        \nonumber       
        \\
        & \hspace*{5em}
        + \mathcal{P} 
        \left(\dfrac{1}{N_\vk}\sum_{KC} (2\chi_{IC}^{AK} - \chi_{CI}^{AK}) t_{KJ}^{CB} - \chi_{IC}^{AK}t_{KJ}^{BC} - \chi_{CJ}^{AK}t_{KI}^{BC}
        \right)
        \Bigg],\qquad \forall I, J, A, B,
        \label{eqn:amplitude_ccd}
\end{align}
where $\varepsilon_{IJAB}^{N_\vk} = \varepsilon_{i\vk_i}^{N_\vk} + \varepsilon_{j\vk_j}^{N_\vk} - \varepsilon_{a\vk_a}^{N_\vk} - \varepsilon_{b\vk_b}^{N_\vk}$, $[\mathcal{A}_{N_\vk}(T)]_{IJAB} = [\mathcal{A}_{N_\vk}(T)]_{ijab}(\vk_i,\vk_j,\vk_a)$, 
and $\mathcal{P}$ is a permutation operator defined as
$
\mathcal{P}(\cdots)_{IJ}^{AB} = (\cdots)_{IJ}^{AB} + (\cdots)_{JI}^{BA}.
$         
This reformulation of CCD amplitude equation is derived from the CCSD amplitude equation in \cite{HirataPodeszwaTobitaEtAl2004} by removing all the terms related to single amplitudes and normalizing  
the involved ERIs and amplitudes (which gives the extra $1/N_\vk$ factor in the equation and the intermediate blocks). 
The intermediate blocks in the equation are defined as
\begin{align*}
        \kappa_C^A & = -\dfrac{1}{N_\vk^2}\sum_{KLD} \left( 2\braket{KL|CD} - \braket{KL|DC} \right) t_{KL}^{AD},
        \\
        \kappa_I^K & = \dfrac{1}{N_\vk^2}\sum_{LCD} \left( 2\braket{KL|CD} - \braket{KL|DC} \right) t_{IL}^{CD},
        \\
        \chi_{IJ}^{KL} & = \braket{KL|IJ} + \dfrac{1}{N_\vk}\sum_{CD} \braket{KL|CD} t_{IJ}^{CD},
        \\
        \chi_{CD}^{AB} & = \braket{AB|CD},
        \\
        \chi_{IC}^{AK} & = \braket{AK|IC} + \dfrac{1}{2N_\vk}\sum_{LD} \left(2\braket{LK|DC} - \braket{LK|CD}\right) t_{IL}^{AD} - \braket{LK|DC} t_{IL}^{DA},
        \\
        \chi_{CI}^{AK} &= \braket{AK|CI} - \dfrac{1}{2N_\vk}\sum_{LD} \braket{LK|CD} t_{IL}^{DA},
\end{align*}
and their momentum vector indices also assume the crystal momentum conservation
\begin{align*}
        \kappa_{P}^{Q} & \rightarrow\vk_p - \vk_q \in \mathbb{L}^*,\\
        \chi_{PQ}^{RS} & \rightarrow \vk_p + \vk_q - \vk_r - \vk_s \in \mathbb{L}^*. 
\end{align*}

In the TDL, the amplitude equation for the exact amplitude $t_* = \{t_{ijab}(\vk_i,\vk_j,\vk_a)\} := \{t_{IJ}^{AB}\}$ 
as functions of $\vk_i, \vk_j, \vk_a \in \Omega^*$ can be formulated by letting $\mathcal{K}$ in \cref{eqn:amplitude_ccd} converge to $\Omega^*$ as 
\begin{align}
        t_{IJ}^{AB}
        &= 
        \dfrac{1}{\varepsilon_{IJAB}^\text{TDL}} [\mathcal{A}_{\text{TDL}}(t_*)]_{IJAB}
        \nonumber \\
        & = 
        \frac{1}{\varepsilon_{IJAB}^\text{TDL}}
        \Bigg[
        \braket{AB|IJ} + 
        \mathcal{P}
        \left(
        \sum_C\kappa_C^A t_{IJ}^{CB}
        -\sum_K\kappa_I^K t_{KJ}^{AB} 
        \right)
        \nonumber
        \\
        & \hspace*{5em}
        + 
        \dfrac{1}{|\Omega^*|}\int_{\Omega^*}\ud\vk_k \sum_{kl}\chi_{IJ}^{KL} t_{KL}^{AB}
        + 
        \dfrac{1}{|\Omega^*|}\int_{\Omega^*}\ud\vk_c \sum_{cd}\chi_{CD}^{AB} t_{IJ}^{CD} 
        \nonumber
        \\
        & \hspace*{5em}
        + \mathcal{P} 
        \left(
                        \dfrac{1}{|\Omega^*|}\int_{\Omega^*}\ud\vk_k \sum_{kc} (2\chi_{IC}^{AK} - \chi_{CI}^{AK}) t_{KJ}^{CB} - \chi_{IC}^{AK}t_{KJ}^{BC} - 
                        \chi_{CJ}^{AK}t_{KI}^{BC}
        \right) 
        \Bigg], \qquad\forall I, J, A, B,
        \label{eqn:amplitude_ccd_tdl}
\end{align} 
where the intermediate blocks in the TDL are defined as 
\begin{align*}
        \kappa_C^A & = -\dfrac{1}{|\Omega^*|^2}\int_{\Omega^*\times \Omega^*}\ud\vk_k\ud\vk_l \sum_{kld} \left( 2\braket{KL|CD} - \braket{KL|DC} \right) t_{KL}^{AD},
        \\
        \kappa_I^K & = \dfrac{1}{|\Omega^*|^2}\int_{\Omega^*\times \Omega^*}\ud\vk_c\ud\vk_d  \sum_{lcd} \left( 2\braket{KL|CD} - \braket{KL|DC} \right) t_{IL}^{CD},
        \\
        \chi_{IJ}^{KL} & = \braket{KL|IJ} + \dfrac{1}{|\Omega^*|}\int_{\Omega^*}\ud\vk_c \sum_{cd} \braket{KL|CD} t_{IJ}^{CD},
        \\
        \chi_{CD}^{AB} & = \braket{AB|CD},
        \\
        \chi_{IC}^{AK} & = \braket{AK|IC} + \dfrac{1}{2|\Omega^*|}\int_{\Omega^*}\ud\vk_l \sum_{ld}  \left(2\braket{LK|DC} - \braket{LK|CD}\right) t_{IL}^{AD} - \braket{LK|DC} t_{IL}^{DA} ,
        \\
        \chi_{CI}^{AK} &= \braket{AK|CI} - \dfrac{1}{2|\Omega^*|}\int_{\Omega^*}\ud\vk_l \sum_{ld} \braket{LK|CD} t_{IL}^{DA}. 
\end{align*}

\section{Proof of Theorem \ref{thm:error_ccdn} }
\label{app:error_ccdn}

\newtheorem*{restateddef2}{Restatement of \cref{thm:error_ccdn}}
\begin{restateddef2}
        In CCD($n$) calculation, the finite-size error in the correlation energy scales as $\Or(N_\vk^{-\frac13})$ in each of the following scenarios 
        (1) there is no finite-size correction, (2) the Madelung constant correction is only applied to the ERI contraction $\mathcal{A}_{N_\vk}$, 
        and (3) the Madelung constant correction is only applied to the orbital energy $\varepsilon_{n\vk}^{N_\vk}$. 

        When the Madelung constant correction is applied to both $\mathcal{A}_{N_\vk}$ and $\varepsilon_{n\vk}^{N_\vk}$ in the CCD$(n)$ calculation, 
        the overall finite-size error scales as $\Or(N_\vk^{-1})$.
        
        As a special case, the same conclusion applies to MP3 calculations.
\end{restateddef2}

\subsection{Proof Outline}
\label{subsec:outline}

The main context of this paper has provided a brief description of the main proof idea. 
In this proof, we will recap some of the equations and concepts discussed before to make it self-contained.
The proofs of \cref{thm:error_ccdn} for CCD$(n)$ calculations with various types of corrections are 
based on the error splitting in \cref{eqn:decomposition_energy} 
and \cref{eqn:decomposition_amplitude}, i.e., 
\begin{align}
        \left|
                E_{\text{CCD}(n)}^\text{TDL} - E_{\text{CCD}(n)}^{N_\vk}
        \right|
        & \leqslant 
        \left|
                \mathcal{G}_\text{TDL}(t_n) - \mathcal{G}_{N_\vk}(\mathcal{M}_{\mathcal{K}}t_n)
        \right|
        +
        C
        \left\|
                \mathcal{M}_\mathcal{K}t_n - T_{n}^{N_\vk}
        \right\|_1,
        \label{eqn:decomposition_energy_app}
        \\
        \mathcal{M}_{\mathcal{K}}t_{n} - T_{n}^{N_\vk}
        & = \dfrac{1}{\varepsilon^\text{TDL}}
        \left[
        \mathcal{M}_{\mathcal{K}}
        \mathcal{A}_\text{TDL}(t_{n-1})
        -
        \mathcal{A}_{N_\vk}(\mathcal{M}_{\mathcal{K}}t_{n-1})
        \right]
        \nonumber\\
        & \hspace*{1em}
        + 
        \mathcal{A}_{N_\vk}(\mathcal{M}_{\mathcal{K}}t_{n-1})
        \left[
        \dfrac{1}{\varepsilon^\text{TDL}}
        -
        \dfrac{1}{\varepsilon^{N_\vk}}
        \right]
        \nonumber\\
        & \hspace*{1em}
        +
        \dfrac{1}{\varepsilon^{N_\vk}}
        \left[
        \mathcal{A}_{N_\vk}( \mathcal{M}_{\mathcal{K}}t_{n-1})
        -
        \mathcal{A}_{N_\vk}(T_{n-1}^{N_\vk})
        \right].
        \label{eqn:decomposition_amplitude_app}
\end{align} 
Note that the amplitude error here is measured in the average norm as 
\[
        \|T\|_1 = \dfrac{1}{N_\vk^3}\sum_{\vk_i,\vk_j,\vk_a\in\mathcal{K}}\sum_{ijab} \left| T_{ijab}(\vk_i,\vk_j,\vk_a) \right|. 
\]
The finite-size error in CCD$(n)$ calculation is thus decomposed into the error in energy calculation 
using exact amplitude, the error in ERI contractions, the error in orbital energies, and the error accumulated from previous iteration.  
The latter three errors together make up of the error in amplitude calculation. 
The error in energy calculation using exact amplitude is studied previously in Ref.~\cite{XingLin2023}. 
For completeness, we provide a brief review of the main results below. 

\subsubsection*{Brief review of error in energy calculation with exact amplitude}
One basic observation is that this error in CCD$(n)$ calculation can be interpreted as the quadrature error of a specific
trapezoidal rule as
\begin{equation}\label{eqn:error_energy_app}
        \left|
                \mathcal{G}_\text{TDL}(t_n) - \mathcal{G}_{N_\vk}(\mathcal{M}_{\mathcal{K}}t_n)
        \right| 
        = 
        \left|
                \dfrac{1}{|\Omega^*|^3}
                \mathcal{E}_{\Omega^*\times\Omega^*\times\Omega^*}
                \left(
                \sum_{ijab}
                (W_{ijab} [t_n]_{ijab}), 
                \mathcal{K}\times\mathcal{K}\times\mathcal{K}
                \right)
        \right|.
\end{equation}
Since both $W_{ijab}$ and $[t_n]_{ijab}$ are periodic with respect to $\vk_i,\vk_j,\vk_a\in\Omega^*$, the asymptotic
scaling of this quadarture error depends on the smoothness of these two components that constitute the integrand. 

A general ERI $\braket{n_1\vk_1, n_2\vk_2 | n_3\vk_3, n_4\vk_4}$ can be viewed as a function of momentum vectors $\vk_1,\vk_2$ and its momentum transfer vector $\vq = \vk_3 - \vk_1$. 
This function is periodic with respect to each variable in $\Omega^*$ and its definition in \cref{eqn:eri} can be decomposed as 
\begin{equation}\label{eqn:eri_decomp}
        \frac{4\pi}{\abs{\Omega}} 
        \frac{
        \hat{\varrho}_{n_1\vk_1,n_3(\vk_1+\vq)}(\bm{0}) \hat{\varrho}_{n_2\vk_2,n_4(\vk_2 - \vq)}(\bm{0})
        }{
                |\vq|^2
        }
        + 
        \frac{4\pi}{\abs{\Omega}}
        \sum_{\vG\in\mathbb{L}^*\setminus \{\bm{0}\}}
        \dfrac{\cdots}{|\vq+\vG|^2}.
\end{equation}
The numerators of all the fractions are smooth with respect to $\vk_1, \vk_2, \vq$ (note the assumption that $\psi_{n\vk}(\vrr)$ is smooth with respect to $\vk$).
Therefore, this ERI example is smooth with respect to $\vk_1, \vk_2$ and has one 
point singularity in $\Omega^*$ with respect to $\vq$ at $\vq = \bm{0}$, which is due to the first fraction term.
For any fixed $\vk_i, \vk_j$, this point singularity can be characterized using the concept of \textit{algebraic singularity of certain orders}.

        \begin{defn}[Algebraic singularity for univariate functions]\label{def:fractional_app}
        A function $f(\vx)$ has \textbf{algebraic singularity of order} $\gamma\in\RR$ at $\vx_0 \in \mathbb{R}^d$  if
        there exists $\delta > 0$ such that 
        \[
        \left|
        \dfrac{\partial^\valpha}{\partial \vx^{\valpha}} f(\vx)  
        \right|
        \leqslant 
        C_{\valpha,\delta} |\vx - \vx_0|^{\gamma - |\valpha|},
        \qquad \forall 0 < |\vx - \vx_0| < \delta, \ \forall \valpha \geqslant 0,
        \]
        where the constant $C_{\valpha,\delta}$ depends on $\delta$ and the non-negative $d$-dimensional derivative multi-index $\valpha$.  
        For brevity, $f$ is also said to be singular at $\vx_0$ of order $\gamma$.
        \end{defn}

The numerator of the first fraction in \cref{eqn:eri_decomp} scales as $\Or(|\vq|^s)$ with $s \in \{0, 1, 2\}$ near $\vq = \bm{0}$ using the orbital orthogonality.  
The value of $s$ depends on the relation between orbital indices $(n_1, n_2)$ and $(n_3, n_4)$. 
As a result, the algebraic singularity of the ERI example above at $\vq = \bm{0}$ with any fixed $\vk_1, \vk_2$ has order in $\{-2, -1, 0\}$. 
In addition, to connect this singularity with varying $\vk_1, \vk_2 \in \Omega^*$, we introduce the algebraic singularity with respect 
to one variable for a multivariate function. 
\begin{defn}[Algebraic singularity for multivariate functions]\label{def:fractional2}
        A function $f(\vx, \vy)$ is smooth with respect to $\vy \in V_Y\subset \mathbb{R}^{d_y}$ for any fixed $\vx$ and has algebraic singularity of order $\gamma$ with respect to $\vx$ at $\vx_0 \in \mathbb{R}^{d_x}$ 
        if there exists $\delta > 0$ such that 
        \[
                \left|
                \dfrac{\partial^\valpha}{\partial \vx^{\valpha}} \left( \dfrac{\partial^\vbeta}{\partial\vy^\vbeta}f(\vx, \vy)  \right)
                \right|
                \leqslant 
                C_{\valpha,\vbeta,\delta} |\vx - \vx_0|^{\gamma - |\valpha|},  \quad \forall 0<|\vx - \vx_0|< \delta, \forall \vy \in V_Y, \forall \valpha,\vbeta\geqslant 0,
        \]
        where constant $C_{\valpha,\vbeta,\delta}$ depends on $\delta$, $\valpha$ and $\vbeta$.
        Compared to the univariate case in \cref{def:fractional}, the key additions are the shared algebraic singularity of partial derivatives over $\vy$ at $\vx = \vx_0$ 
        of order $\gamma$ and the independence of $C_{\valpha, \vbeta,\delta}$ on $\vy\in V_Y$. 
\end{defn}

By this definition, the ERI example is smooth everywhere with respect to $\vk_i, \vk_j, \vq\in\Omega^*$ except at $\vq = \bm{0}$ of order $\gamma \in \{-2, -1, 0\}$.
Specifically, the order $\gamma$ equals to $-2, -1$ and $0$, respectively, when the orbital indices are fully matched, i.e., $n_1=n_3,n_2=n_4$, partially matched, i.e., $n_1=n_3,n_2\neq n_4$ or $n_1\neq n_3, n_2=n_4$, and not matched, i.e., $n_1\neq n_3, n_2\neq n_4$.
If treating the ERI example as a function of $\vk_1, \vk_2, \vk_3$ instead, we equivalently claim that the function is singular at $\vk_1 = \vk_3$ of order $\gamma$. 

One key result in \cite{XingLin2023} is the singularity structure characterization for the exact CCD$(n)$ amplitude $t_n$, which is essential for estimating  
the quadrature error in \cref{eqn:error_energy_app}. 
It turns out that that each exact amplitude entry $[t_n]_{ijab}(\vk_i,\vk_j,\vk_a)$ with any $n > 0$ has one point
of algebraic singularity of order $0$ at $\vk_a - \vk_i = \bm{0}$, sharing a similar singularity structure as the ERI $\braket{i\vk_i,j\vk_j | a\vk_a, b\vk_b}$
or the exact MP2/CCD($1$) amplitude entry $(\varepsilon_{i\vk_i,j\vk_j,a\vk_a,b\vk_b}^\text{TDL})^{-1}\braket{a\vk_a, b\vk_b | i\vk_i,j\vk_j}$.

\newtheorem*{restatedlem0}{Restatement of \cref{lem:thm1_nonsmooth}}
\begin{restatedlem0}[Singularity structure of the amplitude, Lemma 4 in \cite{XingLin2023}]
        In CCD$(n)$ calculation with $n > 0$, each entry of the exact double amplitude $t_n$ belongs to the following function space
        \begin{align*}
        \mathbb{T}(\Omega^*) = \big\{
                f(\vk_i, \vk_j, \vk_a):\  
                & f \text{ is periodic with respect to } \vk_i, \vk_j, \vk_a \in \Omega^*,\\
                & f \text{ is smooth everywhere except at } \vk_a = \vk_i 
                \text{ with algebraic singularity of order } 0,\\
                & f \text{ is smooth with respect to } \vk_i, \vk_j \text{ at the singularity }\vk_a = \vk_i
        \big\}.
        \end{align*}
\end{restatedlem0}

Based on the above singularity structures of ERIs and exact amplitudes, the integrand in the energy calculation in \cref{eqn:error_energy_app} 
consists of products of periodic functions where each has one point singularity of order $0$. 
(Recall that $W_{ijab}$ is the antisymmetrized ERI and consists of two ERIs that can be treated separately.)
We provide a sharp quadrature error bound for trapezoidal rules over periodic functions in such a product form, and 
its application to \cref{eqn:error_energy_app} gives 

\begin{lem}[Energy error with exact amplitude, Lemma 5 in \cite{XingLin2023}]
        \label{lem:error_energy}
        In CCD$(n)$ calculation with any $n > 0$, 
        the finite-size error in the energy calculation using exact amplitude $t_n$ can be estimated as  
        \[
        \left|
                \mathcal{G}_\text{TDL}(t_n) - \mathcal{G}_{N_\vk}(\mathcal{M}_\mathcal{K}t_n)
        \right|
        \leqslant C N_\vk^{-1},
        \]
        where constant $C$ depends on $t_n$. 
\end{lem}

\subsubsection*{Error in amplitude calculation}
Based on the error splitting in \cref{eqn:decomposition_amplitude_app}, analysis of the error in amplitude calculation is reduced to estimating the 
two main error terms $\mathcal{M}_{\mathcal{K}}\mathcal{A}_\text{TDL}(t)-\mathcal{A}_{N_\vk}(\mathcal{M}_{\mathcal{K}}t)$ and $\varepsilon_{n\vk}^\text{TDL} -\varepsilon_{n\vk}^{N_\vk}$, 
and understanding how $\mathcal{A}_{N_\vk}$ amplifies the amplitude error $\mathcal{M}_{\mathcal{K}}t_{n-1} - T_{n-1}^{N_\vk}$ from the previous iteration.

Like the energy calculation, the errors in ERI contractions and orbital energies consist of specific quadrature errors.
In both cases, we can show that the Madelung constant correction is connected to certain singularity subtraction methods and can significantly reduce the corresponding dominant 
quadrature errors.   
\begin{lem}[Error in ERI contractions]
\label{lem:error_amplitude_contraction_app}
        In CCD$(n)$ calculation with any $n > 0$, the finite-size error in the ERI contractions using exact amplitude $t_n$ without any corrections has its entries bounded as  
        \begin{equation}\label{eqn:error_eri_contraction0}
                \left|
                        \left[\mathcal{M}_{\mathcal{K}}\mathcal{A}_\text{TDL}(t_n) - \mathcal{A}_{N_\vk}(\mathcal{M}_{\mathcal{K}}t_n)\right]_{ijab}(\vk_i, \vk_j, \vk_a)
                \right|
                \leqslant 
                CN_\vk^{-\frac13},
        \end{equation}
        The Madelung constant correction reduces this error as 
        \begin{equation}\label{eqn:error_eri_contraction}
                \left|
                        \left[\mathcal{M}_{\mathcal{K}}\mathcal{A}_\text{TDL}(t_n) - \mathcal{A}_{N_\vk,\xi}(\mathcal{M}_{\mathcal{K}}t_n)\right]_{ijab}(\vk_i, \vk_j, \vk_a)
                \right|
                \leqslant 
                C
                \begin{cases}
                        \frac{1}{|\vq|^2} N_\vk^{-1} & \vq \neq \bm{0}\\
                        N_\vk^{-\frac13} & \vq = \bm{0}\\
                \end{cases},  
        \end{equation}
        where $\vq = \vk_ a - \vk_i + \vG_0$ with $\vG_0\in \mathbb{L}^*$ chosen such that $\vq \in\Omega^*$. 
        In both cases, constant $C$ depends on $t_n$ but not on the entry index $(i,j,a,b)$ and $(\vk_i,\vk_j,\vk_a)\in\mathcal{K}\times\mathcal{K}\times\mathcal{K}$. 
\end{lem}
\begin{rem}
        The prefactor $1/|\vq|^2$ in the above estimate is important when $\vq\in\mathcal{K}_\vq$ is 
        near the origin. 
        For example, if $\vk_i, \vk_a\in\mathcal{K}$ are adjacent to each other, $|\vq|$ is of scale $\Or(N_\vk^{-\frac13})$ 
        and the estimate in \cref{lem:error_amplitude_contraction_app} suggests an error bound of scale $\Or(N_\vk^{-\frac13})$.  
\end{rem}

\begin{lem}[Error in orbital energies]\label{lem:error_orbital_energy}
        The finite-size error in orbital energies without any corrections is bounded as 
        \begin{equation*}
                \left|
                        \varepsilon_{n\vk}^\text{TDL} - \varepsilon_{n\vk}^{N_\vk}
                \right|
                \leqslant 
                C
                \begin{cases}
                        N_\vk^{-\frac13} & \text{$n$ is occupied} \\
                        N_\vk^{-1} & \text{$n$ is virtual} 
                \end{cases}.
        \end{equation*}
        The Madelung constant correction reduces this error as 
        \begin{equation*}
                \left|
                        \varepsilon_{n\vk}^\text{TDL} - \varepsilon_{n\vk}^{N_\vk,\xi}
                \right|
                \leqslant 
                C N_\vk^{-1}.
        \end{equation*}
        In both cases, constant $C$ is independent of the entry index $n, \vk\in\mathcal{K}$. 
\end{lem}

Note that there are three multipliers for the three error terms in the amplitude error splitting \cref{eqn:decomposition_amplitude_app}.
These prefactors are bounded by constants independent of $N_\vk$. 
First, the orbital energy difference in \cref{eqn:decomposition_amplitude_app} satisfies $|\varepsilon_{i\vk_i,j\vk_j,a\vk_a,b\vk_b}^\text{TDL}| \geqslant  2\varepsilon_g$ by the assumed indirect gap 
$\varepsilon^\text{TDL}_{a\vk_a} - \varepsilon^\text{TDL}_{i\vk_i}\geqslant \varepsilon_g > 0$.
Second, since $\mathcal{A}_{N_\vk}(T)$ consists of constant, linear and quadratic terms of $T$, a straightforward estimate shows (e.g., using Lemma 7 in \cite{XingLin2023}) that 
\[
                \max_{ijab,\vk_i,\vk_j,\vk_a\in\mathcal{K}}
                \left|
                        [\mathcal{A}_{N_\vk}(\mathcal{M}_\mathcal{K}t_n)]_{ijab}(\vk_i,\vk_j,\vk_a)
                \right|
                \leqslant 
                C \max_{ijab}\|[t_n]_{ijab}\|^2_{L^{\infty}(\Omega^*\times\Omega^*\times\Omega^*)},
\]
where the $L^\infty$-norm of $[t_n]_{ijab}$ is finite according to \cref{lem:thm1_nonsmooth}. 

Based on the above estimates of these multipliers, we have that the summation of the first two error terms in \cref{eqn:decomposition_amplitude_app} is dominated by the error in the ERI contractions 
and the orbital energies discussed in \cref{lem:error_amplitude_contraction_app} and \cref{lem:error_orbital_energy}.
If applying the Madelung constant correction to both orbital energies and ERI contractions, the summation of the first two error terms 
has its entries bounded asymptotically the same as \cref{eqn:error_eri_contraction}.
Otherwise, its entries are bounded asymptotically as in \cref{eqn:error_eri_contraction0}. 
Lastly, for the third error term in \cref{eqn:decomposition_amplitude}, the application of $\mathcal{A}_{N_\vk}$/$\mathcal{A}_{N_\vk, \xi}$ can be proved to maintain the entrywise error scaling obtained in \cref{lem:error_amplitude_contraction_app}.
\begin{lem}\label{lem:error_amplitude_accumulation}
        Consider two arbitrary bounded amplitude tensors $T, S \in 
        \mathbb{C}^{n_\text{occ}\times n_\text{occ} \times n_\text{vir} \times n_\text{vir}\times N_\vk\times N_\vk \times N_\vk}$ and assume their entrywise upperbound independent of $N_\vk$. 
        If the difference between $T$ and $S$ satisfies an entry bound like \cref{eqn:error_eri_contraction0}, i.e., 
        \begin{equation*}                
                \left|
                        [T - S]_{ijab}(\vk_i, \vk_j, \vk_a)
                \right|
                \leqslant 
                C
                N_\vk^{-\frac13}.
        \end{equation*}
        The ERI-contraction map without any corrections $\mathcal{A}_{N_\vk}$ satisfies 
        \begin{equation*}                
                \left|
                        [\mathcal{A}_{N_\vk}(T) - \mathcal{A}_{N_\vk}(S)]_{ijab}(\vk_i, \vk_j, \vk_a)
                \right|
                \leqslant 
                C
                N_\vk^{-\frac13}.
        \end{equation*}

        If the difference between $T$ and $S$ satisfies an estimate similar to that in \cref{eqn:error_eri_contraction}, i.e., 
        \begin{equation*}                
                \left|
                        [T - S]_{ijab}(\vk_i, \vk_j, \vk_a)
                \right|
                \leqslant 
                C
                \begin{cases}
                        \frac{1}{|\vq|^2} N_\vk^{-1} & \vq \neq \bm{0} \\
                        N_\vk^{-\frac13} & \vq = \bm{0}\\
                \end{cases},
        \end{equation*}
        The ERI-contraction map with Madelung constant correction $\mathcal{A}_{N_\vk,\xi}$ satisfies 
        \begin{equation*}                
                \left|
                        [\mathcal{A}_{N_\vk, \xi}(T) - \mathcal{A}_{N_\vk,\xi}(S)]_{ijab}(\vk_i, \vk_j, \vk_a)
                \right|
                \leqslant 
                C
                \begin{cases}
                        \frac{1}{|\vq|^2} N_\vk^{-1} & \vq \neq \bm{0}\\
                        N_\vk^{-\frac13} & \vq = \bm{0}\\
                \end{cases}.  
        \end{equation*}
\end{lem}

Combining the estimates of the three error terms in \cref{eqn:decomposition_amplitude} and the initial condition $\mathcal{M}_\mathcal{K} t_0 = T^{N_\vk}_{0} = \bm{0}$ 
in both cases with and without Madelung constant correction, 
we can obtain the entrywise estimate of the error in the CCD$(n)$ amplitude calculation recursively. 
First, for the CCD$(n)$ calculation without any corrections or with partial Madelung constant corrections to either 
orbital energies or ERI contractions, we have 
\begin{equation}\label{eqn:error_amplitude}
        \left|
                [\mathcal{M}_{\mathcal{K}}t_{n} - T_{n}^{N_\vk}]_{ijab}(\vk_i, \vk_j, \vk_a)
        \right|
        \leqslant 
        C
        N_\vk^{-\frac13}.
\end{equation}
Accordingly, the average entrywise error can be estimated as 
\begin{equation*}
        \left\|
        \mathcal{M}_{\mathcal{K}}t_{n} - T_{n}^{N_\vk}
        \right\|_1 
        = 
        \dfrac{1}{N_\vk^3}
        \sum_{\vk_i,\vk_j,\vk_a\in\mathcal{K}}
        \sum_{ijab}
        |[\mathcal{M}_{\mathcal{K}}t_{n} - T_{n}^{N_\vk}]_{ijab}(\vk_i, \vk_j, \vk_a)|
        \leqslant
        C
        N_\vk^{-\frac13}.
\end{equation*} 

Secondly, for the CCD$(n)$ calculation with the Madelung constant correction to both orbital energies 
and ERI contractions, we have 
\begin{equation}\label{eqn:error_amplitude_xi}
        \left|
                [\mathcal{M}_{\mathcal{K}}t_{n} - T_{n}^{N_\vk}]_{ijab}(\vk_i, \vk_j, \vk_a)
        \right|
        \leqslant 
        C
        \begin{cases}
                \frac{1}{|\vq|^2} N_\vk^{-1} & \vq \neq \bm{0}\\
                N_\vk^{-\frac13} & \vq = \bm{0}\\
        \end{cases}.  
\end{equation}
Accordingly, the average entrywise error can be estimated as 
\begin{align*}
        \left\|
        \mathcal{M}_{\mathcal{K}}t_{n} - T_{n}^{N_\vk}
        \right\|_1 
        = & 
        \dfrac{1}{N_\vk^3}
        \sum_{\vk_i,\vk_j\in\mathcal{K}}
        \sum_{\vq \in \mathcal{K}_\vq}
        \sum_{ijab}
        |[\mathcal{M}_{\mathcal{K}}t_{n} - T_{n}^{N_\vk}]_{ijab}(\vk_i, \vk_j, \vk_i + \vq)|
        \\
        \leqslant
        & 
        C
        \frac{1}{N_\vk}
        \sum_{\vq \in \mathcal{K}_\vq\setminus\{\bm{0}\}}
        \frac{1}{|\vq|^2} N_\vk^{-1}
        + 
        C
        \frac{1}{N_\vk}
        N_\vk^{-\frac13}
        \\
        \leqslant 
        &
        CN_\vk^{-1}.
\end{align*} 

Plugging the above estimate and the error estimate for energy calculation in \cref{lem:error_energy} into \cref{eqn:decomposition_energy} then finishes the proof of \cref{thm:error_ccdn}.

\begin{rem}
        In our previous work \cite{XingLin2023}, we use the maximum entrywise norm to characterize the finite-size error in the amplitude calculation, which 
        loosens the average entrywise norm used in the error splitting in \cref{eqn:decomposition_energy_app} as 
        \[
              \| T\|_1 \leqslant \| T\|_\infty := \max_{ijab,\vk_i,\vk_j,\vk_a\in\mathcal{K}} |T_{ijab}(\vk_i,\vk_j,\vk_a)|.
        \]
        Without the Madelung constant correction, the maximum norm suffices since all entries in the amplitude error are of the same scale, as shown in \cref{eqn:error_amplitude}. 
        However, this norm is no longer sufficient for the calculation with the Madelung constant correction, bounded in \cref{eqn:error_amplitude_xi}. 
        Because the maximum entrywise error is now of scale $\Or(N_\vk^{-\frac13})$, while most entries are of scale $\Or(N_\vk^{-1})$. 
        In this case, the average entrywise norm provides a necessary and tighter estimate of amplitude error.
\end{rem}

\subsection{Proof of Lemma \ref{lem:error_amplitude_contraction_app}: error in ERI contractions}
\label{proof:error_amplitude_contraction}
According to the singularity structure of exact CCD$(n)$ amplitude in \cref{lem:thm1_nonsmooth}, we consider the ERI-contraction using an arbitrary exact amplitude 
$t \in \mathbb{T}(\Omega^*)^{n_\text{occ}\times n_\text{occ}\times n_\text{vir}\times n_\text{vir}}$.
Fixing a set of entry index $(i,j,a,b)$ and $(\vk_i,\vk_j,\vk_a)\in\mathcal{K}\times\mathcal{K}\times\mathcal{K}$, the error in the indexed ERI-contraction entry can be detailed as (by comparing \cref{eqn:amplitude_ccd} and \cref{eqn:amplitude_ccd_tdl})
\begin{align}
        & 
        \left[
                \mathcal{M}_{\mathcal{K}}\mathcal{A}_\text{TDL}(t) - \mathcal{A}_{N_\vk, \xi}(\mathcal{M}_{\mathcal{K}}t) 
        \right]_{ijab, \vk_i\vk_j\vk_a}
        \nonumber \\
        = & 
        \dfrac{1}{|\Omega^*|}
        \left[
        \mathcal{E}_{\Omega^*}
        \left(
                \sum_{kl} \braket{k\vk_k, l\vk_l | i\vk_i, j\vk_j} t_{klab}(\vk_k, \vk_l, \vk_a), 
                \mathcal{K}
        \right)
        + |\Omega^*| \xi t_{ijab}(\vk_i, \vk_j, \vk_a)
        \right]
        \nonumber \\
        & + 
                \dfrac{1}{|\Omega^*|^2}\mathcal{E}_{\Omega^*\times\Omega^*}\left(\sum_{klcd}
                \braket{k\vk_k, l\vk_l | c\vk_c, d\vk_d} t_{ijcd}(\vk_i, \vk_j, \vk_c)t_{klab}(\vk_k, \vk_l, \vk_a), \mathcal{K}\times\mathcal{K}
                \right)
        \nonumber\\
        & + \cdots,
        \label{eqn:error_amp_expansion}
\end{align}
where the constant terms cancel with each other, the first term is the error in the 4h2p linear term calculation with the Madelung constant correction, 
and the second term is the error in the 4h2p quadratic term calculation. 
The neglected ones are the errors in remaining linear and quadratic terms calculations, which can all be similarly formulated as quadrature errors 
of specific trapezoidal rules.

In the error analysis for CCD calculation without corrections in ERI contractions~\cite{XingLin2023}, it has been shown that without the Madelung constant correction 
the error entry in the ERI contractions is uniformly bounded as
\begin{equation}\label{eqn:error_contraction}
        \left|
        \left[
                \mathcal{M}_{\mathcal{K}}\mathcal{A}_\text{TDL}(t) - \mathcal{A}_{N_\vk}(\mathcal{M}_{\mathcal{K}}t) 
        \right]_{ijab, \vk_i\vk_j\vk_a}
        \right| 
        \leqslant CN_\vk^{-\frac13}.
\end{equation}
More specifically, all the quadratic terms and part of the linear terms that only contain ERIs with mismatched orbital indices contribute at most $\Or(N_\vk^{-1})$ errors in both \cref{eqn:error_amp_expansion}
and \cref{eqn:error_contraction}.  
The dominant error in \cref{eqn:error_contraction} comes from the calculation of the following six linear terms
\begin{align*}
\dfrac{1}{N_\vk}&\sum_{\vk_k\in\mathcal{K}}\sum_{kl}     \braket{k\vk_k, l\vk_l | i\vk_i, j\vk_j} t_{klab}(\vk_k, \vk_l, \vk_a), \  
&\dfrac{1}{N_\vk}\sum_{\vk_c\in\mathcal{K}}\sum_{cd}     \braket{a\vk_a, b\vk_b | c\vk_c, d\vk_d} t_{ijcd}(\vk_i, \vk_j, \vk_c),\\
-\dfrac{1}{N_\vk}&\sum_{\vk_k\in\mathcal{K}}\sum_{kc}     \braket{a\vk_a, k\vk_k | c\vk_c, i\vk_i} t_{kjcb}(\vk_k, \vk_j, \vk_c), \
&-\dfrac{1}{N_\vk}\sum_{\vk_k\in\mathcal{K}}\sum_{kc}     \braket{b\vk_b, k\vk_k | c\vk_c, j\vk_j} t_{kica}(\vk_k, \vk_i, \vk_c),\\
-\dfrac{1}{N_\vk}&\sum_{\vk_k\in\mathcal{K}}\sum_{kc}     \braket{a\vk_a, k\vk_k | c\vk_c, j\vk_j} t_{kibc}(\vk_k, \vk_i, \vk_b), \
&-\dfrac{1}{N_\vk}\sum_{\vk_k\in\mathcal{K}}\sum_{kc}     \braket{b\vk_b, k\vk_k | c\vk_c, i\vk_i} t_{kjac}(\vk_k, \vk_j, \vk_a).
\end{align*}
which contain ERIs with fully or partially matched orbital indices. 
The Madelung constant correction in \cref{eqn:error_amp_expansion} is exactly triggered in these six terms. 
In this proof, we focus on the error estimate for the 4h2p linear term (the first term above) with the correction, 
and similar analysis can be applied to all the other five terms. 

Denote the ERI-amplitude product in the 4h2p linear term with orbital indices $(k, l)$ as 
\[
        F^{kl}(\vq_1) 
        = 
        \braket{k(\vk_i - \vq_1), l(\vk_j + \vq_1) |  i\vk_i, j\vk_j} t_{klab}(\vk_i - \vq_1, \vk_j + \vq_1, \vk_a)
        =
        H_\text{eri}^{kl}(\vq_1) H_\text{amp}^{kl}(\vq_1),
\]
where $\vq_1 = \vk_i - \vk_k$ is the momentum transfer vector of the ERI. 
The 4h2p linear term calculation with the Madelung constant correction using a finite mesh $\mathcal{K}$ can be reformulated as 
\[
        \dfrac{1}{N_\vk}\sum_{\vk_k \in \mathcal{K}}\sum_{kl} F^{kl}(\vk_i - \vk_k) - \xi H_\text{amp}^{ij}(\bm{0})
        = 
        \dfrac{1}{N_\vk}\sum_{\vq_1 \in \mathcal{K}_\vq}\sum_{kl} F^{kl}(\vq_1)-  \xi H_\text{amp}^{ij}(\bm{0}),
\]
using the change of variable $\vk_k \rightarrow \vk_i - \vq_1$ and the periodicity of $F^{kl}(\vq_1)$. 
The error of this calculation compared to its TDL value, i.e., the first error term in \cref{eqn:error_amp_expansion}, can be written as 
\begin{equation*}
        \dfrac{1}{|\Omega^*|}
        \left(
        \sum_{kl} 
        \mathcal{E}_{\Omega^*}
        \left(
                F^{kl}(\vq_1), \mathcal{K}_\vq
        \right)
        +  |\Omega^*|\xi H_\text{amp}^{ij}(\bm{0})
        \right).
\end{equation*}

Previously in \cite{XingLin2023}, the quadrature error for $F^{kl}(\vq_1)$ with varying $(k,l)$ is estimated as
\begin{equation}\label{eqn:error_4h2p_linear}
        \left|
        \mathcal{E}_{\Omega^*}
        \left(
                F^{kl}(\vq_1), \mathcal{K}_\vq
        \right) 
        \right|
        \leqslant 
        C
        \begin{cases}
                N_\vk^{-1} & k\neq i, l\neq j \\
                N_\vk^{-\frac23} & k=i, l\neq j \text{ or } k\neq i, l=j\\
                N_\vk^{-\frac13} & k=i, l=j 
        \end{cases}, 
        \quad 
        \forall i,j,a,b,\ \forall \vk_i,\vk_j,\vk_a\in\mathcal{K}.
\end{equation}
As to be demonstrated next, in the case of partially matched orbital indices (e.g., $k=i, l\neq j$), this error estimate turns out to be loose when $\vk_i \neq \vk_a$ and can be further 
improved as
\begin{equation}\label{eqn:error_partial}
        \left| 
                \mathcal{E}_{\Omega^*}
                \left(
                        F^{il}(\vq_1),
                        \mathcal{K}_\vq
                \right) 
        \right|
        \leqslant 
        C
        \begin{cases}
                \frac{1}{|\vq|} N_\vk^{-1} & \vq \neq \bm{0} \\
                N_\vk^{-\frac23} & \vq = \bm{0}\\
        \end{cases},
        \quad 
        \forall i,j,a,b,\ \forall \vk_i,\vk_j,\vk_a\in\mathcal{K},
\end{equation}
where $\vq = \vk_a - \vk_i + \vG_0$ with $\vG_0 \in \mathbb{L}^*$ chosen to make $\vq \in \Omega^*$.
In the case of fully matched orbital indices ($k=i, l=j$), the Madelung constant correction is triggered in the ERI evaluation and can help remove the leading quadrature error when $\vk_i \neq \vk_a$ as
\begin{equation}\label{eqn:error_fully}
        \left| 
                \mathcal{E}_{\Omega^*}
                \left(
                        F^{ij}(\vq_1),
                        \mathcal{K}_\vq
                \right) 
                + 
                |\Omega^*| \xi H_\text{amp}^{ij}(\bm{0})
        \right|
        \leqslant 
        C
        \begin{cases}
                \frac{1}{|\vq|^2} N_\vk^{-1} & \vq \neq \bm{0} \\
                N_\vk^{-\frac13} & \vq = \bm{0}\\
        \end{cases},
        \quad 
        \forall i,j,a,b,\ \forall \vk_i,\vk_j,\vk_a\in\mathcal{K}.
\end{equation}
In the following discussion, we prove these two error estimates \cref{eqn:error_partial} and \cref{eqn:error_fully} at $\vk_i \neq \vk_a$ separately.

\subsubsection*{Error estimate for linear terms with partially matched orbital indices}
Consider a fixed set of $(i,j,a,b)$ and $\vk_i,\vk_j,\vk_a\in \mathcal{K}$ with $\vk_i \neq \vk_a$.
Assume $k = i$ and $l \neq j$. 
Our target is a sharper estimate of $\mathcal{E}_{\Omega^*}\left(F^{il}(\vq_1), \mathcal{K}_\vq\right)$. 
Since $\mathcal{K}_\vq$ and $\Omega^*$ are both inversely symmetric over $\vq_1 = \bm{0}$, the quadrature error can 
be symmetrized as 
\[
        \mathcal{E}_{\Omega^*}
        \left( 
                F^{il}(\vq_1), \mathcal{K}_\vq
        \right) 
        = 
        \frac12
        \mathcal{E}_{\Omega^*}
        \left(
                H_\text{eri}^{il}(\vq_1) H_\text{amp}^{il}(\vq_1) +
                H_\text{eri}^{il}(-\vq_1) H_\text{amp}^{il}(-\vq_1), 
                \mathcal{K}_\vq
        \right). 
\]
This symmetrized integrand can be further decomposed into two terms 
\begin{equation}\label{eqn:partial_decomp}
        \left[
        \left(
                H_\text{eri}^{il}(\vq_1)  +
                H_\text{eri}^{il}(-\vq_1)
        \right)
        H_\text{amp}^{il}(\vq_1)
        \right]
        - 
        \left[
        H_\text{eri}^{il}(-\vq_1)
        \left(
                H_\text{amp}^{il}(\vq_1) - 
                H_\text{amp}^{il}(-\vq_1) 
        \right)
        \right].
\end{equation}

For the first term in \cref{eqn:partial_decomp}, we note that the ERI term can be detailed as  
\begin{equation*}
        H_\text{eri}^{il}(\vq_1) 
        = 
        \frac{4\pi}{\abs{\Omega}} 
        \frac{
        \hat{\varrho}_{i(\vk_i-\vq_1),i\vk_i}(\bm{0}) \hat{\varrho}_{l(\vk_j + \vq_1),j\vk_j}(\bm{0})
        }{
                |\vq_1|^2
        }
        + 
        \frac{4\pi}{\abs{\Omega}}
        \sum_{\vG\in\mathbb{L}^*\setminus \{\bm{0}\}}
        \dfrac{\cdots}{|\vq_1+\vG|^2}.
\end{equation*}
The ERI nonsmoothness with $\vq_1 \in \Omega^*$ comes from the first fraction term whose numerator is smooth and scales as $\Or(|\vq_1|)$ 
near $\vq_1 = \bm{0}$ using orbital orthogonality, 
It can be verified directly that $H_\text{eri}^{il}(\vq_1)$ is singular at $\vq_1 = \bm{0}$ of order $-1$, 
and its symmetrization $H_\text{eri}^{il}(\vq_1) + H_\text{eri}^{il}(-\vq_1)$ is singular at $\vq_1 = \bm{0}$ of order $0$. 
Meanwhile, the amplitude $H_\text{amp}^{il}(\vq_1)$ is smooth everywhere in $\Omega^*$ except at $\vq_1 = -\vq$ of order $0$ according to \cref{lem:thm1_nonsmooth}. 
An error estimate lemma in \cite{XingLin2023} (restated as \cref{lem:quaderror2} in \cref{appendix:quaderror}) provides a quadrature error estimate for periodic functions in such a product form and can show that  
\[
        \left|
        \mathcal{E}_{\Omega^*}
        \left(
                \left(
                        H_\text{eri}^{il}(\vq_1)+H_\text{eri}^{il}(-\vq_1)
                \right)
                H_\text{amp}^{il}(\vq_1), \mathcal{K}_\vq
        \right)
        \right|
        \leqslant 
        C N_\vk^{-1},
\]
where constant $C$ is independent of $i,j,a,b$ and $\vk_i,\vk_j,\vk_a\in\mathcal{K}$. 

For the second term in \cref{eqn:partial_decomp}, direct application of \cref{lem:quaderror2} leads to error estimate of scale $\Or(N_\vk^{-\frac23})$ since 
$H_\text{eri}^{il}(-\vq_1)$ is singular at $\vq_1 = \bm{0}$ of order $-1$.
However, note that $H^{il}_\text{amp}(\vq_1)$ is smooth at $\vq_1 = \bm{0}$ 
(recall that $H^{il}_\text{amp}(\vq_1)$ is only singular at $\vq_1 = -\vq \neq \bm{0}$)
and thus the subtraction $H_\text{amp}^{il}(\vq_1) - H_\text{amp}^{il}(-\vq_1)$ scales as $\Or(|\vq_1|)$ near $\vq_1 = \bm{0}$. 
Multiplication by this extra $\Or(|\vq_1|)$ term improves the algebraic singularity of $H_\text{eri}^{il}(\vq_1)$ at $\vq_1 = \bm{0}$ 
and the overall product can be shown to be singular at $\vq_1 = \bm{0}$ of order $0$. 
Intuitively, this improved algebraic singularity at $\vq_1 = \bm{0}$ can lead to asymptotically smaller quadrature errors. 
To rigorously justify this statement, we generalize \cref{lem:quaderror2} to estimate the quadrature 
error for this special case.

\begin{lem}\label{lem:quaderror2_2}
        Let $f(\vx) = f_1(\vx)f_2(\vx)$ where $f_1(\vx)$ and $f_2(\vx)$ are periodic with respect to $V = [-\frac12,\frac12]^d$ and
        \begin{itemize}
                \item $f_1(\vx)$ is smooth everywhere except at $\vx = \vz_1 = \bm{0}$ of order $\gamma \leqslant -1$,
                \item $f_2(\vx)$ is smooth everywhere except at $\vx = \vz_2 \neq \bm{0}$ of order $0$,
                \item $\partial^{\valpha}_\vx f_2(\bm{0}) = 0$ for any derivative order $|\valpha| \leqslant s$. 
        \end{itemize}
        Assume $\gamma > -d$ for $f(\vx)$ to be integrable in $V$ and $\gamma+s+1 \leqslant 0$ so the leading algebraic singularity of $f(\vx)$ is at $\vx = \bm{0}$. 
        Consider an $m^d$-sized uniform mesh $\mathcal{X}$ in $V$. Assume that $\mathcal{X}$ satisfies that $\vz_1, \vz_2$ are either on the 
        mesh or $\Theta(m^{-1})$ away from any mesh points, and $m$ is sufficiently large that $|\vz_1 - \vz_2| = \Omega(m^{-1})$.
        At $\vx = \vz_1$ and $\vx = \vz_2$, $f(\vx)$ is set to $0$. 
        The trapezoidal rule using $\mathcal{X}$ has quadrature error
        \[
                \left|\mathcal{E}_V(f, \mathcal{X})\right| \leqslant C \mathcal{H}_{V, \vz_1}^{d+1}(f_1)\mathcal{H}_{V, \vz_2}^{d+1}(f_2) |\vz_1 - \vz_2|^{-(s+1)} m^{-(d + \gamma + s + 1)}.
        \]
\end{lem}
\begin{rem}
    The two factors $\mathcal{H}_{V, \vz_1}^{d+1}(f_1)$ and $\mathcal{H}_{V, \vz_2}^{d+1}(f_2)$ characterize the algebraic singularities of the two functions, 
    and their exact definition can be found in \cref{appendix:quaderror}. 
    Proof of \cref{lem:quaderror2_2} is provided in \cref{app:quaderror2_2}. 
\end{rem}

In order to utilize this result, we further decompose the second term in \cref{eqn:partial_decomp} into 
\[
        \left[
                H_\text{eri}^{il}(-\vq_1)
                \left(
                        H_\text{amp}^{il}(\vq_1) - 
                        H_\text{amp}^{il}(\bm{0}) 
                \right)
        \right]
        +
        \left[
                H_\text{eri}^{il}(-\vq_1)
                \left(
                        H_\text{amp}^{il}(0) - 
                        H_\text{amp}^{il}(-\vq_1) 
                \right)
        \right].
\]
Applying \cref{lem:quaderror2_2} to both terms with $\gamma = -1$ and $s = 0$ gives 
\[
        \left|
        \mathcal{E}_{\Omega^*}\left(
        H_\text{eri}^{il}(-\vq_1)
        \left(
                H_\text{amp}^{il}(\vq_1) - 
                H_\text{amp}^{il}(-\vq_1) 
        \right), 
        \mathcal{K}_\vq
        \right)
        \right|
        \leqslant 
        C \frac{1}{|\vq|}
        N_\vk^{-1},
\] 
where constant $C$ is independent of $i,j,a,b$ and $\vk_i, \vk_j, \vk_a \in \mathcal{K}$. 
Combining the above quadrature error estimates for the two terms in \cref{eqn:partial_decomp}, 
we prove a tighter error estimate at $\vq \neq\bm{0}$ shown in \cref{eqn:error_partial} while the previous result in \cref{eqn:error_4h2p_linear} at $\vq = \bm{0}$ still holds.

\subsubsection*{Error estimate for linear terms with fully matched orbital indices}

Let $k=i, l = j$ and $\vk_i \neq \vk_a$ and consider the corresponding calculation in the 4h2p linear term.
The Madelung constant correction is triggered in the ERI evaluation at $\vk_k = \vk_i$ or 
equivalently at $\vq_1 = \bm{0}$ and the calculation can be written as 
\begin{align*}
        & 
        \dfrac{1}{N_\vk}\sum_{\vq_1 \in\mathcal{K}_\vq}
        F^{ij}(\vq_1)
        - 
        \xi
        H_\text{amp}^{ij}(\bm{0})
        \\
        & \qquad =  
        \dfrac{1}{N_\vk}\sum_{\vq_1 \in\mathcal{K}_\vq}
        \left(
        H_{\text{eri}}^{ij}(\vq_1)H_{\text{amp}}^{ij}(\vq_1)
        -
        h_\sigma(\vq_1)      H^{ij}_{\text{amp}}(\bm{0})
        \right)
        +  
        \dfrac{1}{|\Omega^*|}\int_{\Omega^*}\ud\vq_1 h_\sigma(\vq_1)   H_{\text{amp}}^{ij}(\bm{0}) + \Or(N_\vk^{-1}),
\end{align*}
which uses the expansion of the Madelung constant with an arbitrary fixed parameter $\sigma > 0$ in \cref{eqn:madelung_sigma}.
Note that the prefactor in the $\Or(N_\vk^{-1})$ remainder term above only depends on $\sigma$.
The right hand side of the above reformation is equivalent to a singularity subtraction method that decomposes the original integrand into two terms,
\begin{equation*}
        F^{ij}(\vq_1) 
        = 
        H_{\text{eri}}^{ij}(\vq_1)H_{\text{amp}}^{ij}(\vq_1) 
        = 
        \left( H_{\text{eri}}^{ij}(\vq_1)H_{\text{amp}}^{ij}(\vq_1) - h_\sigma(\vq_1)H_\text{amp}^{ij}(\bm{0}) \right) 
        + h_\sigma(\vq_1)H_\text{amp}^{ij}(\bm{0}),
\end{equation*}
and then computes the numerical quadrature of the first term and the exact integral of the second term. 
The Madelung-corrected calculation thus has quadrature error only from the first term as 
\begin{align*}
        \dfrac{1}{|\Omega^*|}
        \mathcal{E}_{\Omega^*}
        \left(
                F^{ij}(\vq_1), \mathcal{K}_\vq
        \right)
        + \xi 
        H_\text{amp}^{ij}(\bm{0})
        & = 
        \dfrac{1}{|\Omega^*|}
        \mathcal{E}_{\Omega^*}\left(
        H_{\text{eri}}^{ij}(\vq_1)H_{\text{amp}}^{ij}(\vq_1)
        -
        h_\sigma(\vq_1)      H_{\text{amp}}^{ij}(\bm{0}), 
\mathcal{K}_\vq
\right) + \Or(N_\vk^{-1}).
\end{align*}

Following a similar approach as in the partially matched case, the effective integrand after the Madelung constant correction above can be split 
into two terms as 
\begin{equation}\label{eqn:madelung_decomposition}
        \left[
        \left(
                H_{\text{eri}}^{ij}(\vq_1) - h_\sigma(\vq_1)
        \right) H_{\text{amp}}^{ij}(\vq_1)
        \right]
        +
        \left[
        h_\sigma(\vq_1)
        \left(
                H_{\text{amp}}^{ij}(\vq_1) - H_{\text{amp}}^{ij}(\bm{0})
        \right)
        \right].
\end{equation}

For the first term in \cref{eqn:madelung_decomposition}, the subtraction part can be detailed as 
\[
H_{\text{eri}}^{ij}(\vq_1)
-
h_\sigma(\vq_1)
= 
\frac{4\pi}{\abs{\Omega}} 
\frac{
\hat{\varrho}_{i(\vk_i-\vq_1),i\vk_i}(\bm{0}) \hat{\varrho}_{j(\vk_j + \vq_1),j\vk_j}(\bm{0}) - e^{-\sigma |\vq_1|^2}
}{
        |\vq_1|^2
}
+ 
\frac{4\pi}{\abs{\Omega}}
\sum_{\vG\in\mathbb{L}^*\setminus \{\bm{0}\}}
\dfrac{\cdots}{|\vq_1+\vG|^2},
\]
where the numerator of the first fraction scales as $\Or(|\vq_1|)$ near $\vq_1 = \bm{0}$. 
This subtraction is thus singular at $\vq_1 = \bm{0}$ of order $-1$ and shares a similar form as the ERI with partially matched orbital indices. 
Using the earlier inverse-symmetry analysis for partially matched case, the quadarture error of the first term in \cref{eqn:madelung_decomposition} when $\vk_i \neq \vk_a$ can be estimated as 
\[
        \left|
        \mathcal{E}_{\Omega^*}
        \left(
        \left(
                H_{\text{eri}}^{ij}(\vq_1) - h_\sigma(\vq_1)
        \right) H_{\text{amp}}^{ij}(\vq_1), 
        \mathcal{K}_\vq
        \right)
        \right|
        \leqslant 
        C \frac{1}{|\vq|}N_\vk^{-1}. 
\]

For the second term in \cref{eqn:madelung_decomposition}, we exploit the inverse symmetry of $\mathcal{K}_\vq$ and $\Omega^*$ again and its quadrature error 
equals to that of its symmetrized version as
\[
        \frac{1}{2}
        h_\sigma(\vq_1)
        \left(
                H_{\text{amp}}^{ij}(\vq_1) - 2H_{\text{amp}}^{ij}(\bm{0}) + H_{\text{amp}}^{ij}(-\vq_1)
        \right).
\]
This formula uses $h_\sigma(\vq_1) = h_\sigma(-\vq_1)$ after symmetrization. 
Note that $h_\sigma(\vq_1)$ is singular at $\vq_1 = \bm{0}$ of order $-2$ while the term in the parenthesis scales as $\Or(|\vq_1|^2)$
near $\vq_1 = \bm{0}$ using the smoothness of $H_\text{amp}^{ij}(\vq_1)$ at $\vq_1 = \bm{0}$.  
Therefore, the overall function above is singular at $\vq_1 = \bm{0}$ of order $0$.
To fit the integrand form in \cref{lem:quaderror2_2}, we further decompose the symmetrized integrand above into 
\[
        \left[
        h_\sigma(\vq_1)
        \left(
                H_{\text{amp}}^{ij}(\vq_1) - H_{\text{amp}}^{ij}(\bm{0}) - \vq_1\cdot \nabla H_{\text{amp}}^{ij}(\bm{0})
        \right)
        \right]
        + 
        \left[
        h_\sigma(\vq_1)
        \left(
                H_{\text{amp}}^{ij}(-\vq_1) - H_{\text{amp}}^{ij}(\bm{0}) - (-\vq_1)\cdot \nabla H_{\text{amp}}^{ij}(\bm{0})
        \right)
        \right].
\] 
Applying \cref{lem:quaderror2_2} to these two terms separately with $\gamma = -2$ and $s = 1$ gives  
\[
        \left|
        \mathcal{E}_{\Omega^*}
        \left(
        \frac{1}{2}
        h_\sigma(\vq_1)
        \left(
                H_{\text{amp}}^{ij}(\vq_1) - 2H_{\text{amp}}^{ij}(\bm{0}) + H_{\text{amp}}^{ij}(-\vq_1)
        \right),
        \mathcal{K}_\vq
        \right)
        \right|
        \leqslant 
        C \frac{1}{|\vq|^2}N_\vk^{-1}. 
\]

Combining the above quadrature error estimates for the two terms in \cref{eqn:madelung_decomposition}, 
we prove a tighter error estimate at $\vq \neq\bm{0}$ shown in \cref{eqn:error_fully} while the previous result in \cref{eqn:error_4h2p_linear} at $\vq = \bm{0}$ still holds,

The same error estimate can be obtained for all the six linear terms that trigger the Madelung constant correction, 
and the remaining linear and quadratic terms contribute at most $\Or(N_\vk^{-1})$ quadrature error.
Gathering the error estimates for all these terms together and plugging into \cref{eqn:error_amp_expansion}, we finish the proof.

\subsection{Proof of \cref{lem:error_orbital_energy}: error in orbital energies}
\label{proof:error_orbital_energy}

In the TDL, the orbital energy $\varepsilon_{n\vk}^{N_\vk}$ with any fixed $n$ and $\vk\in\mathcal{K}$ converges to
\begin{equation*}\label{eqn:orbital_energy_tdl}
        \varepsilon_{n\vk}^\text{TDL} = 
        \braket{n\vk | \hat{H}_0| n\vk} + 
        \dfrac{1}{|\Omega^*|}\int_{\Omega^*}\ud\vk_i \sum_{i} \left(2\braket{i\vk_i, n\vk | i\vk_i, n\vk} - \braket{i\vk_i, n\vk | n\vk, i\vk_i}\right).
\end{equation*} 
Comparing this exact orbital energy with its finite-size calculation in \cref{eqn:orbital_energy}, the finite-size error without any corrections can be written as 
\begin{equation}\label{eqn:error_orbital_expansion}
        \varepsilon_{n\vk}^\text{TDL} - 
        \varepsilon_{n\vk}^{N_\vk}
        = 
        \dfrac{2}{|\Omega^*|}
        \mathcal{E}_{\Omega^*}\left(
                \sum_{i}\braket{i\vk_i, n\vk | i\vk_i, n\vk}, 
                \mathcal{K}
        \right)
        - 
        \dfrac{1}{|\Omega^*|}
        \mathcal{E}_{\Omega^*}\left(
                \sum_{i}\braket{i\vk_i, n\vk | n\vk, i\vk_i}, 
                \mathcal{K}
        \right).
\end{equation}

For the first quadrature error in \cref{eqn:error_orbital_expansion}, i.e., the finite-size error in the direct term, the ERI with each $i$ 
can be specified as 
\[
        \braket{i\vk_i, n\vk | i\vk_i, n\vk}
        = 
        \dfrac{4\pi}{|\Omega|}
        \sum_{\vG \in \mathbb{L}^*\setminus\{\bm{0}\}}
        \dfrac{\hat{\rho}_{i\vk_i,i\vk_i}(\vG) \hat{\rho}_{n\vk,n\vk}(-\vG)}{|\vG|^2}.
\]
Note that the momentum transfer vector of this ERI is always zero and the singular fraction term in this ERI is set to $0$ by definition. 
As a result, the integrand is periodic and smooth with respect to $\vk_i\in\Omega^*$. 
Therefore, the quadrature error of the direct term calculation thus decays super-algebraically according to \cref{lem:quaderror0} as  
\[
        \left|
        \mathcal{E}_{\Omega^*}\left( 
                \sum_{i}\braket{i\vk_i, n\vk | i\vk_i, n\vk}, 
                \mathcal{K}
        \right)
        \right|
        \leqslant C_l N_\vk^{-l}, \quad \forall l > 0.
\]

For the second quadrature error in \cref{eqn:error_orbital_expansion}, i.e., the finite-size error in the exchange term,
the ERI with each $i$ can be written as  
\begin{equation*}
        \braket{i\vk_i, n\vk | n\vk, i\vk_i}
        = 
        \dfrac{4\pi}{|\Omega|}
        \dfrac{\hat{\rho}_{i\vk_i,n\vk}(\bm{0}) \hat{\rho}_{n\vk,i\vk_i}(\bm{0})}{|\vq|^2}
        + 
        \dfrac{4\pi}{|\Omega|}
        \sum_{\vG \in \mathbb{L}^*\setminus\{\bm{0}\}}
        \dfrac{\hat{\rho}_{i\vk_i,i\vk_i}(\vG) \hat{\rho}_{n\vk,n\vk}(-\vG)}{|\vq + \vG|^2},
\end{equation*}
which is singular at $\vq : = \vk - \vk_i = \bm{0}$ of order $-2$ when $n = i$ and $0$ otherwise. 
An error estimate lemma in \cite{XingLin2023} (restated as \cref{lem:quaderror1} in \cref{appendix:quaderror}) gives a tight quadrature error estimate for such periodic functions 
with one point of algebraic singularity, and its application to the above integrand gives 
\begin{equation*}
        \left|
                \mathcal{E}_{\Omega^*}
                \left(
                        \braket{i\vk_i, n\vk | n\vk, i\vk_i}, \mathcal{K}
                \right)
        \right|
        \leqslant
        C
        \begin{cases}
                N_\vk^{-\frac13} & n = i \\
                N_\vk^{-1} & n \neq i \\
        \end{cases}.
\end{equation*}

Combining the estimates of the two error terms in \cref{eqn:error_orbital_expansion}, we obtain the 
overall finite-size error estimate for orbital energies without any corrections as 
\begin{equation*}
        \left|
                \varepsilon_{n\vk}^\text{TDL} - \varepsilon_{n\vk}^{N_\vk}
        \right|
        \leqslant 
        C
        \begin{cases}
                N_\vk^{-\frac13} & \text{$n$ is occupied} \\
                N_\vk^{-1} & \text{$n$ is virtual} 
        \end{cases}.
\end{equation*}

From the above analysis, the dominant finite-size error in an occupied orbital energy lies in the 
calculation of the exchange term with $i = n$. 
In the Madelung-corrected orbital energy $\varepsilon_{n\vk}^{N_\vk,\xi}$, the correction is applied to this exchange term as (ignoring the prefactor $-|\Omega^*|^{-1}$)  
\[
        \mathcal{Q}_{\Omega^*}
        \left(
                \braket{n\vk_i, n\vk | n\vk, n\vk_i}, \mathcal{K}
        \right)
        \rightarrow
        \mathcal{Q}_{\Omega^*}
        \left(
                \braket{n\vk_i, n\vk | n\vk, n\vk_i}, \mathcal{K}
        \right)
        - |\Omega^*|\xi.
\]
Applying the change of variable $\vk_i \rightarrow \vk - \vq$ and using the periodicity of the ERI with respect to $\vk_i$, 
this corrected calculation can be reformulated as 
\begin{align*}
        & \mathcal{Q}_{\Omega^*}
        \left(
                \braket{n\vk_i, n\vk | n\vk, n\vk_i}, \mathcal{K}
        \right)
        - |\Omega^*|\xi
        \\
        & \qquad = 
        \mathcal{Q}_{\Omega^*}
        \left(
                \braket{n(\vk - \vq), n\vk | n\vk, n(\vk - \vq)}
                - 
                h_\sigma(\vq), \mathcal{K}_\vq
        \right)
        +  
        \int_{\Omega^*}\ud\vq_1 h_\sigma(\vq_1)  + \Or(N_\vk^{-1}),
\end{align*}
using the expansion of the Madelung constant with any fixed $\sigma > 0$ in \cref{eqn:madelung_sigma}. 
The quadrature error after the correction can be written as 
\begin{align*}
        &
        \mathcal{E}_{\Omega^*}
        \left(
                \braket{n\vk_i, n\vk | n\vk, n\vk_i}, \mathcal{K}
        \right)
        + |\Omega^*|\xi
        \\
        & \qquad 
        = 
        \mathcal{E}_{\Omega^*}
        \left(
                \braket{n(\vk - \vq), n\vk | n\vk, n(\vk - \vq)}
                - 
                h_\sigma(\vq), \mathcal{K}_\vq
        \right)
        + \Or(N_\vk^{-1}). 
\end{align*}
The effective integrand above after the correction can be detailed as 
\[
        \braket{n(\vk - \vq), n\vk | n\vk, n(\vk - \vq)}
                - 
                h_\sigma(\vq) 
        = \frac{4\pi}{\abs{\Omega}} 
        \frac{
        |\hat{\varrho}_{n(\vk-\vq),n\vk}(\bm{0})|^2 - e^{-\sigma |\vq|^2}
        }{
                |\vq|^2
        }
        + 
        \frac{4\pi}{\abs{\Omega}}
        \sum_{\vG\in\mathbb{L}^*\setminus \{\bm{0}\}}
        \dfrac{\cdots}{|\vq+\vG|^2}.
\]
The integrand singularity comes from the first fraction and is of order $-1$ since the ERI and $h_\sigma(\vq)$ share the 
same leading singular term. 
Similar to the analysis in the ERI contraction, we can then combine this singularity subtraction with the inverse symmetry of $\Omega^*$ and $\mathcal{K}_\vq$ 
to show that 
\[
        \left|
        \mathcal{E}_{\Omega^*}
        \left(
                \braket{n(\vk - \vq), n\vk | n\vk, n(\vk - \vq)}
                - 
                h_\sigma(\vq), \mathcal{K}_\vq
        \right)
        \right| 
        \leqslant C N_\vk^{-1}. 
\]

Combining this estimate with the above estimate for all the remaining direct and exchange terms in $\varepsilon_{n\vk}^{N_\vk, \xi}$, we obtain the finite-size error 
estimate for orbital energies with the Madelung constant correction as 
\[
        \left|
                \varepsilon_{n\vk}^\text{TDL} - \varepsilon_{n\vk}^{N_\vk,\xi}
        \right|
        \leqslant 
        C N_\vk^{-1}. 
\]

\subsection{Proof of \cref{lem:error_amplitude_accumulation}: error from previous iteration}
\label{proof:error_accu_exact}
Fixing a set of entry index $(i,j,a,b)$ and $(\vk_i,\vk_j,\vk_a)\in\mathcal{K}\times\mathcal{K}\times\mathcal{K}$,
the ERI-contraction entry $\left[\mathcal{A}_{N_\vk,\xi}(T)\right]_{ijab,\vk_i\vk_j\vk_a}$ can be detailed as 
\begin{align*}
        & 
        \left[
                \mathcal{A}_{N_\vk, \xi}(T)
        \right]_{ijab}(\vk_i\vk_j\vk_a)
        =
        \braket{a\vk_a, b\vk_b | i\vk_i, j\vk_j}
        \nonumber \\ 
        & 
        \qquad + 
        \dfrac{1}{N_\vk}\sum_{\vk_k \in \mathcal{K}}\sum_{kl} \braket{k\vk_k, l\vk_l | i\vk_i, j\vk_j} T_{klab}(\vk_k, \vk_l, \vk_a)
        \nonumber\\
        & 
        \qquad + 
        \dfrac{1}{N_\vk^2}\sum_{\vk_k\vk_c \in \mathcal{K}} \sum_{klcd} \braket{k\vk_k, l\vk_l | c\vk_c, d\vk_d} 
        T_{ijcd}(\vk_i, \vk_j, \vk_c)T_{klab}(\vk_k, \vk_l, \vk_a)
        \nonumber \\
        & \qquad + \cdots
        + 2\xi T_{ijab}(\vk_i, \vk_j, \vk_a),
        \nonumber 
\end{align*}
where the neglected terms are all the other linear and quadratic terms and the Madelung constant corrections to different terms are collected together at the end. 

In the subtraction $\mathcal{A}_{N_\vk, \xi}(T) - \mathcal{A}_{N_\vk, \xi}(S)$, the constant terms cancel each other.  
The subtraction between the two Madelung constant correction terms can be estimated directly as 
\begin{align*}
        \left|
                2\xi [T - S]_{ijab}(\vk_i,\vk_j,\vk_a)
        \right|
        \leqslant 
        C 
        N_\vk^{-\frac13}
        \begin{cases}
                \frac{1}{|\vq|^2} N_\vk^{-1} & \vq \neq \bm{0}\\
                N_\vk^{-\frac13} & \vq = \bm{0}\\
        \end{cases},  
\end{align*}
using the fact that $\xi = \Or(N_\vk^{-\frac13})$. 

The subtraction between the two 4h2p linear terms can be formulated and bounded as
\begin{align}
        & 
        \left|
        \dfrac{1}{N_\vk}\sum_{\vk_k \in \mathcal{K}}\sum_{kl} \braket{k\vk_k, l\vk_l | i\vk_i, j\vk_j} 
        \left[
                T - S
        \right]_{klab}(\vk_k, \vk_l, \vk_a)
        \right|
        \nonumber \\ 
        \leqslant
        & 
        \left(
                \frac{1}{N_\vk}\sum_{\vk_k\in \mathcal{K}\setminus\{\vk_i,\vk_a\}}
                +
                \frac{1}{N_\vk}\delta_{\vk_k, \vk_i}
                +
                \frac{1}{N_\vk}\delta_{\vk_k, \vk_a}
        \right)
        \left|
        \sum_{kl} \braket{k\vk_k, l\vk_l | i\vk_i, j\vk_j} 
        \left[
                T - S
        \right]_{klab}(\vk_k, \vk_l, \vk_a)
        \right|.
        \label{eqn:subtraction_linear}
\end{align}
For the term with $\vk_k = \vk_i$ (i.e., $\delta_{\vk_k, \vk_i}$), we have 
\begin{equation*}
        \frac{1}{N_\vk}
        \left|
                \sum_{kl} \braket{k\vk_i, l\vk_j | i\vk_i, j\vk_j} 
                \left[
                        T - S
                \right]_{klab}(\vk_i, \vk_j, \vk_a)
        \right| 
        \leqslant 
        \frac{1}{N_\vk}
        \sum_{kl}
        C N_\vk^{-\frac13}
        \leqslant C N_\vk^{-\frac43},
\end{equation*}
where the ERI definition at zero momentum transfer skips the singular fraction and is $\Or(1)$, and according to the assumption on $T -S$ 
it always holds that  
\begin{equation}\label{eqn:entry_bound}
\max_{ijab,\vk_i,\vk_j,\vk_a\in\mathcal{K}} | [T - S]_{ijab}(\vk_i, \vk_j, \vk_a) | \leqslant 
C N_\vk^{-\frac13}.
\end{equation}
For the term with $\vk_k = \vk_a$ (i.e., $\delta_{\vk_k, \vk_a}$), its estimate is the same as the term above when $\vk_a = \vk_i$. 
When $\vk_a \neq \vk_i \in \mathcal{K}$, we have $|\vq| \geqslant C N_\vk^{-\frac13}$ and 
\begin{equation*}
        \frac{1}{N_\vk}
        \left|
                \sum_{kl} \braket{k\vk_a, l\vk_b | i\vk_i, j\vk_j} 
                \left[
                        T - S
                \right]_{klab}(\vk_a, \vk_b, \vk_a)
        \right| 
        \leqslant 
        \frac{1}{N_\vk}
        \sum_{kl}
        \frac{C}{|\vq|^2}
        N_\vk^{-\frac13}
        \leqslant 
        C
        \frac{1}{|\vq|}
        N_\vk^{-1}.
\end{equation*}
For the first summation term in \cref{eqn:subtraction_linear}, we introduce the change of variable $\vk_k \rightarrow \vk_i - \vq_1$ and write it as 
\begin{equation*}
        \frac{1}{N_\vk}\sum_{\vq_1\in \mathcal{K}_\vq\setminus\{\bm{0}, -\vq\}}
        \left|
                \sum_{kl} \braket{k(\vk_i - \vq_1), l(\vk_j + \vq_1) | i\vk_i, j\vk_j} 
                \left[
                        T - S
                \right]_{klab}(\vk_i-\vq_1, \vk_j+\vq_1, \vk_a)
        \right|. 
\end{equation*}
When $\vq = \bm{0}$, this term can be further bounded using \cref{eqn:entry_bound} by 
\begin{align*}
        (*) 
        & \leqslant 
        \frac{1}{N_\vk}\sum_{\vq_1\in \mathcal{K}_\vq\setminus\{\bm{0}\}}
        \frac{C}{|\vq_1|^2}
        N_\vk^{-\frac13}
        \leqslant 
        C
        \int_{\Omega^*}\ud\vq_1
        \frac{1}{|\vq_1|^2}
        N_\vk^{-\frac13}
        \leqslant 
        C N_\vk^{-\frac13}. 
\end{align*}
When $\vq \neq\bm{0}$, this term can be further bounded as
\begin{align*}
        (*)
        & 
        \leqslant 
        C
        \frac{1}{N_\vk}\sum_{\vq_1\in \mathcal{K}_\vq\setminus\{\bm{0}, -\vq\}}
        \frac{1}{|\vq_1|^2}
        \sum_{\vG \in \mathbb{L}_{\bm{0}}^*}
        \frac{1}{|\vq_1 + \vq + \vG|^2} N_\vk^{-1}
        \\
        &\leqslant 
        C
        N_\vk^{-1}
        \int_{\Omega^*}\ud\vq_1
        \sum_{\vG \in \mathbb{L}_{\bm{0}}^*}
        \frac{1}{|\vq_1|^2}
        \frac{1}{|\vq_1 + \vq + \vG|^2}
        \\
        &\leqslant  
        C
        N_\vk^{-1}
        \sum_{\vG \in \mathbb{L}_{\bm{0}}^*}
        \frac{1}{|\vq + \vG|^2}
        \leqslant C \frac{1}{|\vq|^2} N_\vk^{-1},
\end{align*}
where $\mathbb{L}^*_{\bm{0}}$ denotes the set of 27 lattice vectors in $\mathbb{L}^*$ around the origin. 
The first inequality is based on the lemma assumption on $T-S$ that 
\begin{equation*}
        \left|
        \left[
                T - S
        \right]_{klab}(\vk_i-\vq_1, \vk_j+\vq_1, \vk_a) 
        \right|
        \leqslant 
        C \frac{1}{|\vq_1 + \vq + \vG_0|^2} N_\vk^{-1}
        \leqslant 
        C \sum_{\vG\in\mathbb{L}_{\bm{0}}^*}\frac{1}{|\vq_1 + \vq + \vG|^2} N_\vk^{-1},
\end{equation*}
where $\vG_0\in\mathbb{L}_{\bm{0}}^*$ is the unique lattice vector that makes $\vq_1 + \vq + \vG_0\in\Omega^*$. 
The third inequality can be obtained from the nonsmoothness characterization of function 
\[
        f(\vz) = \int_{\Omega^*}\ud\vq_1 \frac{1}{|\vq_1|^2}\frac{1}{|\vq_1 + \vz|^2}. 
\] 
Using Lemma 11 in \cite{XingLin2023}, $f(\vz)$ is singular only at $\vz = \bm{0}$ of order $-2$ and
its value at $\vz \in \Omega^*\setminus\{\bm{0}\}$ is bounded by $C|\vz|^{-2}$.

Based on the above estimates of the first term in \cref{eqn:subtraction_linear}, we obtain the estimate of the error accumulation 
in the 4h2p linear term calculation as 
\begin{equation*}
        \left|
                \dfrac{1}{N_\vk}\sum_{\vk_k \in \mathcal{K}}\sum_{kl} \braket{k\vk_k, l\vk_l | i\vk_i, j\vk_j} 
                \left[
                        T - S
                \right]_{klab}(\vk_k, \vk_l, \vk_a)
        \right|
        \leqslant 
        C
        \begin{cases}
                \frac{1}{|\vq|^2} N_\vk^{-1} & \vq \neq \bm{0}\\
                N_\vk^{-\frac13} & \vq = \bm{0}\\
        \end{cases}.  
\end{equation*}

The same analysis can be applied to all the similar linear terms that contain ERIs with matched orbital indices. 
For other linear terms, the analysis can be done similarly and they all contribute at most $\Or(N_\vk^{-1})$ error to the overall subtraction. 
Taking the subtraction between the following 3h3p linear terms as an example, 
it can be formulated and bounded as 
\begin{align*}
        & 
        \left| 
                \dfrac{1}{N_\vk} \sum_{\vk_k \in \mathcal{K}}\sum_{kc} 
                \braket{a\vk_a, k\vk_k | i\vk_i, c\vk_c } [T - S]_{kjbc}(\vk_k, \vk_j, \vk_b)
        \right|
        \\
        \leqslant & 
        C 
        \left(
                \dfrac{1}{N_\vk} \sum_{\vk_k \in \mathcal{K}\setminus\{\vk_b\}}\sum_{kc} 
                +
                \dfrac{1}{N_\vk} \delta_{\vk_k, \vk_b}
        \right)
        \left|
        [T - S]_{kjbc}(\vk_k, \vk_j, \vk_b)
        \right|
        \\
        \leqslant & 
        C 
        \left(
                \dfrac{1}{N_\vk} \sum_{\vq_1 \in \mathcal{K}_\vq\setminus\{\bm{0}\}} \dfrac{1}{|\vq_1|^2} N_\vk^{-1}
                + 
                \dfrac{1}{N_\vk} N_\vk^{-\frac13}
        \right) 
        \leqslant C N_\vk^{-1}.
\end{align*}

For the subtraction between quadratic terms, we consider the 4h2p quadratic term as an example. 
The subtraction between the two 4h2p quadratic terms can be formulated and bounded as 
\begin{align*}
        & 
        \left|
        \dfrac{1}{N_\vk^2}\sum_{\vk_k\vk_c \in \mathcal{K}} \sum_{klcd} \braket{k\vk_k, l\vk_l | c\vk_c, d\vk_d} 
        \left(
        T_{ijcd}(\vk_i, \vk_j, \vk_c)T_{klab}(\vk_k, \vk_l, \vk_a)
        -
        S_{ijcd}(\vk_i, \vk_j, \vk_c)S_{klab}(\vk_k, \vk_l, \vk_a)
        \right)
        \right|
        \\
        \leqslant & 
        C
        \dfrac{1}{N_\vk^2}\sum_{\vk_k\vk_c \in \mathcal{K}} \sum_{klcd} 
        \left|
                T_{ijcd}(\vk_i, \vk_j, \vk_c)
                [T - S]_{klab}(\vk_k, \vk_l, \vk_a)
                + 
                [T - S]_{ijcd}(\vk_i,\vk_j,\vk_c) S_{klab}(\vk_k, \vk_l, \vk_a)
        \right|
        \\
        \leqslant & 
        C
        \dfrac{1}{N_\vk^2}\sum_{\vk_k\vk_c \in \mathcal{K}} \sum_{klcd} 
        \left|
                [T - S]_{klab}(\vk_k, \vk_l, \vk_a)
        \right|
        +
        C
        \dfrac{1}{N_\vk^2}\sum_{\vk_k\vk_c \in \mathcal{K}} \sum_{klcd} 
        \left|
                [T - S]_{ijcd}(\vk_i, \vk_j, \vk_c)
        \right|
        \\
        \leqslant & 
        C
        N_\vk^{-1}. 
\end{align*}
Similar analysis can be done to all the remaining quadratic terms, and they all contribute at most $\Or(N_\vk^{-1})$ error
to the overall subtraction. 
Gathering all the estimates above together, we finish the proof.

\section{Proof of Theorem \ref{thm:error_ccd_madelung}}
\label{app:error_ccd}

To guarantee the convergence and control the regularity of CCD$(n)$ calculations with $n \rightarrow\infty$, we introduce additional 
technical assumptions similar to those in \cite{XingLin2023}. 
One key difference is that \cite{XingLin2023} measures the error in the amplitude calculation using the maximum entrywise norm. When Madelung constant correction is used, the error of the amplitude should be measured by the average entrywise norm (related to the $L^1$ norm) instead, which is denoted by
\begin{align*}
        \|T\|_1 
        & = 
        \dfrac{1}{N_\vk^3}\sum_{\vk_i,\vk_j,\vk_a\in\mathcal{K}}\sum_{ijab} \left| T_{ijab}(\vk_i,\vk_j,\vk_a) \right|,
        \\
        \|t\|_1 
        & = 
        \dfrac{1}{|\Omega^*|^3}
        \int_{\Omega^*\times\Omega^*\times\Omega^*}\ud\vk_i\ud\vk_j\ud\vk_a
        \sum_{ijab} |t_{ijab}(\vk_i,\vk_j,\vk_a)|,
\end{align*}
where $T$ and $t$ denote a generic amplitude computed using mesh $\mathcal{K}$ and in the TDL, respectively.

Assume that the CCD amplitude equations using a sufficiently fine mesh $\mathcal{K}$ (with or without the Madelung constant correction)  
and in the TDL have unique solutions and denote the solutions as $T_{*}^{N_\vk}$ and $t_*$ respectively. 
Convergence of the CCD$(n)$ amplitudes to the CCD amplitude is defined in the $\|\cdot\|_1$-norm sense as 
\[
        \lim_{n\rightarrow \infty} \|T_n^{N_\vk} - T_*^{N_\vk}\|_1 = 0 \qquad \text{and} \qquad 
        \lim_{n\rightarrow \infty} \|t_n - t_*\|_1 = 0.
\] 

We impose a sufficient condition that guarantees the convergence of fixed point iterations by requiring the target mapping to be contractive 
in a domain that contains both the solution point and the initial guess.

\begin{assumption}\label{assump:contraction}
        For $N_\vk$ sufficiently large, the following statements hold:
\begin{enumerate}

        \item The exact CCD($n$) amplitude $t_n$ converges to the CCD amplitude $t_*$ point-wisely as $n \rightarrow\infty$, i.e., 
        \begin{equation}\label{eqn:assumption_ccd_1}
            \lim_{n\rightarrow\infty} [t_n]_{ijab}(\vk_i, \vk_j, \vk_a) = [t_*]_{ijab}(\vk_i, \vk_j, \vk_a), \qquad \forall i, j, a, b, \vk_i, \vk_j, \vk_a. 
        \end{equation}
        \item $(\varepsilon^{N_\vk})^{-1}\mathcal{A}_{N_\vk}$  is a contraction map in some domain 
        $\mathbb{B}_{N_\vk}\subset \mathbb{C}^{n_\text{occ}\times n_\text{occ} \times n_\text{vir} \times n_\text{vir}\times N_\vk\times N_\vk \times N_\vk}$ that 
        contains $T_*^{N_\vk}$ and the initial guess $\bm{0}$, i.e., 
        \begin{equation}\label{eqn:assumption_ccd_2}
                \begin{split}
                        & (\varepsilon^{N_\vk})^{-1}\mathcal{A}_{N_\vk}(T) \in \mathbb{B}_{N_\vk}, \quad \forall T \in \mathbb{B}_{N_\vk},
                        \\
                        & 
                        \|
                                (\varepsilon^{N_\vk})^{-1}
                                \left(
                                        \mathcal{A}_{N_\vk}(T) - 
                                        \mathcal{A}_{N_\vk}(S) 
                                \right)
                        \|_1
                        \leqslant L \|T - S\|_1,\quad  \forall T,S \in \mathbb{B}_{N_\vk},
                \end{split}
        \end{equation}
        with a Lipschitz constant $L < 1$. 

        \item  The exact CCD$(n)$ amplitude $t_n$ and the domain $\mathbb{B}_{N_\vk}$ above satisfy that 
        \begin{equation}\label{eqn:assumption_ccd_3}
                \mathcal{M}_\mathcal{K}t_n
                \in 
                \mathbb{B}_{N_\vk}, \quad \forall n > 0.
        \end{equation}

        \item 
        For all the amplitudes $\{t_n\}$, there exists a constant $C$ such that 
        \begin{equation}\label{eqn:assumption_ccd_4}
                \begin{split}
                &
                \left\|
                        \mathcal{M}_{\mathcal{K}}\mathcal{A}_\text{TDL}(t_{n}) - 
                        \mathcal{A}_{N_\vk,\xi}(\mathcal{M}_{\mathcal{K}}t_{n})
                \right\|_1
                \leqslant C N_\vk^{-1}, 
                \quad 
                \forall n > 0,
                \\
                &
                \left\|
                        \mathcal{M}_{\mathcal{K}}\mathcal{A}_\text{TDL}(t_{n}) - 
                        \mathcal{A}_{N_\vk}(\mathcal{M}_{\mathcal{K}}t_{n})
                \right\|_1
                \leqslant C N_\vk^{-\frac13}, 
                \quad 
                \forall n > 0.
                \end{split}
        \end{equation}

\end{enumerate}
        Note that when consider the finite-size calculation with Madelung constant correction, 
        the components $\varepsilon^{N_\vk}$ or $\mathcal{A}_{N_\vk}$ in the second assumption 
        need to be changed to $\varepsilon^{N_\vk, \xi}$ or $\mathcal{A}_{N_\vk, \xi}$ accordingly. 
\end{assumption}

\begin{rem}
        The second assumption guarantees that $\{T_n^{N_\vk}\}$ lies in $\mathbb{B}_{N_\vk}$ and converges to $T_*^{N_\vk}$.
        For the third assumption, \cref{thm:error_ccdn} shows that with each fixed $n > 0$ the amplitude $T_n^{N_\vk}$ converges to $t_n$ in the sense of
        \[
                \lim_{N_\vk\rightarrow\infty} \|\mathcal{M}_\mathcal{K}t_n - T_n^{N_\vk}\|_1  = 0,
        \]
        suggesting that $\{T_n^{N_\vk}\} \subset \mathbb{B}_{N_\vk}$ converges to $t_n$
        with $\mathcal{K} \rightarrow \Omega^*$. Therefore $\mathbb{B}_{N_\vk}$ and $\mathcal{M}_\mathcal{K}t_n$ should be related to each other,
        which leads to the third assumption. 

        For the last assumption, we note that the two related error estimates in \cref{lem:error_amplitude_contraction_app} have the prefactor $C$ dependent on the amplitude $t_n$.
        The assumption here is stronger in the sense that $C$ needs to be independent of $t_n$. 
\end{rem}

\newtheorem*{restateddef1}{Rigorous Statement of \cref{thm:error_ccd_madelung}}
\begin{restateddef1}
        Under \cref{assump:contraction}, the finite-size error of the CCD correlation energy scales as $\Or(N_\vk^{-\frac13})$
        in each of the following scenarios 
        (1) the Madelung constant correction is only applied to the ERI contraction $\mathcal{A}_{N_\vk}$ 
        (2) the Madelung constant correction is only applied to the orbital energy $\varepsilon_{n\vk}^{N_\vk}$.
        When the Madelung constant correction is applied to both $\mathcal{A}_{N_\vk}$ and $\varepsilon_{n\vk}^{N_\vk}$ in the CCD calculation, 
        the overall finite-size error scales as $\Or(N_\vk^{-1})$.
\end{restateddef1}

\begin{proof}
The finite-size error in the CCD energy calculation with or without the Madelung constant correction can be estimated as 
\begin{align}
        \left|
                E_{\text{CCD}}^\text{TDL} - E_{\text{CCD}}^{N_\vk}
        \right|
        & = 
        \left|
                \mathcal{G}_\text{TDL}(t_*) - \mathcal{G}_{N_\vk}(T_*^{N_\vk})
        \right|
        \nonumber\\
        & \leqslant 
        \left|
                \mathcal{G}_{N_\vk}(\mathcal{M}_{\mathcal{K}}t_*) - \mathcal{G}_{N_\vk}(T_*^{N_\vk})
        \right|
        +
        \left|
                \mathcal{G}_\text{TDL}(t_*) - \mathcal{G}_{N_\vk}(\mathcal{M}_{\mathcal{K}}t_*)
        \right|
        \nonumber\\
        & \leqslant 
        C \left\|\mathcal{M}_{\mathcal{K}}t_* - T_*^{N_\vk}\right\|_1
        +
        C N_\vk^{-1}
        \label{eqn:error_ccd_converge},
\end{align}
where the last inequality uses the boundedness of $|W_{ijab}|$ in $\mathcal{G}_{N_\vk}$ and \cref{lem:error_energy}. 
To first estimate the above amplitude error when the Madelung constant correction is applied to both orbital energy and ERI contractions, 
we consider the error splitting \cref{eqn:decomposition_amplitude_app} for the amplitude calculation at the $n$-th fixed point iteration as
\begin{align*}
        \left\|
                \mathcal{M}_{\mathcal{K}}t_{n} - T_{n}^{N_\vk}
        \right\|_1
        & \leqslant 
        \left\|
        \dfrac{1}{\varepsilon^\text{TDL}}
        \left[
        \mathcal{M}_{\mathcal{K}}
        \mathcal{A}_\text{TDL}(t_{n-1})
        -
        \mathcal{A}_{N_\vk, \xi}(\mathcal{M}_{\mathcal{K}}t_{n-1})
        \right]
        \right\|_1
        \nonumber\\
        & \hspace*{1em}
        + 
        \left\|
        \left[
        \dfrac{1}{\varepsilon^\text{TDL}}
        -
        \dfrac{1}{\varepsilon^{N_\vk, \xi}}
        \right]
        \mathcal{A}_{N_\vk, \xi}(\mathcal{M}_{\mathcal{K}}t_{n-1})
        \right\|_1
        \nonumber\\
        & \hspace*{1em}
        +
        \left\|
        \dfrac{1}{\varepsilon^{N_\vk, \xi}}
        \left[
        \mathcal{A}_{N_\vk, \xi}( \mathcal{M}_{\mathcal{K}}t_{n-1})
        -
        \mathcal{A}_{N_\vk, \xi}(T_{n-1}^{N_\vk})
        \right]
        \right\|_1
        \\
        &
        \leqslant 
        CN_\vk^{-1}
        + 
        L  
        \left\|
                \mathcal{M}_\mathcal{K}t_{n-1} - T_{n-1}^{N_\vk}
        \right\|_1,
\end{align*}
where the last estimate uses the assumption in \cref{eqn:assumption_ccd_4} for the first term and the 
assumptions of contraction maps in \cref{eqn:assumption_ccd_2} and $\mathcal{M}_\mathcal{K}t_{n-1}\in \mathbb{B}_{N_\vk}$ in \cref{eqn:assumption_ccd_3} for the third term.
Since the initial guesses in the finite and the TDL cases satisfy $\|\mathcal{M}_\mathcal{K}t_0 - T_0^{N_\vk}\|_1= 0$, 
we can recursively derive that 
\[
        \left\|
                \mathcal{M}_{\mathcal{K}}t_{n} - T_{n}^{N_\vk}
        \right\|_1
        \leqslant 
        C \dfrac{1 - L^n}{1 - L} N_\vk^{-1},
\]
and thus the first assumption \cref{eqn:assumption_ccd_1} gives 
\[
        \left\|
                \mathcal{M}_{\mathcal{K}}t_{*} - T_{*}^{N_\vk}
        \right\|_1
        = \lim_{n\rightarrow\infty}
        \left\|
                \mathcal{M}_{\mathcal{K}}t_{n} - T_{n}^{N_\vk}
        \right\|_1
        \leqslant 
        C N_\vk^{-1}.
\]
Plugging this estimate into \cref{eqn:error_ccd_converge} then finishes the proof for the scenario when the Madelung 
constant correction is applied to both orbital energies and ERI contractions. 

For the two scenarios with partial Madelung constant correction, a similar analysis as above gives 
\begin{align*}
        \left\|
                \mathcal{M}_{\mathcal{K}}t_{n} - T_{n}^{N_\vk}
        \right\|_1
        \leqslant 
        CN_\vk^{-\frac13}
        + 
        L  
        \left\|
                \mathcal{M}_\mathcal{K}t_{n-1} - T_{n-1}^{N_\vk}
        \right\|_1,
\end{align*}
where the dominant $CN_\vk^{-\frac13}$ error comes from uncorrected orbital energy or uncorrected ERI contraction term.
Recursively, we can obtain the amplitude error estimate as 
\[
        \left\|
                \mathcal{M}_{\mathcal{K}}t_{*} - T_{*}^{N_\vk}
        \right\|_1
        \leqslant 
        C N_\vk^{-\frac13},
\] 
which finishes the proof for the two scenarios with partial correction. 
\end{proof}

\section{Quadrature error estimate for periodic functions with algebraic singularity}
\label{appendix:quaderror}
This section presents a collection of lemmas that provide quadrature error estimates for trapezoidal rules over periodic functions with specific algebraic singularities, which are used in this paper. 
Most of the lemmas are proven in \cite{XingLin2023} and are restated here for completeness. 
In addition, we introduce and prove a new quadrature error estimate that is critical in describing the efficacy of the Madelung constant correction and inverse symmetry in reducing the quadrature error in orbital energies and ERI contractions. 

The lemmas presented in this section provide the asymptotic scaling of the quadrature errors and also a quantitative relationship between the prefactors in the estimate and the algebraic singularities of the integrand. 
In addition to the singularity order, we further quantitatively characterize the algebraic singularity as follows.
For a univariate function $f(\vx)$ that is smooth everywhere in $V$ except at $\vx = \vx_0$ with algebraic singularity of order $\gamma$, 
we define a constant  
\begin{align}
        \mathcal{H}_{V, \vx_0}^{l}(f) 
& = \min
\left\{
        C: 
        \left|
                \partial^\valpha_\vx f(\vx)
        \right|
        \leqslant 
        C 
        |\vx - \vx_0|^{\gamma - |\valpha|}, 
        \forall |\valpha|\leqslant l, 
        \forall \vx \in V\setminus\{\vx_0\}
\right\}
\nonumber\\
& = 
        \max_{|\valpha|\leqslant l} 
        \left\| 
                \left(\partial^\valpha_\vx f(\vx)\right) / |\vx - \vx_0|^{\gamma - |\valpha|}
        \right\|_{L^{\infty}(V)}.
\nonumber 
\end{align}
For a multivariate function $f(\vx, \vy)$ that is smooth everywhere in $V_X\times V_Y$ except at $\vx = \vx_0$ with algebraic singularity of order $\gamma$,
we define a constant 
\begin{align}
\mathcal{H}_{V_X\times V_Y, (\vx_0,\cdot)}^{l}(f) 
& = \min
\left\{
        C: 
        \left|
                \partial^\valpha_\vx \partial_\vy^\vbeta f(\vx, \vy)
        \right|
        \leqslant 
        C 
        |\vx - \vx_0|^{\gamma - |\valpha|}, 
        \forall |\valpha|, |\vbeta|\leqslant l, 
        \forall \vx \in V_X\setminus\{\vx_0\}, \vy\in V_Y
\right\}
\nonumber\\
& = 
        \max_{|\valpha|\leqslant l} 
        \left\| 
                \left(\partial^\valpha_\vx \partial_\vy^\vbeta f(\vx, \vy)\right)
                / |\vx - \vx_0|^{\gamma - |\valpha|}
        \right\|_{L^{\infty}(V\times V)},
\nonumber
\end{align}
where ``$\cdot$'' in the subscript ``$(\vx_0, \cdot)$'' is a placeholder to indicate the smooth variable.
Using these two quantities, we have following function estimates that will be extensively used in this section 
\begin{align*}
        & 
        \left|
                \partial^\valpha_\vx f(\vx)
        \right|
        \leqslant  
        \mathcal{H}_{V, \vx_0}^{l}(f) |\vx - \vx_0|^{\gamma - |\valpha|}, 
        \qquad
        \forall l \geqslant |\valpha|,\ \forall \vx \in V\setminus\{\vx_0\},
        \\
        & 
        \left|
                \partial^\valpha_\vx \partial^\vbeta_\vy f(\vx, \vy)
        \right|
        \leqslant  
        \mathcal{H}_{V_X\times V_Y, (\vx_0, \cdot)}^{l}(f) |\vx - \vx_0|^{\gamma - |\valpha|}, 
        \qquad
        \forall l \geqslant |\valpha|, |\vbeta|,\ \forall \vx \in V_X\setminus\{\vx_0\},\vy\in V_Y.
\end{align*}

\subsection{Existing results from Ref.~\cite{XingLin2023}}

\begin{lem}\label{lem:quaderror0}
        Let $f(\vx)$ be smooth and periodic in $V = [-\frac12, \frac12]^d$. 
        The quadrature error of a trapezoidal rule using an $m^d$-sized uniform mesh $\mathcal{X}$ in $V$ decays super-algebraically as
        \[
                \left|\mathcal{E}_V(f, \mathcal{X})\right| \leqslant C_l m^{-l},\quad \forall l > 0.
        \]
\end{lem}
\begin{rem}\label{rem:quaderror0}
        If we replace $f(\vx)$ by $f(\vx, \vy)$ defined in $V\times V_Y$ which is smooth and periodic with respect to $\vx$ for each $\vy\in V_Y$ 
        and satisfies $\sup_{\vx\in V, \vy\in V_Y}|\partial_{\vx}^{\valpha}f(\vx, \vy)| < \infty$ for any $\valpha \geqslant 0$, \cref{lem:quaderror0} can be generalized
        as 
        \[
                \left|\mathcal{E}_V(f(\cdot, \vy), \mathcal{X})\right| \leqslant C_l m^{-l},\quad \forall l > 0, \forall \vy \in V_Y,
        \] 
        where constant $C_l$ is independent of $\vy \in V_Y$. 
\end{rem}

\begin{lem}\label{lem:quaderror1}
        Let $f(\vx)$ be periodic with respect to $V = [-\frac12,\frac12]^d$ and smooth everywhere except at $\vx = \bm{0}$ of order $\gamma \geqslant -d+1$.     
        At $\vx = \bm{0}$, $f(\vx)$ is set to $0$.
        The quadrature error of a trapezoidal rule using an $m^d$-sized uniform mesh $\mathcal{X}$ that contains $\vx = \bm{0}$ can be estimated as 
        \[
                \left|\mathcal{E}_V(f, \mathcal{X})\right| \leqslant C \mathcal{H}_{V, \bm{0}}^{d+\max(1,\gamma)}(f) m^{-(d + \gamma)}.
        \]
        If $f(\bm{0})$ is set to an $\Or(1)$ value in the calculation,
        it introduces additional $\Or(m^{-d})$ quadrature error.

\end{lem}
\begin{rem}\label{rem:quaderror1}
        If we replace $f(\vx)$ by $f(\vx, \vy)$ defined in $V \times V_Y$ which is smooth everywhere in $V\times V_Y$ except at $\vx = \bm{0}$ of 
        order $\gamma$, \cref{lem:quaderror1} can be generalized to 
        \[
                \left|\mathcal{E}_V(f(\cdot, \vy), \mathcal{X})\right| \leqslant C \mathcal{H}_{V\times V_Y, (\bm{0},\cdot)}^{d+\max(1,\gamma)}(f) m^{-(d + \gamma)}, \quad \forall \vy \in V_Y,
        \] 
        where the prefactor applies uniformly across all  $\vy\in V_Y$. 
\end{rem}

\begin{lem}\label{lem:quaderror2}
        Let $f(\vx) = f_1(\vx)f_2(\vx)$ where $f_1(\vx)$ and $f_2(\vx)$ are periodic with respect to $V = [-\frac12,\frac12]^d$ and
        \begin{itemize}
                \item $f_1(\vx)$ is smooth everywhere except at $\vx = \vz_1 = \bm{0}$ of order $\gamma\leqslant 0$,
                \item $f_2(\vx)$ is smooth everywhere except at $\vx = \vz_2 \neq \bm{0}$ of order $0$.
        \end{itemize}
        Consider an $m^d$-sized uniform mesh $\mathcal{X}$ in $V$. Assume that $\mathcal{X}$ satisfies that $\vz_1, \vz_2$ are either on the 
        mesh or $\Theta(m^{-1})$ away from any mesh points, and $m$ is sufficiently large that $|\vz_1 - \vz_2| = \Omega(m^{-1})$.
        At $\vx = \vz_1$ and $\vx = \vz_2$, $f(\vx)$ is set to $0$. 
        The trapezoidal rule using $\mathcal{X}$ has quadrature error
        \[
                \left|\mathcal{E}_V(f, \mathcal{X})\right| \leqslant C \mathcal{H}_{V, \vz_1}^{d+1}(f_1)\mathcal{H}_{V, \vz_2}^{d+1}(f_2)m^{-(d + \gamma)}.
        \]
        If $f(\vz_1)$ and $f(\vz_2)$ are set to arbitrary $\Or(1)$ values, it introduces additional $\Or(m^{-d})$ quadrature error. 
\end{lem}

\begin{rem}\label{rem:quaderror2}
        If we replace $f_i(\vx)$ with $i=1,2$ by $f_i(\vx, \vy)$ defined in $V\times V_Y$ which is smooth everywhere in $V\times V_Y$ except at $\vx = \vz_i$ of order $\gamma$ and $0$, respectively,
        \cref{lem:quaderror2} can be generalized to 
        \[
                \left|\mathcal{E}_V(f_1(\cdot, \vy)f_2(\cdot, \vy), \mathcal{X})\right| \leqslant C \mathcal{H}_{V\times V_Y, (\bm{0},\cdot)}^{d+1}(f_1) \mathcal{H}_{V\times V_Y, (\bm{0},\cdot)}^{d+1}(f_2)m^{-(d + \gamma)}, \quad \forall \vy \in V_Y,
        \] 
        where the prefactor applies uniformly across all $\vy \in V_Y$. 
\end{rem}

\subsection{A new quadrature error estimate}
\label{app:quaderror2_2}
Here we prove \cref{lem:quaderror2_2}, which is used in the quadrature error estimate of the ERI contractions with the Madelung constant correction in \cref{proof:error_amplitude_contraction}. 
This lemma is a generalization of the existing result in \cref{lem:quaderror2} with an additional condition.

\newtheorem*{restatedlem}{Restatement of \cref{lem:quaderror2_2}}
\begin{restatedlem}
        Let $f(\vx) = f_1(\vx)f_2(\vx)$ where $f_1(\vx)$ and $f_2(\vx)$ are periodic with respect to $V = [-\frac12,\frac12]^d$ and
        \begin{itemize}
                \item $f_1(\vx)$ is smooth everywhere except at $\vx = \vz_1 = \bm{0}$ of order $\gamma \leqslant -1$,
                \item $f_2(\vx)$ is smooth everywhere except at $\vx = \vz_2 \neq \bm{0}$ of order $0$,
                \item $\partial^{\valpha}_\vx f_2(\bm{0}) = 0$ for any derivative order $|\valpha| \leqslant s$. 
        \end{itemize}
        Assume $\gamma > -d$ for $f(\vx)$ to be integrable in $V$ and $\gamma+s+1 \leqslant 0$ so the leading algebraic singularity of $f(\vx)$ is at $\vx = \bm{0}$. 
        Consider an $m^d$-sized uniform mesh $\mathcal{X}$ in $V$. Assume that $\mathcal{X}$ satisfies that $\vz_1, \vz_2$ are either on the 
        mesh or $\Theta(m^{-1})$ away from any mesh points, and $m$ is sufficiently large that $|\vz_1 - \vz_2| = \Omega(m^{-1})$.
        At $\vx = \vz_1$ and $\vx = \vz_2$, $f(\vx)$ is set to $0$. 
        The trapezoidal rule using $\mathcal{X}$ has quadrature error
        \[
                \left|\mathcal{E}_V(f, \mathcal{X})\right| \leqslant C \mathcal{H}_{V, \vz_1}^{d+1}(f_1)\mathcal{H}_{V, \vz_2}^{d+1}(f_2) |\vz_1 - \vz_2|^{-(s+1)} m^{-(d + \gamma + s + 1)}.
        \]
\end{restatedlem}

\begin{proof}
        For $\vz_1 = \bm{0}$ and any $\vz_2 \in V$, we can introduce a proper translation $f(\vx) \rightarrow f(\vx - \vx_0)$ to move 
        both the singular points $\vz_1$ and $\vz_2$ to the smaller cube $[-\frac14, \frac14]^d$ in $V$.
        The target quadrature error can be correspondingly reformulated as 
        \[
                \mathcal{E}_V(f(\cdot), \mathcal{X}) =   
                \mathcal{E}_{V + \vx_0}(f(\cdot - \vx_0), \mathcal{X} - \vx_0) =  
                \mathcal{E}_V(f(\cdot - \vx_0), \mathcal{X} - \vx_0).
        \]
        which is the quadrature error of the translated function $f(\vx - \vx_0)$. 
        Without loss of generality, we assume such a translation has been applied to $f(\vx)$ and $\mathcal{X}$ and both singular points $\vz_1$ and $\vz_2$ lie in $[-\frac14, \frac14]^d$.

        Define a cutoff function $\psi \in \mathbb{C}_{c}^{\infty}(\mathbb{R}^n)$ satisfying
        \[
                \psi(\mathbf{x})= \begin{cases}1, & |\mathbf{x}| < \frac{1}{4} \\ 0, & |\mathbf{x}| > \frac12\end{cases}.
        \]
        Denote the distance between the two singular points as $\delta\vz = |\vz_2 - \vz_1|$ and define 
        two local cutoff functions that isolate the two singular points as 
        \begin{align*}
                \psi_{\delta\vz,1}(\vx) = \psi\left(\dfrac{\vx- \vz_1}{\delta\vz}\right),
                \qquad 
                \psi_{\delta\vz,2}(\vx) = \psi\left(\dfrac{\vx- \vz_2}{\delta\vz}\right),
        \end{align*}
        whose supports are both inside $V$. 
        The target integrand can be split into three parts as 
        \begin{equation*}
                f(\vx) = f(\vx)\psi_{\delta\vz,1}(\vx) +  f(\vx)\psi_{\delta\vz,2}(\vx) +  f(\vx)\left(1 - \psi_{\delta\vz,1}(\vx) - \psi_{\delta\vz,1}(\vx)\right).
        \end{equation*}
        All three terms satisfy periodic boundary condition on $\partial V$. 
        The first term is smooth everywhere except at $\vx = \vz_1$ of order $\gamma + s + 1$, 
        the second term is smooth everywhere except at $\vx = \vz_2$ of order $0$, and the last term is smooth everywhere. 
        Application of \cref{lem:quaderror0} and \cref{lem:quaderror1} to these terms suggests that 
        \begin{align}
                & 
                \left|
                \mathcal{E}_V\left(
                        f\psi_{\delta\vz, 1}, \mathcal{X}
                \right)
                \right|
                \leqslant C \mathcal{H}_{V,\vz_1}^{d+1}\left(f\psi_{\delta\vz, 1}\right) m^{-(d+\gamma+s+1)}, 
                \label{eqn:lem_error1}\\
                & 
                \left|
                \mathcal{E}_V\left(
                        f\psi_{\delta\vz, 2}, \mathcal{X}
                \right)
                \right|
                \leqslant C \mathcal{H}_{V,\vz_2}^{d+1}\left(f\psi_{\delta\vz, 2}\right) m^{-d},  
                \label{eqn:lem_error2}\\
                & 
                \left|
                \mathcal{E}_V\left(
                        f(1 - \psi_{\delta\vz, 1} - \psi_{\delta\vz, 2}), \mathcal{X}
                \right)
                \right|
                \leqslant C_{|\valpha|} \left(\int_{V}\ud\vx \left| \partial^{\valpha}_\vx f(1 - \psi_{\delta\vz, 1} - \psi_{\delta\vz, 2})\right|\right) m^{-|\valpha|},  \quad \forall |\valpha| > d.
                \label{eqn:lem_error3}
        \end{align}
        Note that the last estimate means a super-algebraic decaying error and constant $C_{|\valpha|}$ only depends on $V$ and $|\valpha|$ (see Eqn.~(H.6) in 
        Ref.~\cite{XingLin2023} for the derivation of this detailed prefactor for \cref{lem:quaderror0}).
        Since $\gamma+s+1\leqslant 0$, the three estimates together prove that the quadrature error of $f(\vx)$ scales asymptotically as $m^{-(d+\gamma+s+1)}$ 
        as $m \rightarrow \infty$. 
        In order to describe the extreme case where the two singular points are only $\Or(m^{-1})$ away from each other for a given mesh $\mathcal{X}$, i.e., $\delta\vz = \Or(m^{-1})$, 
        we provide a precise description of the three prefactors in the estimates above using $\delta\vz$.

        For the prefactor in \cref{eqn:lem_error1}, it is defined by 
        \begin{equation*}
                \mathcal{H}_{V,\vz_1}^{l}\left(f\psi_{\delta\vz, 1}\right)
                =
                \max_{|\valpha|\leqslant l} 
                \left\| 
                        \left(\partial^\valpha_\vx f(\vx)\psi_{\delta\vz, 1}(\vx)\right) / |\vx - \vz_1|^{\gamma + 1 - |\valpha|}
                \right\|_{L^{\infty}(V)}, \quad \forall l > 0.
        \end{equation*}
        For any derivative order $|\valpha| \leqslant l$, we have 
        \begin{align}
                \left|
                        \partial^\valpha_\vx \left(f_1(\vx)f_2(\vx)\psi_{\delta\vz, 1}(\vx)\right)
                \right|
                & \leqslant 
                C \sum_{\valpha_1 + \valpha_2 + \valpha_3 = \valpha}
                \left|
                        \partial^{\valpha_1}_\vx f_1(\vx) 
                \right|
                \left|
                        \partial^{\valpha_2}_\vx f_2(\vx) 
                \right|
                \left|
                        \partial^{\valpha_3}_\vx \psi_{\delta\vz, 1}(\vx)
                \right|
                \nonumber\\
                & \leqslant
                C \sum_{\valpha_1 + \valpha_2 + \valpha_3 = \valpha}
                \mathcal{H}_{V, \vz_1}^{l}(f_1) |\vx - \vz_1|^{\gamma - |\valpha_1|}
                \left|
                        \partial^{\valpha_2}_\vx f_2(\vx) 
                \right|
                \delta\vz^{-|\valpha_3|}
                \delta_{|\vx - \vz_1| < \frac12\delta\vz}
                \nonumber \\
                & \leqslant 
                C \sum_{\valpha_1 + \valpha_2 + \valpha_3 = \valpha}
                \mathcal{H}_{V, \vz_1}^{l}(f_1) |\vx - \vz_1|^{\gamma - |\valpha_1| - |\valpha_3|}
                \left|
                        \partial^{\valpha_2}_\vx f_2(\vx) 
                \right|
                \delta_{|\vx - \vz_1| < \frac12\delta\vz}.
                \label{eqn:lem_estimate1}
        \end{align} 
        The last inequality uses $|\vx - \vz_1| \leqslant \frac12\delta\vz$ by noting that $\partial^{\valpha_3}_\vx \psi_{\delta\vz, 1}(\vx)$ is zero when $|\vx - \vz_1| > \frac12\delta\vz$. 
        Next we estimate $\left|\partial^{\valpha_2}_\vx f_2(\vx)\right|$ in the ball domain $|\vx - \vz_1| \leqslant \frac12\delta\vz$. 
        When $|\valpha_2| \leqslant s$, we have 
        \begin{align*}
                \left|
                        \partial^{\alpha_2}_\vx f_2(\vx) 
                \right|
                & = 
                \left|
                        \sum_{|\vbeta| = s+1 - |\valpha_2|} \frac{s+1-|\valpha_2|}{\vbeta!} (\vx - \vz_1)^{\vbeta} \int_0^1 (1-t)^{s-|\valpha_2|} \partial^{\valpha_2 + \vbeta}_\vx f_2(\vz_1 + t(\vx - \vz_1))\ud t
                \right|
                \\
                & \leqslant 
                C |\vx - \vz_1|^{s+1 - |\valpha_2|} \sup_{|\vx - \vz_1| \leqslant \frac12\delta\vz, |\vbeta| = s+1} \left|\partial^{\vbeta}_\vx f_2(\vx)\right|
                \\
                & \leqslant 
                C |\vx - \vz_1|^{s+1 - |\valpha_2|} \mathcal{H}_{V, \vz_2}^{s+1}(f_2) \delta\vz^{-(s+1)}.
        \end{align*}
        The first equality applies the Taylor expansion of $\partial^{\valpha_2}_\vx f_2(\vx)$ at $\vx = \vz_1$ with the assumption that $\partial^{\vbeta}_\vx f_2(\vz_1) = \bm{0}$ for any $|\vbeta| \leqslant s$. 
        The last inequality uses the algebraic singularity of $f_2(\vx)$ at $\vx = \vz_2$. 
        When $|\valpha_2| > s$, we have 
        \begin{align*}
                \left|
                        \partial^{\alpha_2}_\vx f_2(\vx) 
                \right|
                \leqslant 
                \sup_{|\vx - \vz_1| < \frac12\delta\vz}
                \left|
                        \partial^{\alpha_2}_\vx f_2(\vx) 
                \right|
                \leqslant  
                C \mathcal{H}_{V, \vz_2}^{l}(f_2) \delta\vz^{-|\valpha_2|}
                \leqslant 
                C \mathcal{H}_{V, \vz_2}^{l}(f_2) |\vx - \vz_1|^{s+1-|\valpha_2|}\delta\vz^{-(s+1)}.
        \end{align*}
        Plugging the estimates of $|\partial^{\valpha_2}_\vx f_2(\vx)|$ above into \cref{eqn:lem_estimate1}, we obtain 
        \[
                \left|
                        \partial^\valpha_\vx \left(f_1(\vx)f_2(\vx)\psi_{\delta\vz, 1}(\vx)\right)
                \right| 
                \leqslant 
                C 
                \mathcal{H}_{V, \vz_1}^{l}(f_1) 
                \mathcal{H}_{V, \vz_2}^{l}(f_2)
                \delta\vz^{-(s+1)}
                |\vx - \vz_1|^{\gamma + s + 1 - |\valpha|}, \quad l \geqslant d+1.
        \]
        which suggests that 
        \[
                \mathcal{H}_{V,\vz_1}^{l}\left(f\psi_{\delta\vz, 1}\right)         
                \leqslant 
                C
                \mathcal{H}_{V, \vz_1}^{l}(f_1)
                \mathcal{H}_{V, \vz_2}^{l}(f_2) 
                \delta\vz^{-(s+1)}, \quad \forall l \geqslant d+1.
        \]
        Similar analysis can also be applied to the prefactor in \cref{eqn:lem_error2} to obtain 
        \[
                \mathcal{H}_{V,\vz_2}^{l}\left(f\psi_{\delta\vz, 2}\right)         
                \leqslant 
                C
                \mathcal{H}_{V, \vz_1}^{l}(f_1)
                \mathcal{H}_{V, \vz_2}^{l}(f_2) 
                \delta\vz^{-(s+1)}, \quad \forall l \geqslant d+1.
        \]

        For the estimate of the third term in \cref{eqn:lem_error3} with any $|\valpha| = l \geqslant d+1$, the prefactor can be bounded as
        \begin{align*}
                & 
                \int_{V}\ud\vx \left| \partial^{\valpha}_\vx f(\vx)(1 - \psi_{\delta\vz, 1}(\vx) - \psi_{\delta\vz, 2}(\vx))\right|
                \\
                \leqslant & 
                C\sum_{\valpha_1 + \valpha_2 + \valpha_3 = \valpha}
                \int_{V}\ud\vx 
                \left| 
                        \partial^{\valpha_1}_\vx f_1(\vx)
                \right|
                \left| 
                        \partial^{\valpha_2}_\vx f_2(\vx)
                \right|
                \left| 
                        \partial^{\valpha_3}_\vx 
                        (1 - \psi_{\delta\vz, 1}(\vx) - \psi_{\delta\vz, 2}(\vx))
                \right|
                \\
                \leqslant & 
                C
                \mathcal{H}_{V, \vz_1}^{l}(f_1)
                \mathcal{H}_{V, \vz_2}^{l}(f_2) 
                \sum_{\valpha_1 + \valpha_2 + \valpha_3 = \valpha}
                \int_{V}\ud\vx 
                \left|\vx - \vz_1\right|^{\gamma - |\valpha_1|}
                \left|\vx - \vz_2\right|^{-|\valpha_2|}
                \left| 
                        \partial^{\valpha_3}_\vx 
                        (1 - \psi_{\delta\vz, 1}(\vx) - \psi_{\delta\vz, 2}(\vx))
                \right|.
        \end{align*}
        For the integral with each set of $(\valpha_1, \valpha_2, \valpha_3)$, we consider two cases:
        \begin{itemize}
                \item $\valpha_3 = \bm{0}$. The integral can be bounded by  
                \begin{align*}
                        (*)\leqslant 
                        \int_{V \setminus \left(B(\vz_1, \frac14\delta\vz) \cup B(\vz_2, \frac14\delta\vz)\right)}\ud\vx 
                        \left|\vx - \vz_1\right|^{\gamma - |\valpha_1|}
                        \left|\vx - \vz_2\right|^{|\valpha_1| - |\valpha|}
                        \leqslant C (1 + \delta\vz^{\gamma + d - |\valpha|}),
                \end{align*}
                using the H{\"o}lder-inequality technique developed in the proof of \cref{lem:quaderror2} in \cite{XingLin2023}. 

                \item $\valpha_3 > \bm{0}$. The integral can be bounded by 
                \begin{align*}
                        (*)
                        & 
                        \leqslant 
                        C 
                        \int_{ \left(B(\vz_1, \frac12\delta\vz)\setminus B(\vz_1, \frac14\delta\vz)\right) \cup \left(B(\vz_2, \frac12\delta\vz)\setminus B(\vz_2, \frac14\delta\vz)\right)}\ud\vx 
                        \left|\vx - \vz_1\right|^{\gamma - |\valpha_1|}
                        \left|\vx - \vz_2\right|^{-|\valpha_2|}
                        \delta\vz^{-|\valpha_3|}
                        \\
                        &
                        \leqslant C \delta\vz^{\gamma + d - |\valpha|}.
                \end{align*}
        \end{itemize}
        Collecting these two terms together, the error estimate in \cref{eqn:lem_error3} with any $|\valpha| = l \geqslant d + 1$ can be further bounded  as 
        \begin{align*}
                \left|
                        \mathcal{E}_V\left(
                                f(1 - \psi_{\delta\vz, 1} - \psi_{\delta\vz, 2}), \mathcal{X}
                        \right)
                \right|
                & \leqslant C_{l}
                \mathcal{H}_{V, \vz_1}^{l}(f_1)
                \mathcal{H}_{V, \vz_2}^{l}(f_2)  \delta\vz^{\gamma + d - l}m^{-l}
                \\
                & \leqslant  C_{l}
                \mathcal{H}_{V, \vz_1}^{l}(f_1)
                \mathcal{H}_{V, \vz_2}^{l}(f_2)  \delta\vz^{-(s+1)}m^{-(\gamma + s + 1 + d- l)}m^{-l}
                \\
                & \leqslant  C_{l}
                \mathcal{H}_{V, \vz_1}^{l}(f_1)
                \mathcal{H}_{V, \vz_2}^{l}(f_2)  \delta\vz^{-(s+1)}m^{-(\gamma + s + 1  + d)},
        \end{align*}
        where the second inequality uses the assumption $\delta\vz = \Omega(m^{-1})$ and $\gamma+s+1+d-l \leqslant 0$. 
        Gathering the above prefactor descriptions for \cref{eqn:lem_error1,eqn:lem_error2,eqn:lem_error3}, we obtain 
        \[
        \left|
                \mathcal{E}_V\left(
                        f, \mathcal{X}
                \right)
        \right| 
        \leqslant 
        C
        \mathcal{H}_{V, \vz_1}^{d+1}(f_1)
        \mathcal{H}_{V, \vz_2}^{d+1}(f_2)  \delta\vz^{-(s+1)}m^{-(\gamma + s + 1 + d)},
        \]
        which finishes the proof.
       \end{proof}

\end{document}